\documentclass[11pt,reqno]{amsart}
\usepackage{amssymb,amsfonts,amsthm,amsmath,enumitem}
\usepackage{color}
\usepackage[pagebackref=false]{hyperref}
\usepackage[left=1 in,top=1 in,right=1 in,bottom=1 in]{geometry}
\usepackage{cite}

\newcommand{\beq}{\begin{equation}}
\newcommand{\eeq}{\end{equation}}
\newcommand{\beqs}{\begin{equation*}}
\newcommand{\eeqs}{\end{equation*}}
\newcommand{\ba}{\begin{array}}
\newcommand{\ea}{\end{array}}
\newcommand{\beas}{\begin{eqnarray*}}
\newcommand{\eeas}{\end{eqnarray*}}
\newcommand{\bea}{\begin{eqnarray}}
\newcommand{\eea}{\end{eqnarray}}
\newcommand{\bal}{\begin{align}}
\newcommand{\eal}{\end{align}}

\newcommand{\bals}{\begin{align*}}
\newcommand{\eals}{\end{align*}}

\newcommand{\tnum}{\rm(\roman*)}
\newcommand{\rnum}{\rm(\alph*)}

\newcommand{\R}{\ensuremath{\mathbb R}}

\newcommand{\norm}[1]{\| {#1} \|}

\newcommand{\bds}{\begin{displaystyle}}
\newcommand{\eds}{\end{displaystyle}}

\def\longequals{\mathbin{=\kern-2pt=}}
\def\eqdef{\stackrel{\rm def}{=}}

\def\varep{\varepsilon}

\def\ddt{\frac{\rm d}{\rm dt}}
\def\d{{\rm d}}

\newcommand{\remove}[1]{} 
\renewcommand{\remove}[1]{#1} 

\newtheorem{theorem}{Theorem}[section]

\newtheorem{lemma}[theorem]{Lemma}

\newtheorem{proposition}[theorem]{Proposition}

\newtheorem{assumption}[theorem]{Assumption}
\theoremstyle{remark}
\newtheorem{remark}[theorem]{\bf{Remark}}

\newcommand{\essup}{\mathop{\mathrm{ess\,sup}}}

\definecolor{darkred}{rgb}{.70,.12,.20}

\definecolor{darkgreen}{rgb}{.20,.52,.14}

\newcommand{\eR}{\vec k}

\numberwithin{equation}{section}

\title{On compressible fluid flows of {F}orchheimer-type in rotating heterogeneous porous media}

\date{\today}

\subjclass[2020]{35Q35, 76S05, 35B45, 35K20, 35K65, 35K67}

\keywords{rotating fluids, heterogeneous porous media, Forchheimer flows, compressible fluids, singular/degenerate PDE, weighted Sobolev inequalities, weighted trace theorem, a priori estimates}

\author{Emine Celik$^{1}$}
\address{$^{1}$Department of Mathematics, Sakarya University\\
54050, Sakarya, T\"urkiye}
\email{eminecelik@sakarya.edu.tr}

\author{Luan Hoang$^{2,*}$}
\address{$^2$Department of Mathematics and Statistics,
Texas Tech University\\
1108 Memorial Circle, Lubbock, TX 79409--1042, U. S. A.}
\email{luan.hoang@ttu.edu}

\author{Thinh Kieu$^{3}$}
\address{$^{3}$Department of Mathematics, University of North Georgia, Gainesville Campus\\
3820 Mundy Mill Rd., Oakwood, GA 30566, U. S. A.}
\email{thinh.kieu@ung.edu}
\thanks{$^*$Corresponding author.}

\begin{document}

\begin{abstract}
We study the dynamics of compressible fluids in rotating heterogeneous porous media. The  fluid flow is of {F}orchheimer-type and is subject to a mixed mass and volumetric flux boundary condition. The governing equations are reduced to a nonlinear partial differential equation for the pseudo-pressure. This parabolic-typed equation can be degenerate and/or singular in the spatial variables, the unknown and its gradient. We establish the $L^\alpha$-estimate for the solutions, for any positive number $\alpha$, in terms of the initial and boundary data and the angular speed of rotation. It requires new elliptic and parabolic Sobolev inequalities and trace theorem with multiple weights that are suitable to the nonlinear structure of the equation. The $L^\infty$-estimate is then obtained without imposing any conditions on the $L^\infty$-norms of the weights and the initial and boundary data. 
\end{abstract}

\maketitle

\tableofcontents

\pagestyle{myheadings}\markboth{E. Celik, L. Hoang,  and T. Kieu}
{\sc Compressible fluid flows of {F}orchheimer-type in rotating heterogeneous porous media}

\section{Introduction}\label{intro}

In this paper, we study the following issues simultaneously.
\begin{enumerate}[label=\rnum]
    \item\label{Is1} Generalized Forchheimer flows, i.e., non-Darcy ones, in porous media.
    \item\label{Is2} Rotating fluids.
    \item\label{Is3} Heterogeneous porous media.
    \item\label{Is4} Compressible fluids with the gravity effect.
\end{enumerate}

Previous work may deal with some combinations of the above issues, but not all together. For example, Forchheimer flows of compressible fluids in homogeneous porous media are studied in \cite{ABHI1,HI1,HI2,HIKS1,HKP1,HK1,HK2,CHK1,CHK2,CHIK1,HIK1,HIK2} without rotation, and in \cite{CHK3,CHK4} with rotation. In particular, the case of anisotropic Forchheimer flows is treated in \cite{HK3} using the model from \cite{BarakBear81}. In heterogeneous porous media, they are investigated in \cite{CH1,CH2} for slightly compressible fluids, and in \cite{CHK5} for other gases, all without rotation. 

Regarding \ref{Is1}, we recall that most of the research on fluid flows in porous media often uses the Darcy law. However, there is a deviation from this ubiquitous law when the Reynolds number is large. In that case, the Forchheimer equations \cite{Forchh1901,ForchheimerBook} are commonly utilized instead. 
The reader is referred to \cite{MuskatBook,BearBook,StraughanBook,NieldBook,Scheidegger1974} and also \cite{ABHI1,HI2} for more information about the Forchheimer flows and their generalizations. Compared to the Darcy flows, mathematical analysis of the Forchheimer flows has received much less attention. One of the reasons is the nonlinear nature of the Forchheimer equations in contrast to the linear Darcy's one. Furthermore, the Forchheimer flows in specific settings such as in \ref{Is2}--\ref{Is4} bring even more difficulties into the modeling and analysis.
In addition to the references cited in the second paragraph above, there are treatments of broader non-Darcy flows such as the two-phase Forchheimer ones in \cite{HIK1,HIK2}, and a mixture of the Darcy and non-Darcy regimes in \cite{CHIK1}. This current work continues our research program on the Forchheimer-typed fluid flows in porous media. By studying all four issues \ref{Is1}--\ref{Is4}, we look deeper into the models and properties of these flows.

Regarding \ref{Is2}, fluids in rotating porous media occur in many physical systems in engineering and geophysics.  Some instances are geophysical flows with the effect of the earth's rotation, rotating filtration devices in food process and chemical industries. The rotation introduces the Coriolis and centrifugal forces which influence the dynamics of fluids in porous structures, see \cite[Section 6.22]{NieldBook}, also \cite{capone2014stability,vadasz1998coriolis,VadaszRPM2019},  and  references therein. However, there are not many mathematical studies of compressible Forchheimer flows in rotating porous media. Regarding this matter, our previous work \cite{CHK3,CHK4} only dealt with slightly compressible fluids and homogeneous porous media. By adding other types of compressible fluids and the inhomogeneity of the media to consideration, we greatly increase the complexity of the problem. Nonetheless, we are able to construct a rigorous mathematical framework with the use of nonlinear partial differential equations (PDEs) and derive fundamental estimates for their solutions.

Regarding \ref{Is3}, the properties of porous media can vary from place to place. Some examples are soil's changing types going downward (soil profile), oil and gas reservoirs with fractures, land in places with or without vegetation, different packing of pebbles near or far away from the boundary, media built with layers of different materials. From mathematical point of view, they introduce another level of difficulty when the coefficients depend on the spatial variables and can be singular or degenerate. Their analysis requires sophisticated tools from real and/or harmonic analysis. Some of those techniques will be further developed in this paper.

Regarding \ref{Is4}, most research on compressible fluids in porous media ignores the effect of gravity. This omission gives many models the advantage of having some important scaling properties. One of the most famous consequences is the Barenblatt profiles which play a crucial role in understanding porous medium equations, see e.g. \cite{VazquezPorousBook}.
In contrast, gravity is fully accounted for in this paper.  This inclusion, together with the Forchheimer equations, prompts the Barenblatt profiles to cease to exist. Moreover, it induces another nonlinear term (of the pseudo-pressure) in our model which can potentially has too large growth. Thus, it may hinder the establishment of estimates for solutions in terms of the initial and boundary data. It is worth mentioning that the model and analysis of incompressible fluids, see e.g. \cite{Zabensky2015a,ChadamQin,Payne1999a,Payne1999b,MTT2016,CKU2006,KR2017,Straughan2010} and references therein, are totally different from ours. 

The current work inherits recent progress in the study of Forchheimer flows of compressible fluids in porous media \cite{CHK3,CHK4,CHK5,CH1,CH2}. Combining the methods in \cite{CHK4} and \cite{CHK5}, we are able to derive a new multi-faceted degenerate/singular nonlinear parabolic equation for the pseudo-pressure. We develop analysis tools which are suitable to the multi-factor dependent nonlinear structure of the new model. They are then utilized to obtain estimates of the solutions in terms of the angular speed of rotation, and the initial and boundary data, see the main Theorems \ref{Labound} and \ref{thm45} below.

The paper is organized as follows.
In section \ref{IBVPsec}, we derive a  nonlinear PDE for the pseudo-pressure $u(x,t)$ which can be degenerate/singular in $x$, $u$ and $\nabla u$. It is subject to a mixed boundary condition. The main properties of the nonlinear term are in Lemma \ref{Xlem}. They depend   not only on the multi-factors $x$, $u$ and $\nabla u$, but also explicitly on the angular speed of rotation. Based on these, a more general initial boundary value problem is formulated in \eqref{ibvpg} with Hypotheses \ref{H1}--\ref{H7}. 
In section \ref{Prem},  we establish weighted Sobolev inequalities (elliptic and parabolic versions) in Lemmas \ref{RoS1} and \ref{RoS2}, and a weighted trace theorem in Lemma \ref{Rtrace}. These extend our previous results in \cite{CHK5} to suit the nonlinearity and its $x$-dependent degeneracy/singularity in  problem \eqref{ibvpg}. 
%
In section \ref{Laest}, by utilizing the new inequalities in Section \ref{Prem}, we first derive a differential inequality in Proposition \ref{Diff4u}  for the spatial $L^\alpha$-norm of a solution $u(x,t)$ for a finite positive number $\alpha$. It is then used to obtain the $L^\alpha$-estimate for $u$ in terms of the initial and boundary data as well as the angular speed of rotation  in  Theorem \ref{Labound}.
In section \ref{maxsec}, we derive estimates for the $L^\infty$-norm in $(x,t)$ of the solution $u$. We use Moser's iteration for suitable weighted norms. It is based on the key Lemmas \ref{caccio} and \ref{GLk}. The first $L^\infty-L^\alpha$ norm relation for the solution is proved in Theorem \ref{LinfU}. In Theorem \ref{thm45}, the final spatiotemporal $L^\infty$-estimates, with $t>0$, for the solution  are then established in terms of the angular speed of rotation and the $L^\alpha$-norms, with different finite values of $\alpha$, of the initial and boundary data.
It is worth mentioning that there are no conditions imposed on the $L^\infty$-norms of the weight functions or their reciprocals as well as on the initial and boundary data. 

\section{The initial boundary value problem}\label{IBVPsec}

The equation for the Darcy flows in rotating porous media written in a rotating frame, see \cite{VadaszBook},  is
\beq\label{RDarcy}
av+ c\rho \widetilde\Omega \eR\times v + \rho \widetilde\Omega^2 \eR\times (\eR\times x)=-\nabla p+\rho \vec g,
\eeq
where 
$v\in\R^3$ is the fluid velocity, 
$\rho\in\R_+=[0,\infty)$ is the fluid density, 
$p\in\R$ is the fluid pressure, 
$\vec g\in\R^3$ is the gravitational acceleration,
$\widetilde\Omega \eR$ is the constant angular velocity with
$\widetilde\Omega\in \R_+$ being the constant angular speed of rotation and
$\eR\in\R^3$ being a constant unit vector, 
$x\in\R^3$ is the position in the rotating frame,
$\widetilde\Omega \eR\times v$ is the Coriolis acceleration,
and  $\widetilde\Omega^2 \eR\times (\eR\times x)$ is the centripetal acceleration.

Equation \eqref{RDarcy} is extended to the flows obeying the two-term Forchheimer law \cite{VadaszBook}  by
\beq\label{Fv}
 av+b|v|v + c\rho \widetilde\Omega \eR\times v +\rho \widetilde\Omega^2 \eR\times (\eR\times x)=-\nabla p+\rho \vec g.
\eeq
In this paper, we consider generalized Forchheimer  equations in rotating porous media
\beq\label{GF}
 \sum_{i=0}^N a_i  |v|^{\bar\alpha_i} v+ c\rho \widetilde\Omega \eR\times v +\rho \widetilde\Omega^2 \eR\times (\eR\times x)=-\nabla p+\rho \vec g.
\eeq
In equations \eqref{RDarcy}, \eqref{Fv} and \eqref{GF}, $a$, $b$, $c$ and $a_i$ are positive numbers,  see e.g. \eqref{Ftwo} below,
\beq \label{Nal}
\text{$N\ge 1$ is an integer, and constants $\bar\alpha_0=0<\bar\alpha_1<\bar\alpha_2<\ldots<\bar\alpha_N$ are real numbers.}
\eeq 

Equation \eqref{GF} covers the case \eqref{Fv} as well as other Forchheimer's three-term and power laws \cite{MuskatBook, BearBook,Scheidegger1974}. 
By using dimension analysis of Muskat \cite{MuskatBook} and then Ward \cite{Ward64}, we refine equation \eqref{GF} as follows
 \beq\label{FM}
 \sum_{i=0}^N a_i \rho^{\bar\alpha_i} |v|^{\bar\alpha_i} v+ \frac{2\rho \widetilde\Omega}{\widetilde\phi} \eR\times v+\rho \widetilde\Omega^2 \eR\times (\eR\times x)=-\nabla p +\rho \vec g,
 \eeq
where $\widetilde\phi\in(0,1)$ is the porosity.
In fact, when $N=1$, we have a more detailed version of \eqref{Fv}:
 \beq\label{Ftwo}
\frac{\mu}{k} v+ \frac{c_F \rho}{\sqrt{k}}|v|v + \frac{2\rho \widetilde\Omega}{\widetilde\phi} \eR\times v+\rho \widetilde\Omega^2 \eR\times (\eR\times x)=-\nabla p +\rho \vec g,
 \eeq
where $\mu$ is the dynamic viscosity, $k$ is the permeability, and $c_F$ is the Forchheimer constant \cite{Ward64}.
Even in this case, there is no analysis for these flows of compressible fluids.

Consider a heterogeneous porous medium in a non-empty, open, bounded subset $U$ of $\R^3$. Then the porosity $\widetilde\phi$ is a function $\widetilde\phi(x)$, for $x\in\bar U$,  and all $a_i$ in \eqref{FM} now depend on $x$, i.e., $a_i=a_i(x)$. Thus, equation \eqref{FM} becomes
\beq\label{nGF}
 \sum_{i=0}^N a_i(x) \rho^{\bar\alpha_i} |v|^{\bar\alpha_i} v+ \frac{2\rho}{\widetilde\phi(x)} \widetilde\Omega \eR\times v +\rho \widetilde\Omega^2 \eR\times (\eR\times x)=-\nabla p+\rho \vec g.
\eeq

Define a function $g:\bar U\times \R_+\to\R_+$  by
\beq\label{eq2}
g(x,s)=a_0(x) + a_1(x)s^{\bar\alpha_1}+\cdots +a_N(x)s^{\bar\alpha_N}=\sum_{i=0}^N a_i(x) s^{ \bar\alpha_i}\quad\text{for } x\in\bar U,\ s\ge 0,
\eeq 
 where $N\ge 1$ and $\bar\alpha_i$ are as in \eqref{Nal}, 
 and the coefficients are functions $a_0(x)>0$, $a_N(x)>0$ and $a_i(x)\ge 0$, for $1\le i\le N-1$, for all $x\in\bar U$.

Set $\mathcal R(x,\rho)=2\rho\widetilde\Omega/\widetilde\phi(x)$.
Multiplying both sides of  \eqref{nGF}  by $\rho$ gives
 \beq\label{neweq1}
 g(x,|\rho v|) \rho v + \mathcal R(x,\rho) \eR \times (\rho v)=-\rho\nabla p+ \rho^2 \vec g- \rho^2 \widetilde\Omega^2 \eR \times (\eR \times x) .
 \eeq 
We solve for $\rho v$ from \eqref{neweq1} in terms of the vector on its right-hand side and the $\mathcal R(x,\rho)$. 
To do that, we  define the function $F_{x,z}:\R^3\to\R^3$, for any $x\in\bar U$ and $z\in\R$, by
\beq\label{Fdef}
F_{x,z}(v)= g(x,|v|)v + z  \mathbf J v \quad \text{ for }v\in\R^3,
\eeq
where $\mathbf J$ is the $3\times 3$ matrix for which $\mathbf J x=\eR \times x$ for all $x\in\R^3$. 
Then equation \eqref{neweq1} can be rewritten as
\beq\label{FRv}
F_{x,\mathcal R(x,\rho)}(\rho v)=-(\rho\nabla p - \rho^2( \vec g - \widetilde\Omega^2 \mathbf J^2 x) ).
\eeq
Thanks to \cite[Lemma 1.1]{CHK3}, the function $F_{x,z}$ is odd and bijective for each $x\in\bar U$ and $z\in\R$.
Therefore, we can invert \eqref{FRv} to have
 \beq\label{new2}
\rho v= -(F_{x,\mathcal R(x,\rho)})^{-1}(\rho\nabla p - \rho^2 ( \vec g -\widetilde\Omega^2 \mathbf J^2 x)).
\eeq 

We will combine equation \eqref{new2} with the conservation of mass, see e.g. \cite[Eq. (6.2.4)]{BearBook},
\beq\label{eq5}
\widetilde \phi(x)\frac{\partial \rho}{\partial t} +\nabla\cdot(\rho v)=0. 
\eeq

Two types of fluid flows are considered below: one is the isentropic gas flows and the other is 
flows of slightly compressible fluids.

\medskip\noindent
\emph{Isentropic gas flows.} In this case, the equation of state is
\beq\label{gas}
p=c\rho^\gamma\quad\text{for some positive constants } c,\gamma.
\eeq
It follows from \eqref{gas} that
\beq\label{ru1}
 \rho\nabla p=\nabla u\quad \text{with } u=\frac{c\gamma\rho^{\gamma+1}}{\gamma+1}.
\eeq
For the sake of convenience, we call the new quantity $u$ the pseudo-pressure.
We rewrite the last equation as
\beq\label{maineq1}
\rho=\left(\frac{\gamma+1}{c\gamma}\right)^\frac1{\gamma+1} u^\lambda = \bar c u^\lambda
\text{ with }
\lambda=\frac1{\gamma+1}\in (0,1) \text{ and }\bar c=\left(\frac{\gamma+1}{c\gamma}\right)^\frac1{\gamma+1}.
\eeq 

The continuing calculations below, up to \eqref{uX1}, only require the following two properties
\beq\label{require}
\rho\nabla p=\nabla u\text{ and } \rho=\bar c u^\lambda
\eeq
which, under the current consideration, are guaranteed by \eqref{ru1} and \eqref{maineq1}.

The gravitational field $\vec g$ in the rotating frame is 
\beqs 
\vec g(t)=-\widetilde{\mathcal G} \widetilde e_0(t),
\eeqs 
where $\widetilde{\mathcal G}>0$ is the gravitational constant and $\widetilde e_0\in C^\infty(\R,\R^3)$ with  $|\widetilde e_0(t)|=1$ for all $t\in\R$.
For the sake of convenience, we use the following scaled function and numbers
\beq \label{phGO}
\phi(x) =\bar c\widetilde \phi(x)>0,\quad
\mathcal G=\bar c^2\widetilde{\mathcal G},\quad \Omega=\bar c \widetilde \Omega.
\eeq
Then
\beq\label{Zx1}
- \rho^2 ( \vec g -\widetilde\Omega^2 \mathbf J^2 x)=u^{2\lambda} \mathcal Z(x,t),
\text{ where }
\mathcal Z(x,t)=\mathcal G \widetilde e_0(t)+ \Omega^2  \mathbf J^2 x.
\eeq
Note that
\beq\label{Rstar}
\mathcal R(x,\rho)=R_*(x) u^\lambda, \text{ where }R_*(x)=2\bar c \widetilde\Omega/\widetilde\phi(x)=2\bar c\Omega/\phi(x).
\eeq
With \eqref{require}, \eqref{Zx1} and \eqref{Rstar}, we rewrite \eqref{new2} as 
\beq\label{rvu}
\rho v= -X(x,u,\nabla u+u^{2\lambda} \mathcal Z(x,t)),
\eeq
where
\beq\label{Xdef}
X(x,z,y)=F_{x,R_*(x)z^\lambda}^{-1}(y) \text{ for }x\in\bar U,\ z\in\R_+,\ y\in\R^3.
\eeq
Then we obtain from \eqref{eq5}, \eqref{require}, \eqref{phGO}  and \eqref{rvu} that
\beq\label{uX1}
\phi (u^\lambda)_t=  \nabla \cdot\left (X\left(x,u,\nabla u+u^{2\lambda} \mathcal Z(x,t)\right)\right).  
\eeq

\medskip\noindent
\emph{Flows of slightly compressible fluids.} The equation of state is
   \beq\label{slight}
\frac{1}{\rho}   \frac{d\rho}{dp}=\varpi, 
\quad \text{where the constant compressibility $\varpi>0$ is small}.
   \eeq
Note, by \eqref{slight}, that $\rho\nabla p=\varpi^{-1}\nabla \rho=\nabla u$, where $u=\rho/\varpi$.
Thus, we have \eqref{require} with 
   \beq\label{cb2}
   \lambda=1\text{ and }\bar c=\varpi.
   \eeq
Therefore, we obtain equations \eqref{rvu}, \eqref{Xdef} and \eqref{uX1} again with the parameters in \eqref{cb2}.

\medskip
In summary,  we obtain, for both cases of isentropic gas flows  and slightly compressible fluids, the unified equation \eqref{uX1} with 
$\lambda\in(0,1]$.

We will characterize the effect of rotation by the following parameter, see \eqref{phGO} together with $\bar c$ in \eqref{maineq1} and \eqref{cb2},
\beq \label{chidef}
\chi_*=1+\Omega=1+\bar c \widetilde \Omega.
\eeq 
Using this, we can estimate the function $\mathcal Z$ in \eqref{Zx1} by
\beq\label{Zesto}
| \mathcal Z(x,t)|\le \chi_*^2 (\mathcal G+|x|)\le c_{\mathcal Z}\chi_*^2, \quad c_{\mathcal Z}=\mathcal G+\max\{|x|:x\in \bar U\}.
\eeq

\medskip
We now turn to the boundary condition and consider a mixed mass and volumetric flux  condition. 
Suppose the boundary $\Gamma$ of $U$ is of class $C^1$. 
Denote the outward normal vector on  $\Gamma$ by $\vec{\nu}$. 
Assume $\Gamma=\Gamma_1\cup\Gamma_2$ with $\Gamma_1\cap \Gamma_2=\emptyset$ and
\beq\label{mixedbc}
\rho v\cdot \vec\nu =\widetilde \psi_1\text{ on }\Gamma_1\text{ and }
 v\cdot \vec\nu =\widetilde \psi_2\text{ on }\Gamma_2,
\eeq
where the function $\widetilde \psi_i(x,t)$ is defined on $\Gamma_i\times(0,\infty)$ for $i=1,2$.
Multiplying the second equation in \eqref{mixedbc} by $\rho$ yields
\beq\label{rvbc2}
\rho v\cdot \vec\nu =\widetilde \psi_2\rho\text{ on }\Gamma_2.
\eeq
Using formula \eqref{rvu} for $\rho v$ in both \eqref{mixedbc} and \eqref{rvbc2}, and using $\rho=\bar c u^\lambda$ from \eqref{require} for the right-hand side of \eqref{rvbc2}, we obtain
\beq \label{separate}
\begin{aligned}
X\left(x,u,\nabla u+u^{2\lambda} \mathcal Z(x,t)\right)\cdot \vec\nu +\widetilde \psi_1=0\text{ on }\Gamma_1,\\
 X\left(x,u,\nabla u+u^{2\lambda} \mathcal Z(x,t)\right)\cdot \vec\nu +\bar c \widetilde \psi_2 u^\lambda=0 \text{ on }\Gamma_2.
\end{aligned}
\eeq 
For $i=1,2$, let $\psi_i(x,t)$ to be the extension of $\widetilde \psi_i$ to $\Gamma\times (0,\infty)$ with the zero value when $x\not\in \Gamma_i$. We combine the two  equations in \eqref{separate} into the following one equation
\beq\label{onebc}
X\left(x,u,\nabla u+u^{2\lambda} \mathcal Z(x,t)\right)\cdot \vec\nu + \psi_1 + \bar c \psi_2 u^\lambda=0\text{ on }\Gamma.
\eeq

With the function $X$ being defined by \eqref{Fdef}, \eqref{Rstar} and \eqref{Xdef}, its main properties are the following.
Set
\beqs 
a=\frac{\bar\alpha_N}{1+\bar\alpha_N}\in(0,1),
\eeqs
and, referring to $g(x,s)$ in \eqref{eq2},
\beq\label{chichi}
\chi_0(x)=g(x,1)=\sum_{i=0}^N a_i(x)\text{ and }\chi_1(x,z)=\chi_0(x)+R_*(x)z^{\lambda}
\text{ for }  x\in\bar U,\ z\in\R_+.
\eeq 

\begin{lemma}\label{Xlem}
For all $x\in\bar U$, $z\in\R_+$ and $y\in \mathbb R^3$, one has 
\beq\label{Xo1}
|X(x,z,y)|\le \bar{W}_0(x) |y|^{1-a},
\eeq
\beq\label{Xo2}
X(x,z,y)\cdot y\ge  \frac{\bar{W}_1(x)|y|^{2-a}}{\chi_*^2(1+z)^{2\lambda}}-\bar{W}_2(x), 
\eeq
where
\beqs
\bar{W}_0(x)=(a_N(x))^{a-1},\  
\bar{W}_1(x)=
\frac{2^{-a}\widetilde c_4(x)}{(1+2\bar c)^2(\chi_0(x)+\phi(x)^{-1})^2},\ 
\bar{W}_2(x)=\frac{2^{-a}\widetilde c_4(x)}{\chi_0^2(x)},
\eeqs 
with $\widetilde c_4(x)=(\min\{1,a_0(x),a_N(x)\}/2^{\bar\alpha_N})^{1+a}$.
\end{lemma}
\begin{proof}
 Recall that \cite[Lemma 2.1]{CHK3} contains estimates for $|X(x,z,y)|$ and $X(x,z,y)\cdot y$ when $a_i$ and $\mathcal R$ are constants. We apply those estimates to $a_i(x)$ and $\mathcal R=R_*(x)z^\lambda$, see \eqref{Rstar}, hence,  replace the constants $\chi_0$ and $\chi_1$ in \cite[(2.16)]{CHK3} with $\chi_0(x)$ and $\chi_1(x,z)$ in \eqref{chichi} above.
Specifically, on the one hand, we have from the second inequality of  \cite[Ineq. (2.18)]{CHK3} that 
\beqs
|X(x,z,y)|\le (a_N(x))^{a-1} |y|^{1-a}
\eeqs
which proves \eqref{Xo1}. 
On the other hand, we have from the first inequality of \cite[Ineq. (2.20)]{CHK3} that
\beq\label{X3}
X(z,z,y)\cdot y\ge \frac{2^{-a}\widetilde c_4(x)}{\chi_1^2(x,z)} (|y|^{2-a}-1)
=\frac{2^{-a}\widetilde c_4(x) |y|^{2-a}}{\chi_1(x,z)^2}-\frac{2^{-a}\widetilde c_4(x) }{\chi_1^2(x,z)}.
\eeq 
For the second to last term in \eqref{X3}, we use the following upper bound of $\chi_1(x,z)$, recalling  the formula of  $R_*(x)$ from \eqref{Rstar},
\begin{align*}
\chi_1(x,z)= \chi_0(x)+\frac{2\bar c \Omega}{\phi(x)} z^{\lambda}
&\le (\chi_0(x)+\phi(x)^{-1})(1+2\bar c)(1+\Omega ) (1+z)^\lambda\\
&=(1+2\bar c) \chi_*(\chi_0(x)+\phi(x)^{-1}) (1+z)^\lambda.
\end{align*}
For the last term in \eqref{X3}, we use the lower bound $\chi_1(x,z)\ge  \chi_0(x)$.
With  these bounds of $\chi_1(x,z)$, we obtain inequality \eqref{Xo2} from \eqref{X3}.
\end{proof}

Having in mind the PDE \eqref{uX1}, the boundary condition \eqref{onebc}, Lemma \ref{Xlem} and property \eqref{Zesto},  we formulate our initial boundary value problem with a slight generalization. Although our physical problem is for the spatial dimension $n=2$ or $3$, we consider hereafter any integer $n\ge 2$, 
and the set $U$ now is non-empty, open, bounded in $\R^n$.
We study the initial boundary value problem 
\begin{subequations}\label{ibvpg}
\begin{align}
\label{pde}
& \phi (u^\lambda)_t=\nabla\cdot \Big(X\big(x,u,\nabla u+u^{2\lambda}  \mathcal Z(x,t)\big)\Big) &&\text{in } U\times (0,\infty),\\
\label{ic}
& u(x,0)=u_0(x) &&\text{in } U,\\
\label{bc}
& X\big(x,u,\nabla u+u^{2\lambda}  \mathcal Z(x,t)\big)\cdot \vec\nu + \psi_1(x,t)+ \psi_2(x,t)u^{\lambda} =0 &&\text{on }\Gamma\times (0,\infty),
\end{align}
\end{subequations}

where 
\begin{enumerate}[label={\rm (H\arabic*)}]
\item\label{H1} the unknown  is a function $u:\bar U\times \R_+\to \R_+$, 
    \item $\lambda$ is a positive constant,
    \item function $\phi:U\to(0,\infty)$ is given, 
    \item\label{XH}  function $X:\bar U\times\R_+\times\R^n\to\R^n$ is given that satisfies 
    \beq\label{XX1}
|X(x,z,y)|\le W_0(x) |y|^{1-a},
\eeq
\beq\label{XX2}
X(x,z,y)\cdot y\ge  \frac{W_1(x)|y|^{2-a}}{\chi_*^2(1+z)^{2\lambda}}-W_2(x), 
\eeq
for all $x\in\bar U$, $z\in\R_+$, $y\in \mathbb R^n$, and some positive functions $W_0$, $W_1$ and $W_2$ defined on $\bar U$,
    \item $\mathcal Z:\bar U\times (0,\infty)\to \R^n$ is a given function satisfying
\beq\label{Zest}
| \mathcal Z(x,t)|\le c_{\mathcal Z}\chi_*^2\text{ for all $x\in \bar U$, $t>0$ and some constant }c_{\mathcal Z}>0,
\eeq
    \item $u_0:U\to\R_+$ is a given initial data,
    \item\label{H7} $\psi_1$ and $\psi_2$ are given functions  mapping $\Gamma\times(0,\infty)$ to $\R$.
\end{enumerate}

All the Hypotheses \ref{H1}--\ref{H7} are applicable to the physical problems discussed from \eqref{Fv} to \eqref{slight} and the boundary condition \eqref{mixedbc}.

\begin{remark}\label{difu}
The following remarks on the degeneracy/singularity with respect to $x$, $u$, and $\nabla u$ of the nonlinear PDE \eqref{pde}  are in order. (The degeneracy and singularity here are understood in comparison with the standard parabolic equation.)

\begin{enumerate}[label=\rnum]
    \item\label{rma} \emph{Degeneracy/singularity in $x$.}
On the left hand side of \eqref{pde}, if the general positive function $\phi(x)$ is close to $0$, resp., $\infty$, then the equation is degenerate, resp., singular.
Concerning the right-hand side of equation \eqref{pde}, having in mind  \eqref{XX1} and \eqref{XX2}  in \ref{XH}, the equation is degenerate if $W_1(x)$ is close to $0$, and singular if $W_1(x)$ and/or $W_0(x)$ are close to $\infty$. 

\item\label{rmb}  \emph{Degeneracy/singularity in $u$.}
We formally rewrite \eqref{pde} as
\beq\label{rewr1}
\phi(x)\lambda u^{\lambda-1}u_t =\nabla \cdot X(x,u,\nabla u+u^{2\lambda}\mathcal{Z}(x,t)).
\eeq 
If $\lambda >1$, then the left-hand side of \eqref{rewr1} is degenerate when $u\to 0$.
If $\lambda <1$, then the left-hand side of \eqref{rewr1} is singular, resp., degenerate when $u\to 0$, resp., $u\to \infty$.

Moreover, observe from \eqref{XX2} that the factor $1/(1+z)^{2\lambda}$ goes to zero as $z\to\infty$. Therefore, $X(x,u,\nabla u+u^{2\lambda}\mathcal{Z}(x,t))$ in \eqref{pde} also poses a degeneracy as $ u\to\infty$.

\item\label{rmc} \emph{Degeneracy in $\nabla u$.}  
Regarding $\nabla u$ in $X(x,u,\nabla u+u^{2\lambda}\mathcal{Z}(x,t))$ on the right-hand side of \eqref{pde}, by \eqref{XX1} and \eqref{XX2}, the equation is degenerate as $|\nabla u|\to\infty$. 
\end{enumerate}
\end{remark}

\begin{remark}\label{pLrmk}
We recall that the standard $p$-Laplacian equation, for $p>1$, is
\beq\label{pLa}
u_t=\nabla\cdot X_p(\nabla u),\text{ where } X_p(\xi)=|\xi|^{p-2}\xi \text{ for }\xi\in\R^n.
\eeq
For large  $|\xi|$ the function $X(x,z,\xi)$ in equation \eqref{pde} has a similar growth in $\xi$ to the function $X_p$ for $p=2-a\in (1,2)$,  but not for $p>2$, see \eqref{XX1} and \eqref{XX2}. However, \eqref{pde} is not comparable with \eqref{pLa} because $X$, unlike $X_p$,   is not homogeneous in $\xi$, and it depends also on $x$ and $z$. Moreover,  its third argument in \eqref{pde} is $\nabla u+u^{2\lambda}\mathcal{Z}$ not just only $\nabla u$.
More differences between \eqref{pde} and \eqref{pLa} can bee seen in \ref{rma} and \ref{rmb} of Remark \ref{difu}.
\end{remark}

\section{Preliminaries}\label{Prem}
For $z\in\R$, the positive and negative parts of $z$ are $z^+=\max\{0,z\}$ and $z^-=\max\{0,-z\}$, respectively.
We recall some elementary inequalities that will be used in this paper.
\beq\label{ee3}
(x+y)^p\le 2^{(p-1)^+}(x^p+y^p)\quad  \text{for all } x,y\ge 0,\quad  p>0.
\eeq
Particularly,
\beq\label{ee2}
(x+y)^p\le 2^p(x^p+y^p)\quad  \text{for all } x,y\ge 0,\quad p>0.
\eeq
By the triangle inequality and inequality \eqref{ee3}, we have
\beq\label{ee6}
|x-y|^p\ge 2^{-(p-1)^+}|x|^p-|y|^p \quad \text{for all } x,y\in\R^n,\quad p>0.
\eeq
We will use frequently the following consequences of Young's inequality
\beq\label{ee4}
x^\beta \le x^\alpha+x^\gamma\quad \text{for all } (x>0,\  \gamma\ge \beta\ge\alpha) \text{ or } 
(x=0,\  \gamma\ge \beta\ge\alpha>0),
\eeq
\beq\label{ee5}
x^\beta \le 1+x^\gamma \quad \text{for all } 
x\ge 0,\   \gamma\ge \beta\ge 0 \text{ with }x^2+\beta^2>0.
\eeq

For a Lebesgue measurable function $\varphi:U\to (0,\infty)$, a number $p\ge 1$ and functions $f:U\to \R$, $g:U\times(T_1,T_2)\to \R$ with $T_1<T_2$, we define
\beqs
\|f\|_{L^p_\varphi(U)}=\left(\int_U |f(x)|^p\varphi(x)\d x\right)^{1/p},\ 
\|g\|_{L^p_\varphi(U\times(T_1,T_2))}=\left(\int_{T_1}^{T_2}\int_U |g(x,t)|^p\varphi(x)\d x\d t\right)^{1/p}.
\eeqs

\subsection{Sobolev inequalities with multiple weights}

First, we recall the standard Sobolev inequality without weights.
For a number $p\in[1,n)$, denote $p^*=np/(n-p)$. 
Then there exists a number $\widehat c_p>0$ depending on $U$, $n$ and  $p$, such that
\beq\label{stdSov}
\|f\|_{L^{p^*}(U)}\le \widehat c_p \left(\int_U |\nabla f(x)|^p\d x+\int_U |f(x)|^p\d x\right)^{1/p}\text{ for any  function } f\in W^{1,p}(U).
\eeq

In all Lemmas \ref{WS1}--\ref{RoS2} below, let numbers $p>1$ and $r_1>0$ satisfy 
\beq\label{rone}
\frac{n}{n+p}<r_1<1\le r_1p<n,
\eeq 
and numbers $r\ge 0$ and $s \in\R$ be given.
Denote
\beq\label{rstar}
r_*= 1+\frac p n-\frac1{r_1}\in(0,1).
\eeq

Below, the functions  $\varphi(x)>0$ in all Lemmas \ref{WS1}--\ref{RoS2}, 
$\omega(x)\ge 0$ in Lemmas \ref{WS1} and \ref{RoS1},
$W(x)>0$ in Lemmas \ref{WS1},  \ref{trace} and \ref{WS2}, 
and $W_*(x)>0$ in Lemmas \ref{RoS1}, \ref{Rtrace} and \ref{RoS2}  are defined for $x\in U$.

Next, we recall a Sobolev inequality with weights from our previous work \cite{CHK5}.

\begin{lemma}[{\cite[Lemma 2.1]{CHK5}}]\label{WS1}
 Assume $\alpha\in\R$ satisfies
\beq \label{manya}
\alpha\ge s,\quad \alpha > \frac{p-s}{p-1},\quad  \alpha> \frac{2(r+s-p)}{r_*}.
\eeq  
Denote
\beq\label{powdef}
\begin{aligned}
m&= \frac{\alpha-s+p}{p}\in[1,\alpha), \quad
\theta=\frac{\alpha+2r}{\alpha(1+r_*)+2(p-s)}\in(0,1),\\
 \mu_1&=\frac{r+\theta(s-p)}{1-\theta}\in (-\alpha,\infty).
\end{aligned}
\eeq
Define
\beq \label{Gdef}
\begin{aligned}
    G_1&=\left(\int_U \varphi(x)^{-\frac {\alpha-s+p}{\alpha(p-1)+s-p}} \d x\right)^\frac{\alpha(p-1)+s-p}{\alpha},\quad 
        G_2=\left(\int_U W(x)^{-\frac{r_1}{1-r_1}}\d x\right)^{\frac{1-r_1}{r_1}},  \\
  G_3&=\left(\int_U  \varphi(x)^{-1}\omega(x)^{\frac{1}{(1-\theta)(1+\mu_1/\alpha)}}\d x\right)^{1+\mu_1/\alpha},\quad
  \Phi_1= G_1^\theta G_3^{1-\theta},\quad 
   \Phi_2=  G_2^{\frac{ \theta}{1-\theta}} G_3.
\end{aligned}
\eeq 
Then  one has, for any function $u\in W^{1,{r_1pm}}(U)$ and number $\varepsilon>0$, that
\beq \label{S10}
\begin{aligned}
\int_U |u(x)|^{\alpha+r} \omega(x) \d x 
& \le \varep \int_U |u(x)|^{\alpha-s}|\nabla u(x)|^p W(x)\d x\\ 
&\quad  +D_{1,m,\theta} \Phi_1 J_*^{1+r/\alpha}
 +\varep^{-\frac\theta{1-\theta}}D_{2,m,\theta} \Phi_2 J_*^{1+\mu_1/\alpha},
\end{aligned}
\eeq 
where $J_*=\int_U |u(x)|^\alpha\varphi(x)\d x$ and 
\beq \label{dthe}
D_{1,z,\eta}=(c_4 2^{z})^{\eta p} ,\quad 
D_{2,z,\eta}=(c_3z 2^{z})^\frac{\eta p}{1-\eta}
\text{ for $z>0$, $\eta\in(0,1)$}
\eeq 
with $c_3$, $c_4$ being some positive numbers depending on  $U$, $n$, $r_1$ and $p$.
\end{lemma}

Note from the first two inequalities in \eqref{manya}, by considering $s>1$ and $s\le 1$, that
\beq\label{ag1}
\alpha>1.
\eeq

Although already complicated, the above Lemma \ref{WS1} cannot be used directly to treat the new nonlinear structure appearing in equation \eqref{uX1}. However, based on this lemma, we derive a new form of inequality which is more appropriate to our problem.

\begin{lemma}\label{RoS1}
We consider the same assumptions and use the same notation as in Lemma \ref{WS1}. 
Let $\beta$ be any positive number  and assume additionally that 
\beq\label{abeta}
\alpha >\frac{r_1\beta}{1-r_1}.
\eeq 
Define 
\beq \label{Calpha}
\begin{aligned}
{\widetilde G}_2&=\left(\int_U \varphi(x)^{-\frac{r_1\beta}{\alpha(1-r_1)-r_1\beta}}W_*(x)^{-\frac{r_1\alpha}{\alpha(1-r_1)-r_1\beta}}\d x\right)^{\frac{\alpha(1-r_1)-r_1\beta}{\alpha r_1}},\\
{\mathcal E}_1&= 1+\int_U \varphi(x) \d x,
\quad 
\widetilde\Phi_2 ={\widetilde G}_2^\frac\theta{1-\theta}G_3 {\mathcal E}_1^{\frac{\beta\theta}{\alpha(1-\theta)}}.
\end{aligned}
\eeq
Set the numbers
\beq\label{hatmu1}
\widehat\mu_1=\mu_1+\frac{\beta\theta}{1-\theta},\quad 
\widehat\mu_2=\min\{r,\mu_1\},\quad 
\widehat\mu_3=\max\{r,\widehat\mu_1\}.
\eeq
Then one has, for any number $\varep>0$,
\beq \label{RSi2}
\begin{aligned}
\int_U |u(x)|^{\alpha+r} \omega(x) \d x 
& \le \varep \int_U \frac{|u(x)|^{\alpha-s}}{(1+|u(x)|)^\beta}|\nabla u(x)|^p W_*(x)\d x \\
&\quad  +{\widehat\Phi}_1
\min\left\{ J_*^{1+\widehat\mu_2/\alpha}+ J_*^{1+\widehat\mu_3/\alpha},\left(1+J_* \right)^{1+\widehat\mu_3/\alpha}
\right\},
\end{aligned}
\eeq 
where $J_*=\int_U |u(x)|^\alpha\varphi(x)\d x$, and
\beq\label{phi00}
{\widehat\Phi}_1=D_{1,m,\theta} \Phi_1 
+2^{1+\widehat\beta_1}  \varep^{-\frac\theta{1-\theta}}D_{2,m,\theta} \widetilde\Phi_2  
\text{ with }\widehat\beta_1=\frac{\beta\theta}{1-\theta}\left(1+\frac1\alpha\right).
\eeq
\end{lemma}

The proof of Lemma \ref{RoS1} will be given in Appendix \ref{Appex}.


\subsection{Trace theorem with weights}
Recall the standard trace theorem without weights
\beq\label{firstrace}
\int_\Gamma |f(x)|\d S \le c_5\int_U |f(x)|\d x +c_6\int_U |\nabla f(x)|\d x
\eeq
for all $f \in W^{1,1}(U)$, where $c_5$ and $c_6$ are positive constants depending on $U$.

Below is a version with two different weights and general powers which was proved in \cite{CHK5}.

\begin{lemma}[{\cite[Lemma 2.3]{CHK5}}]\label{trace}
 Assume $\alpha$ is a number that satisfies 
 \beq \label{manyb}
\alpha\ge s\text{ and }\alpha > \frac{p-s}{p-1}.
\eeq 
Let the numbers $m$, $\theta$ and $\mu_1$ be defined as in \eqref{powdef}, and set
\beq\label{rtilde} 
\widetilde r=r+\frac{r+s-p}{p-1}=\frac{rp+s-p}{p-1}.
\eeq
Let $c_5, c_6$ be as in \eqref{firstrace} and 
 $D_{1,z,\eta}, D_{2,z,\eta}$ as in \eqref{dthe}.
Define  $G_1$, $G_2$ as in \eqref{Gdef}, and 
\beq\label{Ph34}
    G_4=\left(\int_U  \varphi(x)^{-1} \d x\right)^{1+\mu_1/\alpha},\quad 
    \Phi_3=G_1^\theta G_4^{1-\theta}, \quad
    \Phi_4=G_2^{\frac{ \theta}{1-\theta}} G_4.
\eeq

Let $\varep>0$ and $u(x)$ be any function on $U$ that satisfies $|u|^\alpha\in W^{1,1}(U)$ and $|u|^m\in W^{1,p}(U)$.
Denote $J_*=\int_U |u(x)|^\alpha \varphi(x)\d x$.

\begin{enumerate}[label=\tnum]
    \item\label{trneg} Case $\widetilde r <0$. Suppose $\alpha>|\widetilde r|$.
Let 
\beq \label{Ph5}   
  \Phi_5=\left(\int_U W(x)^{\frac {\alpha}{(p-1)\widetilde r}}\varphi(x)^{\frac{\alpha+\widetilde r}{\widetilde r}}\d x\right)^{-\widetilde r/\alpha}.
\eeq
Then one has
\beq \label{trace0}
\begin{aligned}
&\int_\Gamma |u(x)|^{\alpha+r} \d S 
\le 2\varepsilon \int_U |u(x)|^{\alpha-s}|\nabla u(x)|^p W(x) \d x\\
&\quad  +z_1 \Phi_3 J_*^{1+ r/\alpha} +\varep^{-\frac\theta{1-\theta}} z_2 \Phi_4J_*^{1+\mu_1/\alpha}
  + \varepsilon^{-\frac 1{p-1}}  z_3 \Phi_5 J_*^{1+ \widetilde r/\alpha},
\end{aligned}
\eeq 
where 
\beq \label{zzz}
z_1=c_5 D_{1,m,\theta},\quad 
z_2= c_5^{\frac{1}{1-\theta}}D_{2,m,\theta},\quad 
z_3= (c_6 (\alpha+r))^\frac p{p-1}.
\eeq 

\item\label{trpos} Case $\widetilde r \ge 0$.
Suppose 
\beq\label{only}  
\alpha> \frac{2(\widetilde r+s-p)}{r_*} =\frac{2p(r+s-p)}{r_*(p-1)}.
\eeq 
Let 
\beq\label{tmtilde} 
\widetilde\theta=\frac{\alpha+2\widetilde{r}}{\alpha(1+r_*)+2(p-s)}\in(0,1), \quad 
 \widetilde\mu_1=\frac{\widetilde r+\widetilde\theta(s-p)}{1-\widetilde\theta}\in(-\alpha,\infty),
\eeq
\beq \label{Ph67}
G_5=\left(\int_U  \varphi(x)^{-1}W(x)^{-\frac{1}{(p-1)(1-\widetilde\theta) (1+\widetilde\mu_1/\alpha)}}\d x\right)^{1+\widetilde\mu_1/\alpha},\  
\Phi_6=  G_1^{\widetilde\theta} G_5^{1-\widetilde\theta}, \ 
\Phi_7=  G_2^{\frac{\widetilde\theta}{1-\widetilde\theta}} G_5.
\eeq
Then one has
\begin{align} \label{trace1}
&\int_\Gamma |u(x)|^{\alpha+r} \d S 
\le 3\varepsilon \int_U |u(x)|^{\alpha-s}|\nabla u(x)|^p W(x) \d x 
+ z_1 \Phi_3 J_*^{1+r/\alpha} \notag \\
&\quad + \varep^{-\frac\theta{1-\theta}}z_2 \Phi_4 J_*^{1+\mu_1/\alpha}
+ \varepsilon^{-\frac 1{p-1}}z_4 \Phi_6
  J_*^{1+\widetilde r/\alpha}+ \varep^{-(\frac 1 {p-1}+\frac p{p-1}\cdot \frac{\widetilde\theta}{1-\widetilde\theta} )}z_5 \Phi_7
  J_*^{1+\widetilde\mu_1/\alpha},
\end{align}
where $z_1$, $z_2$, $z_3$ are from \eqref{zzz}, $z_4=z_3 D_{1,m,\widetilde\theta}$ and $z_5=z_3^\frac1{1-\widetilde\theta} D_{2,m,\widetilde\theta}$.
\end{enumerate}
\end{lemma}

The following is a new trace theorem which, again, takes into account the new nonlinearity in \eqref{ibvpg}.

\begin{lemma}\label{Rtrace}
We use the same notation as in Lemma \ref{trace}, and consider only case \ref{trpos}, i.e. when $\widetilde r\ge 0$ with all the needed assumptions up to \eqref{Ph67}.     
Assume, additionally, that
\beq\label{bcond}
\alpha>\beta_*\eqdef \frac{\beta}{(p-1)(1-\widetilde\theta) (1+\widetilde\mu_1/\alpha)}.
\eeq
Set
\beq \label{muh45}
\begin{aligned}
    \widehat\beta_2&=\frac{\beta}{p-1}\left(1+\frac1\alpha\right), &&\widehat\beta_3=\frac{\beta}{1-\widetilde\theta}\left(\frac1{p-1}+\widetilde\theta\right),\\
\widehat\mu_4&=\min\{r,\mu_1,\widetilde r, \widetilde\mu_1\},&&
\widehat\mu_5=\max\left \{r, \widehat\mu_1, \widetilde r+\frac{\beta}{p-1}, \widetilde\mu_1+\widehat\beta_3\right\},
\end{aligned}
\eeq 
and
    \beq\label{tfifour} \begin{aligned}
    \widetilde \Phi_4&= {\mathcal E}_1^{\frac{\beta\theta}{\alpha(1-\theta)}}{\widetilde G}_2^\frac\theta{1-\theta} G_4,\quad  
    &{\widetilde G}_5&=\left(\int_U \varphi(x)^{-\frac{\alpha+\beta_*}{\alpha-\beta_*}} W_*(x)^{-\frac{\beta_*}{\beta}\cdot \frac{\alpha}{\alpha-\beta_*}}\d x\right)^{\frac{\alpha-\beta_*}{\alpha}(1+\widetilde\mu_1/\alpha)},\\ 
   {\widetilde\Phi}_6&={\mathcal E}_1^{\frac{\beta}{\alpha(p-1)}}G_1^{\widetilde\theta} {\widetilde G}_5^{1-\widetilde\theta},\  
    &{\widetilde\Phi}_7&= {\mathcal E}_1^{\frac{\widehat\beta_3}{\alpha}} {\widetilde G}_2^{\frac{\widetilde\theta}{1-\widetilde\theta}}{\widetilde G}_5.
\end{aligned}\eeq 
Then one has 
\beq\label{tru3}
\begin{aligned}
\int_\Gamma |u(x)|^{\alpha+r} \d S 
&\le 3\varepsilon \int_U \frac{|u(x)|^{\alpha-s}}{(1+|u(x)|)^\beta}|\nabla u(x)|^p W_*(x)\d x 
 \\
&\quad + {\widehat\Phi}_2 \min\left \{ J_*^{1+\widehat\mu_4/\alpha}+ J_*^{1+\widehat\mu_5/\alpha},\left (1+  J_* \right)^{1+\widehat\mu_5/\alpha}\right\},
\end{aligned}
\eeq
where $J_*=\int_U |u(x)|^\alpha\varphi(x)\d x$, 
\beq\label{phi20}
{\widehat\Phi}_2
=z_1 \Phi_3+ \varep^{-\frac\theta{1-\theta}}\widetilde z_2 {\widetilde\Phi}_4 
+ \varepsilon^{-\frac 1{p-1}}\widetilde z_4 \widetilde\Phi_6
+ \varep^{-(\frac 1 {p-1}+\frac p{p-1}\cdot \frac{\widetilde\theta}{1-\widetilde\theta} )}\widetilde z_5 \widetilde\Phi_7
\eeq 
with
\beq\label{ztils}
\widetilde z_2=2^{1+\widehat\beta_1}  z_2,\ 
\widetilde z_4=2^{1+\widehat\beta_2} z_4,\ 
\widetilde z_5=2^{1+\widehat\beta_3(1+\frac1\alpha)} z_5.
\eeq 
\end{lemma}
The proof of Lemma \ref{Rtrace} will be given in Appendix \ref{Appex}.

\begin{remark}\label{absok}
Concerning condition \eqref{bcond}, it is equivalent to 
\beq\label{abeuqiv}
\alpha> \frac{\beta}{(p-1)(1-\widetilde\theta)}-\widetilde\mu_1.
\eeq
Note from \eqref{tmtilde} that
\beqs 
\lim_{\alpha\to\infty}\widetilde\theta=\frac{1}{1+r_*}\text{ and }
 \lim_{\alpha\to\infty}\widetilde\mu_1=\frac{\widetilde r+(s-p)/(1+r_*)}{1-1/(1+r_*)}.
\eeqs
Therefore, \eqref{abeuqiv} and, hence, \eqref{bcond} are satisfied for sufficiently large $\alpha$. 
\end{remark}

\begin{remark}
In Lemma \ref{Rtrace}, a similar inequality to \eqref{tru3} can be obtained for the case $\widetilde r<0$ based on \eqref{trace0}. However, in its proof, after the substitution \eqref{subW} in the term $\Phi_5$ in \eqref{Ph5}, in order to apply H\"older's inequality, see the treatment of $G_5$ from \eqref{Ph4est} to \eqref{Gfiveest}, we need an extra condition $\beta< (p-1)|\widetilde r|$. To shorten  the paper, we omitted it and will choose $r$ such that $\widetilde r\ge 0$ in later applications.
\end{remark}

\subsection{Parabolic Sobolev inequality with multiple weights}

We recall our previous result in \cite{CHK5}.

\begin{lemma}[{\cite[Lemma 2.4]{CHK5}}]\label{WS2}
Define
\beq \label{Jss}
 \Phi_8=\left ( \int_U  \varphi(x)^{\frac{2}{r_*}-1}\d x\right )^\frac{r_*}{2}\left[ \left(\int_U W(x)^{-\frac{r_1}{1-r_1}}\d x \right)^{1-r_1}
+\left(\int_U \varphi(x)^{-\frac{r_1}{1-r_1}}\d x\right)^{1-r_1}\right]^\frac{1}{r_1}.
\eeq 
Assume $\alpha$ is a number such that 
\beq \label{alcond}
    \alpha\ge s,\ 
    \alpha> \frac{p-s}{p-1} \text{ and }
    \alpha >\frac{2(s-p)}{r_*}.
\eeq
Let the number $m$ be as in \eqref{powdef}, and define
\beq
 \label{defkappa}
 \kappa = 1+ \frac{r_*}{2}+\frac{p-s}{\alpha}\in (1,\infty), \ 
\theta_0 =\frac 1{1+\frac{r_*\alpha}{2(\alpha-s+p)} }\in (0,1).
\eeq 
Let $u(x,t)$ be any function on $U\times (0,T)$ such that for almost every $t\in (0,T)$ the function $x\in U\mapsto u(x,t)$ belongs to $W^{1,r_1pm}(U)$.
Then one has
\beq\label{ppsi1}
\begin{aligned}
\|u\|_{L^{\kappa \alpha}_\varphi(U\times(0,T))}
&\le  (c_7^p m^\frac1{r_1} \Phi_8)^\frac{1}{\kappa\alpha}  \Bigg(\int_0^T\int_U |u(x,t)|^{\alpha-s}|\nabla u(x,t)|^p W(x)\d x\d t \\
&\quad +\int_0^T\int_U |u(x,t)|^{\alpha-s+p} \varphi(x)  \d x\d t\Bigg)^\frac{1}{\kappa\alpha}
 \cdot  \essup_{t\in(0,T)} \|u(\cdot,t)\|_{L_\varphi^{\alpha}(U)}^{1-\theta_0},
\end{aligned}
\eeq
where  $c_7$ is the positive  constant $\widehat c_{r_1 p}$ in  \eqref{stdSov}.
\end{lemma}

Regarding $\kappa$ and $\theta_0$ in Lemma \ref{WS2}, there is a special relation \cite[Formula (2.87)]{CHK5} which reads as
\beq \label{tk}
(1-\theta_0)\kappa=r_*/2.
\eeq

Next, we prove a new version of the above weighted parabolic Sobolev inequality.

\begin{lemma}\label{RoS2}
Let the assumptions in Lemma \ref{WS2} hold and we use the same notation there. Assume, additionally, that
\beq\label{albe}
\alpha\ge  \frac{\beta}{1-r_1}.
\eeq
Let 
\beq \label{EEdef}
{\mathcal E}_2=\left ( \int_U  \varphi(x)^{\frac{2}{r_*}-1}\d x\right )^\frac{r_*}{2},\  
{\mathcal E}_3=\left(1+\int_U (1+ \varphi(x)^{-1})^{\frac{ r_1}{1-r_1}} (1+W_*(x)^{-1})^{\frac{r_1}{(1-r_1)^2}}\d x 
\right)^{\frac{1-r_1}{r_1}}.
\eeq 
Then one has
\beq \label{pssi6}
\begin{aligned}
&\|u\|_{L^{\kappa \alpha}_\varphi(U\times(0,T))}
\le {\widehat\Phi}_3^\frac{1}{\kappa\alpha}
 \Bigg(\int_0^T\int_U \frac{|u(x,t)|^{\alpha-s}}{(1+|u(x,t)|)^\beta}|\nabla u(x,t)|^p W_*(x)\d x\d t \\
&\quad +\int_0^T\int_U |u(x,t)|^{\alpha-s+p} \varphi(x)  \d x\d t\Bigg)^\frac{1}{\kappa\alpha}
 \cdot \left(  \essup_{t\in(0,T)} \|u(\cdot,t)\|_{L_\varphi^{\alpha}(U)}^{1-\theta_0}
 + \essup_{t\in(0,T)} \|u(\cdot,t)\|_{L_\varphi^{\alpha}(U)}^{1-\widehat\theta_0}\right),
\end{aligned}
\eeq 
where 
\beq\label{tilph10}
   {\widehat\Phi}_3= 2^{1+\beta+\frac1{r_1}} c_7^p m^\frac{1}{r_1} {\mathcal E}_1^\frac\beta\alpha {\mathcal E}_2 {\mathcal E}_3
   \text{ and } \widehat \theta_0=\theta_0-\frac{\beta}{\kappa\alpha}.
\eeq 
\end{lemma}

The proof of Lemma \ref{RoS2} will be given in Appendix \ref{Appex}.

\section{Estimates of the spatial $L^\alpha$-norm}\label{Laest}
In this section, we focus on a non-negative solution $u(x,t)$ of the initial boundary value problem \eqref{ibvpg} and estimate its spatial $L^\alpha$-norm  for a finite number $\alpha>0$ and small $t\ge 0$. For the rest of this paper, the solution $u(x,t)\ge 0$ has enough regularity so that our calculations are valid.
In particular, we assume that the weak formulation
\beq \label{weakf}
\begin{aligned}
\int_U \phi (u^\lambda)_t\cdot \varphi(x,t)\d x
&=-\int_U \Big(X\big(x,u,\nabla u+u^{2\lambda}  \mathcal Z(x,t)\big)\Big)\cdot\nabla \varphi(x,t)\d x\\
&\quad -\int_\Gamma (\phi_1(x,t)+\psi_2(x,t)u^\lambda)\cdot \varphi(x,t)\d S    
\end{aligned}
\eeq 
holds for any test function $\varphi(x,t)$ we later choose.

Below, we will apply Lemmas \ref{RoS1}, \ref{Rtrace} and \ref{RoS2} to 
\beq\label{psbe}
p=2-a,\quad s=\lambda+1\text{ and }\beta=2\lambda.
\eeq
We prepare for these applications with needed parameters and quantities below.
Let $r_1$ be any given number satisfying 
\beq\label{newro}
\frac{n}{n+2-a}<r_1<1<r_1(2-a). \eeq
Because $n\ge 2>r_1(2-a)=r_1p$, the condition \eqref{rone} is met. 
The number $r_*$ in \eqref{rstar} becomes 
\beq\label{newrs}
r_*=\frac{nr_1+r_1(2-a)-n}{nr_1}=1+\frac{2-a}{n}-\frac{1}{r_1}\in(0,1).
\eeq 
Let $r$ be a fixed number such that  
\beq \label{newrr}
r>\max\left\{0, \lambda(5-4a)-1,\frac{1-a-\lambda}{2-a}\right\}.
\eeq
Denote
\beq \label{newtr}
 \widetilde r=r+\frac{r-1+\lambda+a}{1-a}=\frac{(2-a)r-1+\lambda+a}{1-a}>0.
\eeq 

Recall that the functions $W_0$, $W_1$ and $W_2$ are from Hypothesis \ref{XH} in section \ref{IBVPsec}.
Define the two combined weights that will be used later
\beq\label{W3}
W_3(x)=W_0^{2-a}(x)W_1^{a-1}(x)+W_1(x)
\text{ and }
W_4(x)=W_2(x)+W_3(x).
\eeq
For a sufficiently large number $\alpha$, define 
\begin{align}\label{trtm0}
\theta&=\Theta(\alpha,r)\eqdef\frac{\alpha+2r}{\alpha(1+r_*)+2(1-a-\lambda)}\in(0,1), && 
\widetilde\theta=\Theta(\alpha,\widetilde r)\in(0,1), \\ 
\label{trtm1}
 \mu_1&= \Lambda(r,\theta)\eqdef  \frac{r+\theta(a+\lambda-1)}{1-\theta}\in(-\alpha,\infty),&&
 \widetilde\mu_1=\Lambda(\widetilde r,\widetilde \theta)\in(-\alpha,\infty),\\
 \label{trtm2}
\beta_*&=\frac{2\lambda}{(1-a)(1-\widetilde\theta)(1+\widetilde\mu_1/\alpha)},
\end{align}
and let 
\begin{align}\label{Kz}
{\mathcal K}_0&=1+\int_U \phi(x) \d x,\quad 
{\mathcal K}_1= \left(\int_U  \phi(x)^{-1}\d x\right)^{1+\mu_1/\alpha},\\
{\mathcal K}_2&= \left(\int_U \phi(x)^{-\frac {\alpha+1-\lambda-a}{\alpha(1-a)-1+\lambda+a}} \d x\right)^{\frac{\alpha(1-a)-1+\lambda+a}{\alpha}},\notag
\end{align}
together with the $W_j$-related quantities
\begin{align}\notag
 {\mathcal K}_3&= \left(\int_U  \phi(x)^{-\frac{2\lambda r_1}{\alpha(1-r_1)-2\lambda r_1}}W_1(x)^{-\frac{r_1\alpha}{\alpha(1-r_1)-2\lambda r_1}}\d x \right)^\frac{
 \alpha(1-r_1)-2\lambda r_1}{\alpha r_1},\\
 \notag
 {\mathcal K}_4&= \left(\int_U  \phi(x)^{-\frac{\alpha+\beta_*}{\alpha-\beta_*}}W_1(x)^{-\frac{\beta_*\alpha}{2\lambda(\alpha-\beta_*)}}\d x \right)^\frac{(\alpha-\beta_*)(1+\widetilde\mu_1/\alpha)}{\alpha},\\
 {\mathcal K}_5&= \int_U \phi(x)^{-\frac{\alpha-\lambda-1}{\lambda+1}}W_4^{\frac{\alpha}{\lambda+1}}\d x, \quad
\mathcal{K}_6=\int_U W_3(x)^{\frac{\alpha+r}{r-\lambda(5-4a)+1}} \d x.\label{Kf}
\end{align}

\begin{proposition}\label{Diff4u}
There exist positive  numbers $\alpha_*$, $\mu_*$, $\gamma_*$  such that the following statement holds true. If $\alpha>\alpha_*$ and the ${\mathcal K}_j$, for $0\le j\le 6$, are finite numbers, then there is a positive number $\bar C_\alpha$ depending on $\alpha$  and ${\mathcal K}_j$ such that, for $t>0$,
\beq\label{J0dbest}
\begin{aligned}
  & \ddt\int_U u(x,t)^\alpha\phi(x) \d x  +\frac{\alpha(\alpha-\lambda)\chi_*^{-2}}{2^{3-a}\lambda}\int_U \frac{u(x,t)^{\alpha-\lambda-1}}{(1+u(x,t))^{2\lambda}}|\nabla u(x,t)|^{2-a} W_1(x)\d x\\
   &\le \bar C_\alpha\chi_*^{\gamma_*}\left\{ \left(1 +\int_U u(x,t)^\alpha\phi(x) \d x\right)^{1+\mu_*/\alpha}
 +M_\alpha(t)\right\},
 \end{aligned}
\eeq
where 
\beq\label{Mtdef}
M_\alpha(t) =1+\int_\Gamma \left[ (\psi_1^-(x,t))^{\frac{\alpha+r}{r+\lambda}}+ (\psi_2^-(x,t))^{\frac{\alpha+r}{r}} \right]\d S.
\eeq
\end{proposition}
\begin{proof}
Throughout the proof, the number $\alpha$ satisfies $\alpha>\alpha_*$, where $\alpha_*$ will be explicitly defined in \eqref{alstar} of Step 7 below.

\medskip\noindent\emph{Step 1.} Suppose 
\beq\label{coal0}
\alpha>\lambda+1.
\eeq
Define two quantities 
$$ I_0(t)=\int_U \frac{u(x,t)^{\alpha-\lambda-1}}{(1+u(x,t))^{2\lambda}}|\nabla u(x,t)|^{2-a} W_1(x)\d x
\text{ and }
J_0(t)=\int_U u(x,t)^\alpha\phi(x) \d x.$$

For our convenience, we denote
$Q(x,t)=\nabla u(x,t)+u^{2\lambda}(x,t) \mathcal Z(x,t)$.
We will use the following estimates for $Q$.
Applying inequality \eqref{ee3} with $p=1-a<1$ yields
\beq\label{Psine1}
| Q|^{1-a}=|\nabla u+u^{2\lambda} \mathcal Z(x,t)|^{1-a}\le |\nabla u|^{1-a}+u^{2\lambda(1-a)}|\mathcal Z(x,t)|^{1-a}.
\eeq
Applying inequality \eqref{ee6} with $p=2-a>1$ yields
\beq\label{Psine}
| Q|^{2-a}=|\nabla u+u^{2\lambda} \mathcal Z(x,t)|^{2-a}\ge 2^{a-1}|\nabla u|^{2-a}-u^{2\lambda(2-a)}| \mathcal Z(x,t)|^{2-a}.
\eeq

Note that
\beq\label{ez}
\ddt\int_U \phi u^\alpha \d x
=\ddt\int_U \phi (u^\lambda)^{\alpha/\lambda} \d x
=\frac{\alpha}{\lambda}\int_U \phi(u^\lambda)_t u^{\alpha-\lambda} \d x.
\eeq
Choosing the test function $\varphi(x,t)=u^{\alpha-\lambda}$ in \eqref{weakf} gives
\beqs
\int_U \phi(u^\lambda)_t u^{\alpha-\lambda} \d x
= -(\alpha-\lambda)\int_U  X (x, u, Q)\cdot(\nabla u)  u^{\alpha-\lambda-1} \d x
-\int_\Gamma  (\psi_1 +\psi_2  u^\lambda)u^{\alpha-\lambda} \d S.
\eeqs
Utilizing \eqref{ez} for the left-hand side,  we obtain
\beq\label{dJ1}
\frac\lambda\alpha J_0' = -(\alpha-\lambda) I_1 -\int_\Gamma  (\psi_1 u^{\alpha-\lambda}+\psi_2  u^\alpha)  \d S \le -(\alpha-\lambda) I_1 +I_2,
\eeq
where
\beqs
I_1=\int_U  X (x,u,Q)\cdot(\nabla u)  u^{\alpha-\lambda-1} \d x
\text{ and }
I_2=\int_\Gamma  (\psi_1^-u^{\alpha-\lambda}+\psi_2^- u^\alpha)  \d S.
\eeqs

\medskip\noindent\emph{Step 2.}
In $I_1$, rewriting  the $\nabla u= Q -u^{2\lambda} \mathcal Z(x,t)$, we have 
\beq \label{Ionea}
I_1=\int_U X (x, u, Q)\cdot Q  u^{\alpha-\lambda-1} \d x
-\int_U  X (x,u,  Q)\cdot  \mathcal Z(x,t)\, u^{\alpha+\lambda-1} \d x
\eqdef I_{1,1}+I_{1,2}.
\eeq 

For the integrand of the term $I_{1,1}$  in \eqref{Ionea}, we use the  inequality \eqref{XX2} from Hypothesis \ref{XH} in section \ref{IBVPsec} and inequality \eqref{Psine} to have
\begin{align*}
J_1 &\eqdef   X (x, u,  Q)\cdot Q  
    \ge \frac{W_1| Q|^{2-a}}{\chi_*^2(1+u)^{2\lambda}}-W_2\\
 &\ge \frac{2^{a-1}W_1|\nabla u|^{2-a}}{{\chi_*^2(1+u)^{2\lambda}}}
 -\frac{W_1u^{2\lambda(2-a)}|\mathcal Z(x,t)|^{2-a}}{\chi_*^2(1+u)^{2\lambda}}-W_2.
\end{align*}
Simply estimating $u^{2\lambda(2-a)}/(1+u)^{2\lambda}\le u^{2\lambda(1-a)}$, we obtain 
\beq \label{J1i}
J_1\ge \frac{2^{a-1}W_1|\nabla u|^{2-a}}{{\chi_*^2(1+u)^{2\lambda}}}
 -\frac{W_1u^{2\lambda(1-a)}| \mathcal Z(x,t)|^{2-a}}{\chi_*^2}-W_2.
\eeq 

Denote $c_0=2^{a-1}\in(0,1)$.  Multiplying inequality \eqref{J1i} by $u^{\alpha-\lambda-1}$ and integrating over $U$ give 
\beq\label{Ioneone}
    I_{1,1}\ge c_0 \chi_*^{-2}I_0 - \chi_*^{-2}I_3 -I_4,
\eeq 
where
\begin{align*}
I_3=\int_U W_1| \mathcal Z(x,t)|^{2-a} u^{\alpha+\lambda(1-2a)-1}\d x
\text{ and }
I_4=\int_U W_2 u^{\alpha-\lambda-1} \d x.
\end{align*}

For the term $I_{1,2}$ in \eqref{Ionea}, applying the Cauchy--Schwarz inequality gives
\beq\label{I12}
|I_{1,2}| \le \int_U J_2 \d x,\text{ where }J_2= |X (x,u,  Q)|\cdot | \mathcal Z(x,t)|\, u^{\alpha+\lambda-1}.
\eeq
Using property   \eqref{XX1} from Hypothesis \ref{XH} in section \ref{IBVPsec} for $X$ and inequality \eqref{Psine1} for $|Q|^{1-a}$, we estimate $J_2$ by 
\begin{align*}
 J_2  &\le W_0| Q|^{1-a}\cdot |\mathcal Z(x,t)|   u^{\alpha+\lambda-1}\\
    &\le  W_0\cdot|\nabla u|^{1-a}\cdot| \mathcal Z(x,t)| u^{\alpha+\lambda-1}+ W_0\cdot u^{2\lambda(1-a)}| \mathcal Z(x,t)|^{1-a}\cdot | \mathcal Z(x,t)| u^{\alpha+\lambda-1}.
\end{align*}
Thus,
\beq \label{J2est}
J_2\le J_{2,1} + W_0 |\mathcal Z(x,t)|^{2-a}u^{\alpha+\lambda(3-2a)-1},
\text{ where } J_{2,1}=W_0|\nabla u|^{1-a} |\mathcal Z(x,t)| u^{\alpha+\lambda-1}.
\eeq 

Given a number $\varep_1>0$. Rewriting $J_{2,1}$ as
\begin{align*}
 J_{2,1} 
   &=u^{\alpha-\lambda-1} \left\{ \varepsilon_1^{\frac{1-a}{2-a}}W_1^{\frac{1-a}{2-a}}|\nabla u|^{1-a}(1+u)^{-2\lambda (1-a)/(2-a)}\right\}\\
 &\quad   \cdot\left\{ \varepsilon_1^{-\frac{1-a}{2-a}} W_1^{-\frac{1-a}{2-a}}   W_0  |\mathcal Z(x,t)| u^{2\lambda}(1+u)^{2\lambda (1-a)/(2-a)} \right\}.
\end{align*}
Then applying Young's inequality with powers $(2-a)/(1-a)$ and $(2-a)$ to the last two factors grouped by the parentheses, we have
\beqs 
 J_{2,1} 
\le u^{\alpha-\lambda-1}  \left\{ \varep_1\frac{ W_1|\nabla u|^{2-a}}{(1+u)^{2\lambda}}
 + \varepsilon_1^{-(1-a)} W_1^{a-1}W_0^{2-a}|\mathcal Z(x,t)|^{2-a}  u^{2\lambda(2-a)}(1+u)^{2\lambda(1-a)}
 \right\}.
\eeqs 
We use \eqref{ee2} to bound the last term $(1+u)^{2\lambda(1-a)}$ by 
\beqs
(1+u)^{2\lambda(1-a)}\le 2^{2\lambda(1-a)}(1+u^{2\lambda(1-a)}).
\eeqs
Thus,
\begin{align*}
 J_{2,1} 
   &\le  \varep_1 W_1|\nabla u|^{2-a}(1+u)^{-2\lambda}u^{\alpha-\lambda-1}\\
&\quad 
 + 2^{2\lambda(1-a)}\varepsilon_1^{-(1-a)} W_1^{a-1}W_0^{2-a}|\mathcal Z(x,t)|^{2-a} \cdot  \left\{ u^{\alpha+\lambda(3-2a)-1}
 + u^{\alpha+\lambda(5-4a)-1} \right\}.
\end{align*}
Combining this with \eqref{I12} and \eqref{J2est} gives
\beq\label{Ionetwo}
    |I_{1,2}|
\le \varep_1 I_0 +2^{2\lambda(1-a)}\varepsilon_1^{-(1-a)} I_5+I_6,
\eeq
where
\begin{align*}
I_5&= \int_U W_0^{2-a}W_1^{a-1} |\mathcal Z(x,t)|^{2-a}( u^{\alpha+\lambda(3-2a)-1}
 + u^{\alpha+\lambda(5-4a)-1})\d x,\\
I_6&=\int_U  W_0|\mathcal Z(x,t)|^{2-a} u^{\alpha+\lambda(3-2a)-1}
 \d x.
\end{align*}
Combining \eqref{Ionea}, \eqref{Ioneone} and \eqref{Ionetwo},  we have 
\beq\label{prediff}
\begin{aligned}
 -(\alpha-\lambda)I_1
 &\le  -(\alpha-\lambda)(c_0\chi_*^{-2}-\varep_1)I_0
     +(\alpha-\lambda) \chi_*^{-2}I_3+(\alpha-\lambda)I_4\\
     &\quad +2^{2\lambda(1-a)}(\alpha-\lambda)\varepsilon_1^{-(1-a)} I_5 +(\alpha-\lambda)I_6.
\end{aligned}
\eeq

\medskip\noindent\emph{Step 3.}
In dealing with the integral $I_6$ in \eqref{prediff}, we note, by Young's inequality with powers $2-a$ and $(2-a)/(1-a)$, that 
\beq\label{W0W1}
W_0=(W_0 W_1^{-\frac{1-a}{2-a}})\cdot W_1^{\frac{1-a}{2-a}}
\le W_0^{2-a}W_1^{a-1}+W_1=W_3.
\eeq
Observe that, in formulas of $I_3$, $I_5$, $I_6$, the powers of $u$ have the  minimum  $\alpha-\lambda-1$ and the maximum $\alpha+\lambda(5-4a)-1$.
Then, by inequality \eqref{ee4},
\beq\label{III}
I_3\le I_7,\ 
I_5\le 2 I_7,\text{ and with the aid of \eqref{W0W1}, } 
I_6\le I_7,
\eeq
where 
\begin{align*}
I_7
&=\int_U (W_0^{2-a}W_1^{a-1}+W_1)|\mathcal Z(x,t)|^{2-a}\left (u^{\alpha-\lambda-1} +u^{\alpha+\lambda(5-4a)-1}\right)\d x\\
&=\int_U W_3|\mathcal Z(x,t)|^{2-a}\left (u^{\alpha-\lambda-1} +u^{\alpha+\lambda(5-4a)-1}\right)\d x.
\end{align*}
Combining the estimates in \eqref{III} with \eqref{prediff} gives
\beq\label{JoIs}
 -(\alpha-\lambda)I_1\le -(\alpha-\lambda)(c_0\chi_*^{-2}-\varep_1)I_0
     + (\alpha-\lambda)I_4+ (\alpha-\lambda)C_0 I_7,
\eeq
where
\beq\label{C0def}
C_0=\chi_*^{-2}+2^{2\lambda(1-a)+1}\varep_1^{-(1-a)}+1.
\eeq
Using \eqref{Zest} to estimate $\mathcal Z(x,t)$ in the definition of $I_7$, we obtain 
\beqs
I_7\le c_{\mathcal Z}^{2-a}\chi_*^{2(2-a)} \left(\int_U W_3u^{\alpha-\lambda-1} \d x+\int_U W_3u^{\alpha+\lambda(5-4a)-1}\d x\right).
\eeqs
Together with the formula of $I_4$ and the fact $\chi_*\ge 1$, this implies
\beq\label{I47}
I_4+C_0 I_7\le C_1\left(\int_U (W_2+W_3)u^{\alpha-\lambda-1} \d x+\int_U W_3u^{\alpha+\lambda(5-4a)-1}\d x\right)
=C_1(I_8+I_9),
\eeq
where 
\begin{align}\label{C1def}
C_1&=\max\{1,C_0c_{\mathcal Z}^{2-a}\}\chi_*^{2(2-a)},\\
I_8&=\int_U W_4 u^{\alpha-\lambda-1} \d x,\quad 
I_9=\int_U W_3u^{\alpha+\lambda(5-4a)-1}\d x. \notag  
\end{align}
Thus, it follows from \eqref{JoIs} and \eqref{I47} that 
\beq\label{negI1}
 -(\alpha-\lambda)I_1\le -(\alpha-\lambda)(c_0\chi_*^{-2}-\varep_1)I_0
 + C_1(\alpha-\lambda)(I_8+I_9).
\eeq
Combining \eqref{dJ1} with \eqref{negI1} yields
\beq\label{J0dif0}
 \frac\lambda\alpha J_0' +(\alpha-\lambda)(c_0\chi_*^{-2}-\varep_1)I_0
 \le C_1(\alpha-\lambda)(I_8+I_9)   +I_2.
\eeq

\medskip\noindent\emph{Step 4.}
We estimate $I_8, I_9$ and $I_2$ in \eqref{J0dif0}. 

\noindent\emph{Estimates of $I_8$.} By rewriting $I_8$ and, thanks to \eqref{coal0},  applying Young's inequality for the powers $\alpha/(\lambda+1)$ and $\alpha/(\alpha-\lambda-1)$, we have
\beq\label{I48}
    I_8=\int_U (W_4 \phi^{-\frac{\alpha-\lambda-1}{\alpha}})\cdot  (\phi^{\frac{\alpha-\lambda-1}{\alpha}}u^{\alpha-\lambda-1}) \d x\\
    \le {\mathcal K}_5+J_0.
\eeq

\noindent\emph{Estimate of $I_9$.} 
Note from \eqref{newrr} and \eqref{coal0} that  
\beqs 
r-\lambda(5-4a)+1>0 \text{ and }\alpha+r>\alpha+\lambda(5-4a)-1>\alpha-1>0.
\eeqs
Then applying Young's inequality with the powers $\frac{\alpha+r}{r-\lambda(5-4a)+1}$  and $\frac{\alpha+r}{\alpha+\lambda(5-4a)-1}$, we have
\beq\label{b2esti}
 I_9 \le \int_U W_3^{\frac{\alpha+r}{r-\lambda(5-4a)+1}} \d x+ \int_U u^{\alpha+r}\d x
={\mathcal K}_6 + \int_U u^{\alpha+r}\d x.
\eeq
Let us estimate the last integral in \eqref{b2esti}. Let $\varep_2>0$.
We apply  Lemma \ref{RoS1} to
\beq \label{apply}
p, s, \beta \text{ in \eqref{psbe}, } \varep=\varep_2,\   \omega(x)=1,\ \varphi(x)=\phi(x),\  W_*(x)=W_1(x).
\eeq 
Conditions \eqref{manya} and \eqref{abeta} from Lemma~\ref{RoS1} become
\beq \label{coal1}
\alpha\ge 1+\lambda,\  \alpha>\frac{1-a-\lambda}{1-a},\ 
\alpha>\frac{2(r-1+\lambda+a)}{r_*}\text{ and }
\alpha>\frac{2\lambda r_1}{1-r_1}.
\eeq 
It follows from inequality \eqref{RSi2} that
\beq\label{ualphar}
\int_U u^{\alpha+r}  \d x 
 \le \varep_2 I_0 +{\widehat\Phi}_1
\left( 1 + J_0\right)^{1+\widehat\mu_3/\alpha},
\eeq 
where $\widehat\mu_3$ and ${\widehat\Phi}_1$ are given by formulas \eqref{hatmu1} and \eqref{phi00}, respectively,  applied to \eqref{apply}.
Then we have 
\beq\label{CI9}
 I_9
 \le {\mathcal K}_6+\varep_2 I_0 +{\widehat\Phi}_1 \left( 1 + J_0\right)^{1+\widehat\mu_3/\alpha}.
\eeq

\noindent\emph{Estimate of $I_2$.}
Applying  Young's inequality with the powers $\frac{\alpha+r}{r+\lambda }$  and $\frac{\alpha+r}{\alpha-\lambda}$  to the product $\psi_1^-u^{\alpha-\lambda}$, and the powers $\frac{\alpha+r}r$ and $ \frac{\alpha+r}\alpha$ to the product $\psi_2^- u^\alpha$ gives
\beq\label{I2ineq}
\begin{aligned}
I_2
&\le \left(\int_\Gamma (\psi_1^-)^{\frac{\alpha+r}{r+\lambda}} \d S
+\int_\Gamma u^{\alpha+r}\d S\right)
+\left(\int_\Gamma (\psi_2^-)^{\frac{\alpha+r}{r}} \d S+\int_\Gamma u^{\alpha+r}\d S\right)\\
&= 2 \int_\Gamma u^{\alpha+r}\d S+K,
\end{aligned}
\eeq
 where 
\beqs 
K=K(t)\eqdef \int_\Gamma \left[ (\psi_1^-(x,t))^{\frac{\alpha+r}{r+\lambda}}+ (\psi_2^-(x,t))^{\frac{\alpha+r}{r}} \right] \d S.
\eeqs 
Let $\varep_3>0$. We estimate the surface integral of $u^{\alpha+r}$ over $\Gamma$, with $r>0$ from \eqref{newrr},  by applying Lemma~\ref{Rtrace}
to
\beq \label{traceapp}
p, s, \beta \text{ in \eqref{psbe}, }
\varep=\varep_3,\    \varphi(x)=\phi(x),\    W_*(x)=W_1(x).
\eeq 
The numbers $\widetilde\theta$, $\widetilde\mu_1$, $\widehat\mu_5$ and ${\widehat\Phi}_2$ in \eqref{tmtilde}, \eqref{muh45} and \eqref{phi20} are now computed using \eqref{traceapp}.
Because of the assumption $r>\frac{1-a-\lambda}{2-a}$ from \eqref{newrr}, we have $\widetilde{r}$ in \eqref{newtr} is positive.  The conditions \eqref{manyb},  \eqref{only}  from Lemma~\ref{Rtrace} become 
\beq \label{coal2}
\alpha\ge 1+\lambda,\  \alpha>\frac{1-a-\lambda}{1-a}, 
\ 
\alpha>\frac{2(2-a)(r-1+\lambda+a)}{r_*(1-a)},\ 
\eeq 
while \eqref{bcond}, using the equivalent form \eqref{abeuqiv}, becomes 
\beq \label{coal3}
\alpha>\frac{2\lambda}{(1-a)(1-\widetilde\theta)}-\widetilde\mu_1.
\eeq 
Then it follows from inequalities \eqref{tru3} and \eqref{I2ineq} that
\beq\label{I2ineqq}
I_2 \le 6\varepsilon_3 I_0  + 2{\widehat\Phi}_2 \left( 1 + J_0\right)^{1+\widehat\mu_5/\alpha}+K.
\eeq

\medskip\noindent\emph{Step 5.}
Combining \eqref{I48}, \eqref{CI9} and \eqref{I2ineqq} with \eqref{J0dif0}, 
and noticing that $\alpha-\lambda\ge 1$, $C_1\ge 1$, we now have
\begin{align*}
   & \frac\lambda\alpha J_0' +\left\{(\alpha-\lambda) (c_0\chi_*^{-2}-\varep_1-C_1\varep_2)-6\varep_3\right\}I_0\\
   & \le   
 (\alpha-\lambda)\left[ C_1J_0+C_1\widehat\Phi_1(1+ J_0)^{1+\widehat\mu_3/\alpha}+ 2\widehat\Phi_2(1+ J_0)^{1+\widehat\mu_5/\alpha}\right]+C_1(\alpha-\lambda)({\mathcal K}_5+{\mathcal K}_6)+K\\
 &\le  C_2(\alpha-\lambda)(1+ J_0)^{1+\max\{\widehat\mu_3,\widehat\mu_5\}/\alpha}+C_1(\alpha-\lambda)({\mathcal K}_5+{\mathcal K}_6)+K,
\end{align*}
where
\beq\label{Hmu}
C_2=C_1(1+\widehat\Phi_1)+2\widehat\Phi_2.
\eeq
Comparing $\widehat\mu_3$ in \eqref{hatmu1} with $\widehat\mu_5$ in \eqref{muh45} gives $\widehat\mu_3\le \widehat\mu_5$.
Then by taking 
\beq \label{eee}
\varep_1=c_0\chi_*^{-2}/4,\ 
\varep_2=c_0\chi_*^{-2}/(4C_1)\text{ and }
\varep_3=c_0\chi_*^{-2}(\alpha-\lambda)/24,
\eeq 
we obtain 
\beq\label{J0dif2}
   \frac\lambda\alpha J_0' +\frac{c_0(\alpha-\lambda)\chi_*^{-2}}{4}I_0
   \le  C_2(\alpha-\lambda)(1+ J_0)^{1+\widehat\mu_5/\alpha}+C_1(\alpha-\lambda)({\mathcal K}_5+{\mathcal K}_6)+K.
\eeq

\medskip\noindent\emph{Step 6.} We calculate the powers and quantities appearing in \eqref{J0dif2}.

For inequality \eqref{RSi2} in Lemma \ref{RoS1}, which was used with the choice \eqref{apply} to obtain \eqref{ualphar} above, the number $m$ in \eqref{powdef} becomes 
$m=\frac{\alpha+1-\lambda-a}{2-a}$, and numbers $\theta$, $\mu_1$ there become those in \eqref{trtm0}, \eqref{trtm2}, while the numbers $\widehat\mu_1$  in \eqref{hatmu1} and $\widehat\beta_1$ in \eqref{phi00} become
\beq \label{mhmhb}
 \widehat\mu_1=\mu_1+\frac{2\lambda\theta}{1-\theta}\text{ and }   \widehat\beta_1=\frac{2\lambda\theta}{1-\theta}(1+\frac1\alpha).
\eeq 
From  \eqref{Gdef},  \eqref{Calpha} and with the choice \eqref{apply},
\begin{align*}
{\mathcal E}_1&= 1+\int_U \phi(x) \d x={\mathcal K}_0,\quad  
G_1 = \left(\int_U \phi(x)^{-\frac {\alpha+1-\lambda-a}{\alpha(1-a)-1+\lambda+a}} \d x\right)^\frac {\alpha(1-a)-1+\lambda+a}{\alpha}={\mathcal K}_2, \\ 
G_3 &= \left(\int_U \phi(x)^{-1} \d x\right)^{1+\mu_1/\alpha} ={\mathcal K}_1,  \quad 
\Phi_1 =G_1^{\theta}G_3^{1-\theta}=\mathcal K_2^{\theta}\mathcal K_1^{1-\theta},\\ 
{\widetilde G}_2&=\left(\int_U \phi(x)^{-\frac{2\lambda r_1}{\alpha(1-r_1)-2\lambda r_1}}W_1(x)^{-\frac{r_1\alpha}{\alpha(1-r_1)-2\lambda r_1}}\d x\right)^{\frac{\alpha(1-r_1)-2\lambda r_1}{\alpha r_1}}={\mathcal K}_3.
\end{align*}
Therefore, 
${\widetilde\Phi}_2$ from \eqref{Calpha} becomes
${\widetilde\Phi}_2={\mathcal K}_3^\frac\theta{1-\theta} \mathcal K_1 {\mathcal K}_0^{\frac{2\lambda\theta}{\alpha(1-\theta)}}$, 
and
${\widehat\Phi}_1 $ from \eqref{phi00} becomes
\beq\label{Phat1}
{\widehat\Phi}_1 
=D_{1,m,\theta}\mathcal K_2^{\theta}\mathcal K_1^{1-\theta}
+2^{1+\widehat\beta_1} D_{2,m,\theta}\varep_2^{-\frac\theta{1-\theta}}{\mathcal K}_3^\frac\theta{1-\theta}\mathcal K_1 {\mathcal K}_0^{\frac{2\lambda\theta}{\alpha(1-\theta)}}.
\eeq 

For Lemma~\ref{Rtrace}, which was used with the choice \eqref{traceapp} to obtain \eqref{I2ineqq} above,  we compute $\widetilde r$ in \eqref{rtilde} by \eqref{newtr}, $\widetilde \theta$ and $\widetilde \mu_1$ in \eqref{tmtilde} by \eqref{trtm0} and \eqref{trtm1}, number $\beta_*$ in \eqref{bcond} by \eqref{trtm2}, 
powers $ \widehat\beta_2$, $\widehat\beta_3$ and $\widehat\mu_5$ in \eqref{muh45} by
\beq \label{muh5}
\widehat\beta_2=\frac{2\lambda}{1-a}\left(1+\frac1\alpha\right),
   \widehat\beta_3=\frac{2\lambda}{1-\widetilde\theta}\left(\frac1{1-a}+\widetilde\theta\right),
\widehat\mu_5=\max\left \{r, \widehat\mu_1, \widetilde r+\frac{2\lambda}{1-a}, \widetilde\mu_1+\widehat\beta_3\right\}.
\eeq 
Then, one has from \eqref{Ph34}
\beqs
G_4 = \left(\int_U \phi(x)^{-1} \d x\right)^{1+\mu_1/\alpha} ={\mathcal K}_1,  \quad
\Phi_3=G_1^{\theta}G_4^{1-\theta}=\mathcal K_2^\theta\mathcal K_1^{1-\theta}, 
\eeqs 
and from \eqref{tfifour}
\begin{align*}
{\widetilde\Phi}_4 &={\mathcal E}_1^{\frac{2\lambda\theta}{\alpha(1-\theta)}}{\widetilde G}_2^\frac\theta{1-\theta} G_4={\mathcal K}_0^{\frac{2\lambda\theta}{\alpha(1-\theta)}}{\mathcal K}_3^\frac\theta{1-\theta} \mathcal K_1,\\
{\widetilde G}_5&=\left(\int_U \phi(x)^{-\frac{\alpha+\beta_*}{\alpha-\beta_*}} W_1(x)^{-\frac{\beta_*}{2\lambda}\cdot \frac{\alpha}{\alpha-\beta_*}}\d x\right)^{\frac{\alpha-\beta_*}{\alpha}(1+\widetilde\mu_1/\alpha)}={\mathcal K}_4,\\
{\widetilde\Phi}_6&={\mathcal E}_1^{\frac{2\lambda}{\alpha(1-a)}}G_1^{\widetilde\theta} {\widetilde G}_5^{1-\widetilde\theta}={\mathcal K}_0^{\frac{2\lambda}{\alpha(1-a)}}\mathcal K_2^{\widetilde\theta} {\mathcal K}_4^{1-\widetilde\theta},\quad
{\widetilde\Phi}_7
 = {\mathcal E}_1^{\frac{\widehat\beta_3}{\alpha}} {\widetilde G}_2^{\frac{\widetilde\theta}{1-\widetilde\theta}}\widetilde G_5={\mathcal K}_0^{\frac{\widehat\beta_3}{\alpha}} {\mathcal K}_3^{\frac{\widetilde\theta}{1-\widetilde\theta}}{\mathcal K}_4.
 \end{align*}
Let $z_1,\ldots,z_5$, $\widetilde z_2$, $\widetilde z_4$, $\widetilde z_5$,   be the numbers in \eqref{zzz}, \eqref{trace1} and \eqref{ztils} with the choice \eqref{psbe}.
 Thus, from \eqref{phi20},
\beq\label{Phat2}
\begin{aligned}
{\widehat\Phi}_2&=  z_1 \Phi_3+ \varep_3^{-\frac\theta{1-\theta}}\widetilde z_2  
{\widetilde\Phi}_4 
+\varep_3^{-\frac 1{1-a}}\widetilde z_4 \widetilde\Phi_6 + \varep_3^{-(\frac 1 {1-a}+\frac {2-a}{1-a}\cdot \frac{\widetilde\theta}{1-\widetilde\theta} )}\widetilde z_5 \widetilde\Phi_7\\
&=z_1\mathcal K_2^\theta\mathcal K_1^{1-\theta}
+\varep_3^{-\frac\theta{1-\theta}}\widetilde  z_2{\mathcal K}_0^{\frac{2\lambda\theta}{\alpha(1-\theta)}}{\mathcal K}_3^{\frac\theta{1-\theta}} \mathcal K_1\\
&\quad +\varep_3^{-\frac 1{1-a}}\widetilde z_4 {\mathcal K}_0^{\frac{2\lambda}{\alpha(1-a)}}\mathcal K_2^{\widetilde\theta} {\mathcal K}_4^{1-\widetilde\theta}
+\varep_3^{-(\frac 1 {1-a}+\frac {2-a}{1-a}\cdot \frac{\widetilde\theta}{1-\widetilde\theta} )}\widetilde z_5 {\mathcal K}_0^{\frac{\widehat\beta_3}{\alpha}} {\mathcal K}_3^{\frac{\widetilde\theta}{1-\widetilde\theta}}{\mathcal K}_4.
\end{aligned}
\eeq 

\medskip\noindent\emph{Step 7.} We verify the conditions for $\alpha$ and specify the number $\alpha_*$ mentioned at the beginning of the proof. 
The requirements are \eqref{coal0}, \eqref{coal1}, \eqref{coal2} and \eqref{coal3}.

It is clear that under assumption \eqref{coal0}, the first two inequalities in each of \eqref{coal1} and \eqref{coal2} are satisfied. Hence, the  conditions \eqref{coal0}, \eqref{coal1}, \eqref{coal2}  are equivalent to 
\beq\label{newaa2}
\alpha> \max\left\{ \lambda+1,\frac{2(r-1+\lambda+a)}{r_*},\frac{2\lambda r_1}{1-r_1}, 
\frac{2(2-a)(r+a+\lambda-1)}{r_*(1-a)}\right\}.
\eeq
Consider the numbers in the list on the right-hand side of \eqref{newaa2}.
If $r+a+\lambda-1\le 0$, then we ignore both the second and fourth numbers for the condition \eqref{newaa2}. Otherwise, the fourth number is larger than the second one. Thus, 
we can ignore the second number and, hence, \eqref{newaa2} is equivalent to 
\beq\label{newaa3}
\alpha> \max\left\{ \lambda+1,\frac{2\lambda r_1}{1-r_1}, 
\frac{2(2-a)(r+a+\lambda-1)}{r_*(1-a)}\right\}.
\eeq

We now examine condition \eqref{coal3}. The last number in \eqref{coal3} can be computed and estimated by 
\beq \label{lmm}
 \frac{2\lambda}{(1-a)(1-\widetilde\theta)}-\widetilde\mu_1
 =\frac{2\lambda-(1-a)[\widetilde r+\widetilde \theta(a+\lambda-1)]}{(1-a)(1-\widetilde\theta)}
 <\frac{2\lambda+(1-a)\widetilde \theta }{(1-a)(1-\widetilde\theta)}.
\eeq 
To estimate the last number further, we find an upper bound for the number $\widetilde\theta$ first.
Set 
\beqs 
\theta_{**}=\frac{1}{1+r_*/2}.
\eeqs 
\noindent\textbf{Claim.}\textit{  One has $\widetilde\theta\le \theta_{**}$ provided
\beq\label{biga2} 
\alpha\ge 2\widetilde r\left(1+\frac2{r_*}\right)+\frac{4\lambda}{r_*}.
\eeq
}
\noindent\textit{Proof of the Claim.} Because of $\widetilde r > 0$, thanks to \eqref{newtr},  and \eqref{biga2}, one has 
\beq\label{biga1} 
\alpha>\frac{4\lambda}{r_*}> \frac{2\lambda}{1+r_*}.
\eeq
Then
\beqs
\widetilde\theta=\frac{\alpha+2\widetilde r}{\alpha(1+r_*)+2(1-a-\lambda)}
<\frac{\alpha+2\widetilde r}{\alpha(1+r_*)-2\lambda}
\le \frac{1}{1+r_*/2}=\theta_{**},
\eeqs 
where the last inequality comes from \eqref{biga1} and \eqref{biga2}.
Thus the Claim is true.

Under such a condition \eqref{biga2}, we have, thanks to the Claim, $\widetilde\theta\le \theta_{**}$, and, together with \eqref{lmm}, it implies 
\beqs
\frac{2\lambda}{(1-a)(1-\widetilde\theta)} -\widetilde\mu_1
< \frac{2\lambda+(1-a)\theta_{**} }{(1-a)(1-\theta_{**})}
=\alpha_{**}\eqdef \frac{2\lambda(2+r_*)}{(1-a)r_*}+\frac{2}{r_*}.
\eeqs 
For the purpose of later estimates, we also bound $\theta$ in a similar way here.
By replacing $\widetilde r$ with $r$ in  the Claim, one has    
\beqs
\theta=\frac{\alpha+2r}{\alpha(1+r_*)+2(1-a-\lambda)}
<\frac{\alpha+2 r}{\alpha(1+r_*)-2\lambda}\le \frac{1}{1+r_*/2}=\theta_{**},
\eeqs
provided
\beqs
\alpha\ge 2r\left(1+\frac2{r_*}\right)+\frac{4\lambda}{r_*}.
\eeqs
Therefore, a sufficient condition for \eqref{newaa3}, \eqref{coal3} and, additionally, 
\beq \label{tte}
\theta\le \theta_{**}\text{ and }\widetilde\theta\le \theta_{**}
\eeq 
is 
\beq\label{newaa1}
\begin{aligned}
\alpha>\alpha_*, \text{ where }\alpha_*=  \max & \Big\{  \lambda+1,\frac{2\lambda r_1}{1-r_1}, 
\frac{2(2-a)(r+a+\lambda-1)}{r_*(1-a)},\alpha_{**},\\
&\quad 
2r\left(1+\frac2{r_*}\right)+\frac{4\lambda}{r_*},
2\widetilde r\left(1+\frac2{r_*}\right)+\frac{4\lambda}{r_*}\Big\}.
\end{aligned}
\eeq
Taking into account the fact $a,r_*\in(0,1)$, one has 
$$\alpha_{**}=\frac{2\lambda(2+r_*)}{(1-a)r_*}+\frac{2}{r_*}>4\lambda+2>\lambda+1.$$ Therefore, we obtain
\beq\label{alstar}
\begin{aligned}
\alpha_*= \max \Big\{ & \frac{2\lambda r_1}{1-r_1}, 
\frac{2(2-a)(r+a+\lambda-1)}{r_*(1-a)},\frac{2\lambda(2+r_*)}{(1-a)r_*}+\frac{2}{r_*},\\
&\quad 
2r\left(1+\frac2{r_*}\right)+\frac{4\lambda}{r_*},
2\widetilde r\left(1+\frac2{r_*}\right)+\frac{4\lambda}{r_*}\Big\}.
\end{aligned}
\eeq 
With $\alpha_*$ defined in \eqref{alstar}, the requirement $\alpha>\alpha_*$ from \eqref{newaa1} was assumed from the beginning of the proof which enabled all the calculations and estimates in the previous steps.

\medskip\noindent\emph{Step 8.} 
We estimate $\widehat\mu_5$ , $C_1$ and $C_2$ appearing in \eqref{J0dif2}. Below, $C$ denotes a generic constant depending on $\alpha$ but independent of $\chi_*$ and $K_j$; its value may change from one place to another.

\noindent\emph{Estimate of $\widehat\mu_5$.} 
Note from \eqref{mhmhb} and \eqref{trtm1}  that 
\beqs
\widehat\mu_1=\frac{r+\theta(a+\lambda-1)}{1-\theta}+\frac{2\lambda\theta}{1-\theta}
\le \frac{r+\theta(a+\lambda)+2\lambda\theta}{1-\theta}=\frac{r+\theta(a+3\lambda)}{1-\theta}.
\eeqs
Similarly, it follows from \eqref{trtm1} and \eqref{muh5} that 
\beqs
 \widetilde\mu_1+\widehat\beta_3 =\frac{\widetilde r+\widetilde\theta(a+\lambda-1)}{1-\widetilde\theta}
 +\frac{2\lambda \widetilde\theta}{1-\widetilde\theta}+\frac{2\lambda}{(1-\widetilde\theta)(1-a)}
 \le \frac{\widetilde r+\widetilde\theta(a+3\lambda)}{1-\widetilde\theta}+\frac{2\lambda}{(1-\widetilde\theta)(1-a)}.
\eeqs 
Moreover, 
\beqs
r\le \frac{r+\theta(a+3\lambda)}{1-\theta}\text{ and }
\widetilde r+\frac{2\lambda}{1-a}\le \frac{\widetilde r+\widetilde\theta(a+3\lambda)}{1-\widetilde\theta}+\frac{2\lambda}{(1-\widetilde\theta)(1-a)}.
\eeqs
Combining these estimates with the formula of $\widehat\mu_5$ in \eqref{muh5} and then utilizing \eqref{tte} yield 
\beq\label{mums}
\begin{split}
\widehat\mu_5&\le \max\left \{\frac{r+\theta(a+3\lambda)}{1-\theta}, \frac{\widetilde r+\widetilde\theta(a+3\lambda)}{1-\widetilde\theta}+\frac{2\lambda}{(1-\widetilde\theta)(1-a)}\right\}\\
&\le \mu_*\eqdef \max\left \{\frac{r+\theta_{**}(a+3\lambda)}{1-\theta_{**}}, \frac{\widetilde r+\theta_{**}(a+3\lambda)}{1-\theta_{**}}+\frac{2\lambda}{(1-\theta_{**})(1-a)}\right\}.
\end{split}
\eeq

\noindent\emph{Estimate of $C_1$.} 
By the definitions of $C_1$ in \eqref{C1def}, $C_0$ in \eqref{C0def} and $\varep_1$ in \eqref{eee}, we have
\beq \label{C1ineq}
C_1=\max\left\{1,[\chi_*^{-2}+2^{2\lambda(1-a)+1}(4c_0^{-1}\chi_*^2)^{1-a}+1]c_{\mathcal Z}^{2-a}\right\}\chi_*^{2(2-a)}
\le C \chi_*^{2(1-a)+2(2-a)}=C \chi_*^{2(3-2a)} .
\eeq 

\noindent\emph{Estimate of $C_2$.} We start with the numbers ${\widehat\Phi}_1$ and ${\widehat\Phi}_2$ in the formula \eqref{Hmu} of $C_2$.
Observe, in their definitions \eqref{Phat1} and \eqref{Phat2},  that all the numbers  $D_{1,m,\theta}$, $2^{1+\widehat\beta_1}$,  $D_{2,m,\theta}$, $z_1$, $\widetilde z_2$, $\widetilde z_4$, $\widetilde z_5$ depend on $\alpha$, but are independent of $\chi_*$ and the $\mathcal K_j$. Hence, we can estimate, from \eqref{Phat1} and \eqref{Phat2},
\beq \label{PPest}
{\widehat\Phi}_1\le C {\widehat K}_1\text{ and } {\widehat\Phi}_2\le C  {\widehat K}_2,
\eeq 
where
\begin{align*}
{\widehat K}_1&=\mathcal K_2^{\theta}\mathcal K_1^{1-\theta}
+ \varep_2^{-\frac\theta{1-\theta}}{\mathcal K}_3^\frac\theta{1-\theta}\mathcal K_1 {\mathcal K}_0^{\frac{2\lambda\theta}{\alpha(1-\theta)}},\\
{\widehat K}_2&=\mathcal K_2^\theta\mathcal K_1^{1-\theta}
+\varep_3^{-\frac\theta{1-\theta}}{\mathcal K}_0^{\frac{2\lambda\theta}{\alpha(1-\theta)}}{\mathcal K}_3^{\frac\theta{1-\theta}} \mathcal K_1
 +\varep_3^{-\frac 1{1-a}}{\mathcal K}_0^{\frac{2\lambda}{\alpha(1-a)}}\mathcal K_2^{\widetilde\theta} {\mathcal K}_4^{1-\widetilde\theta}
  +\varep_3^{-(\frac 1 {1-a}+\frac {2-a}{1-a}\cdot \frac{\widetilde\theta}{1-\widetilde\theta} )} {\mathcal K}_0^{\frac{\widehat\beta_3}{\alpha}} {\mathcal K}_3^{\frac{\widetilde\theta}{1-\widetilde\theta}}{\mathcal K}_4.
\end{align*}
With the choice of $\varep_2$ and $\varep_3$ in \eqref{eee} and the estimate \eqref{C1ineq} of $C_1$, together with $\alpha-\lambda\ge 1$, one has 
\beqs
\varep_2^{-1}=4C_1 c_0^{-1}\chi_*^2\le C\chi_*^{2(3-2a)+2}=C\chi_*^{4(2-a)},\  
\varep_3^{-1}=24c_0^{-1}\chi_*^2 (\alpha-\lambda)^{-1}\le C \chi_*^2.
\eeqs
Therefore, we can estimate
\begin{align*}
{\widehat K}_1
&\le C\left (\mathcal K_2^{\theta}\mathcal K_1^{1-\theta}+\chi_*^{\frac{4(2-a)\theta}{1-\theta}}{\mathcal K}_3^\frac\theta{1-\theta}\mathcal K_1 {\mathcal K}_0^{\frac{2\lambda\theta}{\alpha(1-\theta)}}\right)\le C\chi_*^{\frac{4(2-a)\theta}{1-\theta}}(\mathcal K_2^{\theta}\mathcal K_1^{1-\theta}+{\mathcal K}_3^\frac\theta{1-\theta}\mathcal K_1 {\mathcal K}_0^{\frac{2\lambda\theta}{\alpha(1-\theta)}}),\\
{\widehat K}_2
&\le C\left (
\mathcal K_2^\theta\mathcal K_1^{1-\theta}
+\chi_*^{\frac{2\theta}{1-\theta}}{\mathcal K}_0^{\frac{2\lambda\theta}{\alpha(1-\theta)}}{\mathcal K}_3^{\frac\theta{1-\theta}} \mathcal K_1 
+\chi_*^{\frac 2{1-a}}{\mathcal K}_0^{\frac{2\lambda}{\alpha(1-a)}}\mathcal K_2^{\widetilde\theta} {\mathcal K}_4^{1-\widetilde\theta}
+\chi_*^{\frac{2}{1-a}+\frac{2(2-a)}{1-a}\cdot \frac{\widetilde\theta}{1-\widetilde\theta} } {\mathcal K}_0^{\frac{\widehat\beta_3}{\alpha}} {\mathcal K}_3^{\frac{\widetilde\theta}{1-\widetilde\theta}}{\mathcal K}_4\right)\\
&\le C \chi_*^{\max\{\frac{2\theta}{1-\theta},\frac{2}{1-a}+\frac{2(2-a)}{1-a}\cdot \frac{\widetilde\theta}{1-\widetilde\theta} \}}
\left (
\mathcal K_2^\theta\mathcal K_1^{1-\theta}
+{\mathcal K}_0^{\frac{2\lambda\theta}{\alpha(1-\theta)}}{\mathcal K}_3^{\frac\theta{1-\theta}} \mathcal K_1 
+{\mathcal K}_0^{\frac{2\lambda}{\alpha(1-a)}}\mathcal K_2^{\widetilde\theta} {\mathcal K}_4^{1-\widetilde\theta}
+ {\mathcal K}_0^{\frac{\widehat\beta_3}{\alpha}} {\mathcal K}_3^{\frac{\widetilde\theta}{1-\widetilde\theta}}{\mathcal K}_4
\right).
\end{align*}
Above, we used the fact $\chi_*\ge 1$ from \eqref{chidef}.
From \eqref{Hmu}, \eqref{C1ineq}, \eqref{PPest} and the last two estimates of ${\widehat K}_1$ and ${\widehat K}_2$, it follows that
\beq \label{Hest}
C_2
=C_1(1+\widehat\Phi_1)+2\widehat\Phi_2\le C\chi_*^{\gamma_1}C_3,
\eeq  
where
\beqs 
C_3=1+ \mathcal K_2^{\theta}\mathcal K_1^{1-\theta}+{\mathcal K}_3^\frac\theta{1-\theta}\mathcal K_1 {\mathcal K}_0^{\frac{2\lambda\theta}{\alpha(1-\theta)}}
+{\mathcal K}_0^{\frac{2\lambda}{\alpha(1-a)}}\mathcal K_2^{\widetilde\theta} {\mathcal K}_4^{1-\widetilde\theta}
+{\mathcal K}_0^{\frac{\widehat\beta_3}{\alpha}} {\mathcal K}_3^{\frac{\widetilde\theta}{1-\widetilde\theta}}{\mathcal K}_4\ge 1,
\eeqs 
\beqs 
\gamma_1=\max\left\{2(3-2a)+\frac{4(2-a)\theta}{1-\theta},  \frac{2\theta}{1-\theta},\frac{2}{1-a}+\frac{2(2-a)}{1-a}\cdot \frac{\widetilde\theta}{1-\widetilde\theta} \right\}>0.
\eeqs 

\medskip\noindent\emph{Step 9.} 
Combining \eqref{mums}, \eqref{C1ineq}, \eqref{Hest}  with \eqref{J0dif2} yields
\beqs
   \frac\lambda\alpha J_0' +\frac{c_0(\alpha-\lambda)\chi_*^{-2}}{4}I_0
   \le  C C_3\chi_*^{\gamma_1}(1+ J_0)^{1+\mu_*/\alpha}+C\chi_*^{2(3-2a)}({\mathcal K}_5+{\mathcal K}_6)+K.
\eeqs
Note that the last power of $\chi_*$ is $2(3-2a)\le \gamma_1$.
Using the fact  
\beqs 
2(3-2a)+\frac{4(2-a)\theta}{1-\theta}>
\frac{4(2-a)\theta}{1-\theta}>\frac{2\theta}{1-\theta},
\eeqs
we can drop the number $\frac{2\theta}{1-\theta}$ in the definition of $\gamma_1$.
Then  using the bounds of $\theta$ and $\widetilde\theta$ in \eqref{tte}, we have 
\begin{align*}
\gamma_1&=\max\left\{2(3-2a)+\frac{4(2-a)\theta}{1-\theta},\frac{2}{1-a}+\frac{2(2-a)}{1-a}\cdot \frac{\widetilde\theta}{1-\widetilde\theta} \right\}\\
&\le \gamma_*\eqdef \max\left\{2(3-2a)+\frac{4(2-a)\theta_{**}}{1-\theta_{**}},\frac{2}{1-a}+\frac{2(2-a)}{1-a}\cdot \frac{\theta_{**}}{1-\theta_{**}} \right\}.
\end{align*}
Thus,  we obtain 
\beq\label{J0dif3}
   J_0' +\frac{\alpha(\alpha-\lambda)\chi_*^{-2}}{2^{3-a}\lambda}I_0\le C C_4\chi_*^{\gamma_*}\left\{ (1 +J_0)^{1+\mu_*/\alpha} + 1+K\right\},
\eeq
where $C_4=C_3+{\mathcal K}_5+{\mathcal K}_6$.
Then the desired inequality \eqref{J0dbest} follows from \eqref{J0dif3} with the constant $\bar C_\alpha=C C_4$.
\end{proof}

\begin{remark}\label{3weight}
Below are alternative estimates for $I_9$ instead of \eqref{CI9}.
If $\lambda(5-4a)-1< 0$, then  by H\"older's inequality for the powers $-\alpha/(\lambda(5-4a)-1)$ and $\alpha/(\alpha+\lambda(5-4a)-1)$, one has
\beq\label{HI9}
\begin{split}
    I_9&=\int_U \left(W_3 \phi^{-\frac{\alpha+\lambda(5-4a)-1}{\alpha}}\right)\cdot  \left(\phi^{\frac{\alpha+\lambda(5-4a)-1}{\alpha}}u^{\alpha+\lambda(5-4a)-1}\right) \d x\\
    &\le \int_U W_3^{\frac{\alpha}{-\lambda(5-4a)+1}}\phi^{-\frac{\alpha+\lambda(5-4a)-1}{-\lambda(5-4a)+1}}\d x +J_0.
    \end{split}
\eeq 
If $\lambda(5-4a)-1\ge 0$, then applying inequality \eqref{RSi2} in Lemma \ref{RoS1} to
\beq \label{app}
p,s,\beta \text{ in \eqref{psbe}, } 
r= \lambda(5-4a)-1,\ 
\varep=\varep_2,\
\omega(x)=W_3(x),\  \varphi(x)=\phi(x),\    W_*(x)=W_1(x),
\eeq 
we have 
\beq\label{EstI9}
 I_9 \le \varep_2 I_0 +{\widehat\Phi}'_1
(1+ J_0)^{1+\widehat\mu'_3/\alpha},
\eeq
where ${\widehat\Phi}'_1$ and $\widehat\mu'_3$  are the numbers ${\widehat\Phi}_1$ and $\widehat\mu_3$ in formulas \eqref{phi00} and \eqref{hatmu1}, respectively, applied to \eqref{app}.
Above estimates \eqref{HI9} and \eqref{EstI9} are natural taking into account the available tools in Lemma \ref{RoS1}. However, we chose to present \eqref{CI9} for the sake of simplicity. 
Indeed, \eqref{CI9} was conveniently combined the estimate \eqref{I2ineqq} of $I_2$, see the same terms $\mathcal K_2^{\theta}\mathcal K_1^{1-\theta}$ and ${\mathcal K}_3^\frac\theta{1-\theta}\mathcal K_1 {\mathcal K}_0^{\frac{2\lambda\theta}{\alpha(1-\theta)}}$ in both \eqref{Phat1} and \eqref{Phat2}.
\end{remark}

We now obtain estimates for the solution $u(x,t)$ in terms of the initial and boundary data with an explicit dependence on $\chi_*$.

\begin{theorem}\label{Labound}
Let the numbers $\alpha$, $\mu_*$, $\gamma_*$
and  $\bar C_\alpha$ be as in Proposition \ref{Diff4u}, and $M_\alpha(t)$ be defined by \eqref{Mtdef}. 
Set $V_{\alpha,0}=1+\int_U \phi(x) u_0(x)^{\alpha}\d x$.
Suppose there is a number $T>0$  such that
 \beq\label{tsmall1}
  \int_0^T M_\alpha(\tau)\d\tau \le  \frac{\alpha}{2\mu_* \chi_*^{\gamma_*}\bar C_\alpha  } V_{\alpha,0}^{-\mu_*/\alpha} .
\eeq 
For $t\in[0,T)$, define 
$\delta(t)=1-2\mu_* \chi_*^{\gamma_*} \alpha^{-1}\bar C_\alpha   V_{\alpha,0}^{\mu_*/\alpha} \int_0^t M_\alpha(\tau) \d\tau >0.$
\begin{enumerate}[label=\tnum]
    \item\label{Le1} One has, for all $t\in[0,T)$,
  \beq\label{uwest}
    \left(\int_U u^\alpha(x,t)\phi(x) \d x \right)^{1/\alpha}\le 
    V_{\alpha,0}^{1/\alpha} \delta(t)^{-1/\mu_*}.
    \eeq
\item\label{Le2} If $\eta$ is a number in the interval $(0,\alpha)$ such that  
$$C_{\alpha,\eta}\eqdef \int_U \phi(x)^{-\frac\eta{\alpha-\eta}} \d x\text{ is finite,}$$
then one has, for all $t\in[0,T)$, 
\beq\label{ube}
    \left(\int_U u^\eta(x,t) \d x\right)^{1/\eta} 
    \le  C_{\alpha,\eta}^{1/\eta-1/\alpha} V_{\alpha,0}^{1/\alpha}\delta(t)^{-1/\mu_*}.
    \eeq
\end{enumerate}
\end{theorem}
\begin{proof}
\ref{Le1} For $t\in(0.T)$, define $V_\alpha(t)= 1+\int_U \phi(x) u(x,t)^\alpha \d x$. 
Neglecting the second integral on the left-hand side of of \eqref{J0dbest}, we obtain
\beqs 
    V_\alpha'(t) \le \bar C_\alpha\chi_*^{\gamma_*}\left( V_\alpha(t)^{1+\mu_*/\alpha}
 +M_\alpha(t)\right) .
\eeqs
Since $V_\alpha(t)\ge 1$ and $M_\alpha(t)\ge 1$, it follows that 
\beqs 
    V_\alpha'(t)  \le 2\bar C_\alpha\chi_*^{\gamma_*} M_\alpha(t) V_\alpha(t)^{1+\mu_*/\alpha}.
\eeqs 
This inequality yields 
\begin{align*}
  V_\alpha(t) &\le \left (V_{\alpha,0}^{-\mu_*/\alpha} -\frac{2\mu_* \bar C_\alpha\chi_*^{\gamma_*} }{\alpha} \int_0^t M_\alpha(\tau)\d \tau\right)^{-\alpha/\mu_*}
  =V_{\alpha,0}\delta(t)^{-\alpha/\mu_*}\text{ for } t\in[0,T).
\end{align*}
Combining this estimate with $\|u(\cdot,t)\|_{L^\alpha_\phi(U)}\le (V_\alpha(t))^{1/\alpha}$ implies  \eqref{uwest}.

\ref{Le2} Applying H\"older's inequality with powers $\alpha/\eta$ and $\alpha/(\alpha-\eta)$, one has 
\beqs 
\left(\int_U u(x,t)^\eta  \d x \right)^{1/\eta}
= \left(\int_U (u^\eta \phi^{\eta/\alpha} )\cdot \phi^{-\eta/\alpha}   \d x \right)^{1/\eta}
\le \left( \int_U u(x,t)^{\alpha}\phi(x)  \d x \right )^{1/\alpha} C_{\alpha,\eta}^{1/\eta-1/\alpha} .
\eeqs
Together with the estimate \eqref{uwest}, this implies \eqref{ube}.
\end{proof}

\section{Estimates of the spatiotemporal $L^\infty$-norm}\label{maxsec}

Let $u(x,t)$ be a non-negative solution of the problem \eqref{ibvpg} as in Section \ref{Laest}.
We estimate the  $L^\infty$-norm in $(x,t)$ of $u$ by using a weighted Moser iteration. 
We briefly describe such a method below.

Find a strictly increasing sequence $\beta_j\to\infty$ as $j\to\infty$ and a decreasing sequence of sets $\mathcal Q_j\subset \R^{n+1}$ which converges to a set $\mathcal Q_*\subset \R^{n+1}$ as $j\to\infty$. 
Let $w(x)>0$ be a fixed weight function, and consider the following weighted $L_{x,t}^{\beta_j}$-norm   
$y_j=\left(\int_{\mathcal Q_j} |u(x,t)|^{\beta_j}w(x) \d x\d t \right)^\frac{1}{\beta_j}.$
For each $j$, we will establish the inequality 
\beq\label{Yjineq}
y_{j+1}\le A^\frac{\omega_j}{\kappa_j}(y_j^{r_j}+y_j^{s_j})^\frac1{\kappa_j}.
\eeq
Under certain conditions on $\omega_j$, $r_j$, $s_j$ and $\kappa_j$, we can iterate the above inequality \eqref{Yjineq} to have $\limsup_{j\to\infty} y_j$ is finite. In fact, we will have an estimate for this limit in terms of $y_0$, $A$, $\omega_j$, etc. Doing so, we can find an estimate for
$\|u\|_{L_{x,t}^{\infty}(\mathcal Q_*)}$.
The inequality \eqref{Yjineq} will be obtained by using suitable trace and Sobolev inequalities with weights. The original iteration in \cite{moser1971pointwise} has $w(x)\equiv 1$, $r_j=s_j=\kappa_j$, and $\beta_j=\kappa^j\beta_0$ with $\kappa$  simply being $1+2/n$ for $n>2$, and $5/3$ for $n=1,2$.

Recall that $c_5$, $c_6$ are the positive constants from Lemma \ref{trace} and  the weights  $W_3(x)$ and $W_4(x)$ are defined by \eqref{W3}. 
For $T>0$, denote $Q_T =U\times(0,T)$.

\begin{lemma} \label{caccio} 
Let numbers $\widetilde \kappa>1$ and $p_i>1$ for $i=1,2,\ldots,6$  satisfy
\beq\label{ktil}
p_1,p_2,p_3p_4<\widetilde \kappa\text{ and } p_6\eqdef \frac{p_5[p_3(2-a)-1]}{1-a}<\widetilde \kappa .
\eeq 
For $i=1,2\ldots, 6$, let $q_i$  be the H\"older conjugate exponent of $p_i$. 
Assume the integrals 
\beq\label{wcond1}
\int_U\left\{ \phi(x)+ \phi(x)^{-\frac{q_4}{p_4}}+W_4(x)^{q_1} \phi(x)^{-\frac{q_1}{p_1}}  + W_3(x)^{q_2}\phi(x)^{-\frac{q_2}{p_2}} 
 + W_1(x)^{-\frac{q_5}{1-a}}\phi(x)^{-\frac{q_5}{p_5}}
\right\} \d x
\eeq 
and 
\beq\label{Psicond}
\int_0^T \int_{\Gamma} (\psi_1^-(x,t))^{q_3}+ (\psi_2^-(x,t))^{q_3} \d S\d t
\eeq 
are finite. Then there are constants $c_8>1$ and $c_9>1$ depending on  $U$, $a$, $\lambda$, $c_{\mathcal Z}$, $p_3$ such that the following statement holds true. 
If $T>T_2>T_1\ge 0$ and  $\alpha$ is a number satisfying 
\beq\label{alp-Large}
\alpha>\max\left\{\lambda+1,\frac{p_2(\lambda(5-4a)-1)}{\widetilde\kappa-p_2},   \frac{p_5(a+3\lambda-1)}{(1-a)(\widetilde\kappa -p_6)}\right\},
\eeq
then one has  
\beq\label{newS1}
\begin{aligned}
&\essup_{t\in(T_2,T)} \int_U u^\alpha(x,t)\phi(x)\d x \\
&\le  c_8\chi_*^\gamma \alpha^{\frac{3-2a}{1-a}}(1+T)\left(1+\frac1{T_2-T_1}\right) \mathcal N_1 \Psi_T 
\Big(\norm {u}_{L_\phi^{\widetilde \kappa\alpha}(U\times(T_1,T))}^{\alpha-h_1}+\norm {u}_{L_\phi^{\widetilde \kappa\alpha}(U\times(T_1,T))}^{\alpha+h_2}\Big),
\end{aligned}
\eeq
and
\beq\label{newS2}
\begin{aligned}
&\int_{T_2}^T\int_U \frac{u^{\alpha-\lambda-1}(x,t)}{(1+u(x,t))^{2\lambda}}|\nabla u(x,t)|^{2-a}W_1(x)\d x \d t\\
&\le c_9\chi_*^{2+\gamma}\alpha^{\frac{2-a}{1-a}}(1+T)\left(1+\frac1{T_2-T_1}\right) \mathcal N_1 \Psi_T  
\Big(\norm {u}_{L_\phi^{\widetilde \kappa\alpha}(U\times(T_1,T))}^{\alpha-h_1}+\norm {u}_{L_\phi^{\widetilde \kappa\alpha}(U\times(T_1,T))}^{\alpha+h_2}\Big),
\end{aligned}
\eeq
where 
 \beq\label{gamdef}
 \gamma=2\max\left\{3-2a,\frac1{p_3(2-a)-1}\right\},
 \eeq
 \beq \label{h12}
h_1=\lambda+1>1, \  
h_2=\max \left\{0,\lambda(5-4a)-1, \frac{a+3\lambda-1}{p_3(2-a)-1}\right\}\ge 0,
\eeq 
\begin{align}
   \mathcal N_1
&=\left(1+\int_U\phi(x) \d x\right)\Bigg\{ 1+ \left(\int_U W_4(x)^{q_1} \phi(x)^{-\frac{q_1}{p_1}} \d x \right)^\frac{1}{q_1} + \left(\int_U W_3(x)^{q_2}\phi(x)^{-\frac{q_2}{p_2}} \d x \right)^\frac{1}{q_2}\notag \\
&\quad+\left(\int_U \phi(x)^{-\frac{q_4}{p_4}} \d x \right)^\frac{1}{p_3q_4}
+\left(\int_U W_1(x)^{-\frac{q_5}{1-a}}\phi(x)^{-\frac{q_5}{p_5}} \d x\right)^\frac{1-a}{q_5(p_3(2-a)-1)}\Bigg\}, \label{N1} \\
\label{bcquant}
\Psi_T&= 1+\left(\int_0^T \int_{\Gamma} (\psi_1^-(x,t))^{q_3}+ (\psi_2^-(x,t))^{q_3} \d S\d t\right)^{\frac{p_3(2-a)}{q_3(p_3(2-a)-1)}}.
\end{align}
\end{lemma}
\begin{proof}
We proceed in multiple steps.

\medskip\noindent\emph{Step 1.}
Let $\xi=\xi(t)$ be a $C^1$-function  on $[0,T]$ with    
\begin{align}\label{xiprop} 
&0\le \xi(t)\le 1 \text{  on }[0,T],\ 
\xi(t)=0 \text{ on } [0,T_1], \ 
\xi(t)=1 \text{ on } [T_2,T],\\
 \label{xiprime}
 &0\le \xi'(t)\le \frac{2}{T_2-T_1} \text{ on } [0,T].
\end{align} 
Choosing the test function $\varphi(x,t)=u^{\alpha-\lambda}\xi^2$ in \eqref{weakf} gives
\begin{align*}
\int_U \phi(u^\lambda)_t u^{\alpha-\lambda}\xi^2 \d x
&= -(\alpha-\lambda)\int_U  X (x, u, \nabla u+u^{2\lambda}\mathcal Z(x,t))\cdot(\nabla u)  u^{\alpha-\lambda-1}\xi^2 \d x\\
&\quad -\int_\Gamma  (\psi_1 +\psi_2  u^\lambda)u^{\alpha-\lambda}\xi^2 \d S.    
\end{align*}
Same as \eqref{ez}, we have
\beqs 
\ddt\int_U \phi u^\alpha \xi^2 \d x
=\int_U  \phi \left\{\frac{\partial}{\partial t}\left[(u^\lambda)^{\alpha/\lambda}\right]\xi^2 + 2 u^\alpha \xi\xi'\right\} \d x
=\frac{\alpha}{\lambda}\int_U \phi(u^\lambda)_t u^{\alpha-\lambda}\xi^2 \d x
+2\int_U \phi u^{\alpha} \xi \xi' \d x.
\eeqs 
Therefore, we obtain
\beq \label{dJxi}
 \begin{aligned}
\frac{\lambda}{\alpha}\ddt \int_U \phi u^{\alpha} \xi^2 \d x
&=-(\alpha-\lambda) \int_U X(x, u, \nabla u+u^{2\lambda}\mathcal Z(x,t)) \cdot (\nabla u)u^{\alpha-\lambda-1}\xi^2 \d x\\
&\quad - \int_\Gamma (\psi_1 u^{\alpha-\lambda}+\psi_2u^\alpha)  \xi^2 \d S
+\frac{2\lambda}{\alpha}\int_U \phi u^{\alpha} \xi \xi' \d x.
\end{aligned}
\eeq 

Let $\varep_1>0$.
We estimate the first term on the right-hand side of \eqref{dJxi} by multiplying inequality \eqref{negI1} by $\xi^2(t)$.
Then we obtain from this and \eqref{dJxi} a similar inequality to \eqref{J0dif0}, namely,
\beq \label{bint4}
\begin{aligned}
&\frac{\lambda}{\alpha}\ddt\int_U \phi u^\alpha \xi^2 \d x
+(\alpha-\lambda)( c_0 \chi_*^{-2}-\varep_1) \int_U \frac{u^{\alpha-\lambda-1}}{(1+u)^{2\lambda}}|\nabla u|^{2-a} W_1\xi^2\d x\\
&\le C_1(\alpha-\lambda)\left[\int_U W_4 u^{\alpha-\lambda-1}\xi^2\d x+\int_U W_3u^{\lambda(5-4a)+\alpha-1} \xi^2\d x\right]\\
&\quad - \int_\Gamma \psi_1u^{\alpha-\lambda} \xi^2 \d S- \int_\Gamma\psi_2u^{\alpha} \xi^2 \d S
+\frac{2\lambda}{\alpha}\int_U \phi u^{\alpha} \xi \xi' \d x,
\end{aligned}
\eeq
where $c_0=2^{a-1}$ and $C_1$ is defined in \eqref{C1def}.
Denote
\beq\label{BAo}
B_0=\essup_{t\in(0,T)}\int_U \phi(x) u(x,t)^\alpha \xi^2(t) \d x,\  
A_0=\iint_{Q_T} \frac{u(x,t)^{\alpha-\lambda-1}}{(1+u(x,t))^{2\lambda}}|\nabla u(x,t)|^{2-a} W_1(x)\xi^2(t) \d x\d t.
\eeq 
Note that
$\int_U \phi(x) u^\alpha(x,0) \xi^2(0) \d x=0$.
Then integrating inequality \eqref{bint4} in time from $0$ to $t$, for any $t\in(0,T)$, and  taking  the essential supremum with respect to $t$ yield
\beq \label{Bzest}
\begin{aligned}
\frac{\lambda}{\alpha}B_0\le C_1(\alpha-\lambda)(A_1+A_2)+ A_3 +A_4+\frac{2\lambda}{\alpha}A_5,
\end{aligned}
\eeq 
where 
\begin{align*}
    A_1&=\iint_{Q_T} W_4(x) u(x,t)^{\alpha-\lambda-1}\xi^2(t)\d x \d t,\  
     &&A_2=\iint_{Q_T}  W_3(x) u(x,t)^{\lambda(5-4a)+\alpha-1} \xi^2(t)\d x \d t,\\ 
 A_3&=\int_0^T\int_\Gamma \psi_2^-(x,t)u(x,t)^\alpha \xi^2(t) \d S \d t, \ 
 &&A_4=\int_0^T\int_\Gamma\psi_1^-(x,t)u(x,t)^{\alpha-\lambda} \xi^2(t) \d S \d t,\\ 
 A_5&= \iint_{Q_T} \phi(x) u(x,t)^{\alpha} \xi \xi' \d x \d t.
\end{align*}

Now, just integrating inequality \eqref{bint4} in time from $0$ to $T$ gives  
\beq \label{Azest}
(\alpha-\lambda)( c_0 \chi_*^{-2}-\varep_1)A_0\le C_1(\alpha-\lambda)(A_1+A_2)+ A_3 +A_4+\frac{2\lambda}{\alpha}A_5.
\eeq
Adding \eqref{Bzest} to \eqref{Azest} yields
\beq\label{MainD}
\frac{\lambda}{\alpha}B_0
+(\alpha-\lambda)( c_0 \chi_*^{-2}-\varep_1)A_0
\le   2C_1(\alpha-\lambda)(A_1+A_2)+ 2(A_3 +A_4)+\frac{4\lambda}{\alpha}A_5.
\eeq

\medskip\noindent\emph{Step 2.} We estimate the right-hand side of \eqref{MainD}.
Denote $$E(\alpha)=\iint_{Q_T} \phi(x)u^{\alpha}(x,t)\xi(t)  \d x \d t\text{ for }\alpha>0.$$  
In estimates below, we will often use the fact $0\le \xi^2\le \xi\le 1$.

\medskip
\noindent\textit{Estimation of $A_1$.}
Let 
\beqs 
\alpha_1=p_1(\alpha-\lambda-1)\text{ and } F_1=\left(\iint_{Q_T} W_4(x)^{q_1} \phi(x)^{-q_1/p_1}\xi^2(t) \d x\d t\right)^{1/q_1}.
\eeqs 
By applying H\"older's inequality with the powers $p_1$ and $q_1$, we have
\beq\label{A1est}
\begin{aligned}
A_1& =\iint_{Q_T}  (\phi^{1/p_1}u^{\alpha-\lambda-1})\cdot  (W_4(x) \phi^{-1/p_1}) \xi^2\d x \d t  \\
&\le \left(\iint_{Q_T} \phi(x) u^{\alpha_1}\xi^2 \d x \d t\right)^{1/p_1} F_1 \le E(\alpha_1)^{1/p_1}F_1.
\end{aligned}
\eeq

\medskip
\noindent\textit{Estimation of $A_2$.}
Let 
\beqs 
\alpha_2=p_2(\lambda(5-4a)+\alpha-1)\text{ and }  F_2=\left(\iint_{Q_T} W_3(x)^{q_2}\phi(x)^{-q_2/p_2} \xi^2(t)\d x\d t\right)^{1/q_2}.
\eeqs 
By applying H\"older's inequality with the powers $p_2$ and $q_2$, we have, same as \eqref{A1est}, 
\beq\label{A3est}
\begin{aligned}
A_2& = \iint_{Q_T} (\phi^{1/p_2}u^{\lambda(5-4a)+\alpha-1})\cdot (W_3(x)\phi^{-1/p_2}) \xi^2\d x\d t\\
&\le \left(\iint_{Q_T} \phi(x) u^{\alpha_2}\xi^2 \d x \d t\right)^{1/p_2} F_2 \le E(\alpha_2)^{1/p_2}F_2.  
\end{aligned}
\eeq

\medskip
\noindent\textit{Estimation of $A_3$.}
Let 
\beq\label{alpha3}
\alpha_3=p_3 \alpha \text{ and }F_3= \left(\int_0^T \int_{\Gamma} (\psi_2^-(x,t))^{q_3} \xi^2(t)\d S\d t\right)^{1/q_3}.
\eeq
We first apply H\"older's inequality with the powers $q_3$ and $p_3$ to obtain
\beq\label{A4est}
    A_3
        \le  \left(\int_0^T \int_{\Gamma} (\psi_2^-(x,t))^{q_3} \xi^2(t)\d S\d t\right)^{1/q_3} \left(\int_0^T \int_{\Gamma} u^{p_3\alpha} \xi^2\d S\d t\right)^{1/p_3} 
    =F_3  J_1^{1/p_3} ,
\eeq 
 where $ J_1= \int_0^T \int_{\Gamma} u^{\alpha_3} \xi^2\d S\d t$.

Applying the trace theorem \eqref{firstrace} to the  function $f:=u^{\alpha_3}$, we have
\begin{equation}\label{i4est}
    J_1 \le c_5 J_2+c_6\alpha_3 J_3, \text{  where } J_2=\iint_{Q_T} u^{\alpha_3}\xi^2 \d x \d t,\quad J_3=\iint_{Q_T} u^{\alpha_3-1}|\nabla u|\xi^2 \d x \d t.
\end{equation}
Raising both sides of inequality \eqref{i4est} to the power $1/p_3$ and applying inequality \eqref{ee2}  to $p=1/p_3<1$ yields
\begin{equation}\label{J1p3}
    J_1^{1/p_3}\le (c_5 J_2+c_6\alpha_3 J_3)^{1/p_3} \le c_5^{1/p_3} J_2^{1/p_3}+(c_6\alpha_3)^{1/p_3} J_3^{1/p_3}.
\end{equation}
For $J_2$, letting 
\beqs 
\alpha_4=p_4\alpha_3=p_4p_3\alpha\text{ and }  F_4=\left(\iint_{Q_T} \phi(x)^{-q_4/p_4}\xi^2(t) \d x \d t\right)^{1/q_4},
\eeqs 
and applying H\"older's inequality with the powers $p_4$ and $q_4$, we have, same as \eqref{A1est}, 
\beq\label{I4part1}
J_2= \iint_{Q_T} (u^{\alpha_3}\phi^{1/p_4})\cdot \phi^{-1/p_4}\cdot \xi^2 \d x \d t 
\le E(\alpha_4)^{1/p_4} F_4.
\eeq
To  estimate $J_3$, we set the numbers
\beq\label{mfour}
\begin{aligned}
m_4&=\frac{(\alpha_3-1)(2-a)-(\alpha-\lambda-1)}{1-a}, \ 
\alpha_5=p_5 m_4,\ 
\alpha_6=p_5\left(m_4+\frac{2\lambda}{1-a}\right).
\end{aligned}
\eeq
Applying H\"older's inequality with the powers $(2-a)$ and $(2-a)/(1-a)$, we have 
\begin{align*}
J_3&
=\iint_{Q_T}\left( \frac{u^\frac{\alpha-\lambda-1}{2-a}}{(1+u)^\frac{2\lambda}{2-a}} |\nabla u| W_1^\frac1{2-a}\right)
\cdot \left( W_1^{-\frac1{2-a}} u^{\alpha_3-1-\frac{\alpha-\lambda-1}{2-a}} (1+u)^\frac{2\lambda}{2-a}\right)   \xi^2 \d x \d t\\
&
 \le  \left(\iint_{Q_T} \frac{u^{\alpha-\lambda-1}}{(1+u)^{2\lambda}} |\nabla u|^{2-a}W_1 \xi^2 \d x \d t\right)^\frac{1}{2-a}
    \left(\iint_{Q_T}W_1^{-\frac1{1-a}} u^{m_4}(1+u)^\frac{2\lambda}{1-a}\xi^2 \d x \d t\right)^\frac{1-a}{2-a}.
\end{align*}
Using \eqref{ee2} to estimate $(1+u)^\frac{2\lambda}{1-a}\le 2^\frac{2\lambda}{2-a}(1+u^\frac{2\lambda}{1-a})$ in the last integral,
we obtain
\begin{align*}
J_3&
\le A_0^\frac{1}{2-a} 2^\frac{2\lambda}{2-a}\left(\iint_{Q_T} W_1^{-\frac1{1-a}}
    u^{m_4}  \xi^2 \d x \d t+\iint_{Q_T} W_1^{-\frac1{1-a}}
    u^{m_4+\frac{2\lambda}{1-a}}  \xi^2 \d x \d t\right)^\frac{1-a}{2-a}\\
    &= 2^\frac{2\lambda}{2-a} A_0^\frac{1}{2-a} \Bigg\{ \iint_{Q_T} \left(W_1^{-\frac1{1-a}}\phi^{-\frac{1}{p_5}}\right)
\cdot \left(\phi^\frac{1}{p_5} u^{m_4} \right) \xi^2 \d x \d t\\
&\quad +\iint_{Q_T} \left(W_1^{-\frac1{1-a}}\phi^{-\frac{1}{p_5}}\right)\cdot 
    \left(\phi^\frac{1}{p_5}u^{m_4+\frac{2\lambda}{1-a}}\right)  \xi^2 \d x \d t\Bigg\}^\frac{1-a}{2-a}.
\end{align*}
Applying H\"older's inequality again with the powers $q_5$ and $p_5$ to estimate the  last two integrals yields 
\beqs
J_3   \le  2^\frac{2\lambda}{2-a}A_0^\frac{1}{2-a} \left( F_5 E(\alpha_5)^\frac{1}{p_5} 
     +F_5 E(\alpha_6)^\frac{1}{p_5}\right)^\frac{1-a}{2-a}
     =2^\frac{2\lambda}{2-a}A_0^\frac{1}{2-a}F_5^\frac{1-a}{2-a}  \left(  E(\alpha_5)^\frac{1}{p_5} 
     + E(\alpha_6)^\frac{1}{p_5}\right)^\frac{1-a}{2-a},
\eeqs 
 where 
\beqs 
F_5=\left(\iint_{Q_T} W_1(x)^{-\frac{q_5}{1-a}}\phi(x)^{-\frac{q_5}{p_5}} \xi^2(t) \d x \d t\right)^{1/q_5}.
\eeqs 
Raising this inequality to the power $1/p_3$, and applying inequality \eqref{ee3} to $p=\frac{1-a}{p_3(2-a)}<1$
 to estimate $(  E(\alpha_5)^\frac{1}{p_5} 
     + E(\alpha_6)^\frac{1}{p_5})^\frac{1-a}{p_3(2-a)}$, we obtain
\beq\label{I4part2}
J_3^{1/p_3}  \le  2^\frac{2\lambda}{p_3(2-a)}  A_0^\frac{1}{p_3(2-a)} F_5^\frac{1-a}{p_3(2-a)}\left( E(\alpha_5)^\frac{1-a}{p_3p_5(2-a)}
     +E(\alpha_6)^\frac{1-a}{p_3p_5(2-a)}\right).
\eeq 
Combining \eqref{J1p3} with \eqref{I4part1} and  \eqref{I4part2} gives
\begin{equation}\label{i4estf}
     J_1^{1/p_3} \le c_5^{1/p_3} E(\alpha_4)^\frac{1}{p_3p_4} F_4^{1/p_3}     +2^\frac{2\lambda}{p_3(2-a)}(c_6\alpha_3 )^{1/p_3}A_0^\frac{1}{p_3(2-a)} J_4,
\end{equation}
where
\beq\label{J4} 
J_4=F_5^{\frac{1-a}{p_3(2-a)}}\left(E(\alpha_5)^{\frac{1-a}{p_3p_5(2-a)}} +E(\alpha_6)^{\frac{1-a}{p_3p_5(2-a)}}\right).
\eeq 

Given $\varep>0$. 
Now utilizing estimate \eqref{i4estf} in \eqref{A4est} yields
\begin{align*}
    A_3 &\le c_5^{1/p_3} E(\alpha_4)^\frac{1}{p_3p_4} F_4^{1/p_3}F_3  +2^\frac{2\lambda}{p_3(2-a)}(c_6\alpha_3 )^{1/p_3}A_0^\frac{1}{p_3(2-a)}
J_4F_3.
\end{align*}
Rewriting the last term  and applying Young's inequality with powers $p_3(2-a)$ and $\frac{p_3(2-a)}{p_3(2-a)-1}$ give 
\begin{align*}
&2^\frac{2\lambda}{p_3(2-a)}(c_6\alpha_3 )^{1/p_3}A_0^\frac{1}{p_3(2-a)}
J_4F_3
  =\left(\varep^\frac{1}{p_3(2-a)} A_0^\frac{1}{p_3(2-a)}\right)\cdot 
\left( \varep^{-\frac{1}{p_3(2-a)}}  2^\frac{2\lambda}{p_3(2-a)}(c_6\alpha_3 )^{1/p_3}F_3 J_4\right) \\
    &\le  \varep A_0
    +\left[\varep^{-1}2^{2\lambda} (c_6 \alpha_3)^{2-a}\right]^{\frac{1}{p_3(2-a)-1}} F_3^\frac{p_3(2-a)}{p_3(2-a)-1} J_4^\frac{p_3(2-a)}{p_3(2-a)-1}.
\end{align*}
In order to estimate $J_4^\frac{p_3(2-a)}{p_3(2-a)-1}$, we use formula \eqref{J4} and apply inequality \eqref{ee3} with $p=\frac{p_3(2-a)}{p_3(2-a)-1}$ to have 
\begin{align*}
 J_4^\frac{p_3(2-a)}{p_3(2-a)-1}
 &= F_5^{\frac{1-a}{p_3(2-a)-1}} \left(E(\alpha_5)^{\frac{1-a}{p_3p_5(2-a)}} +E(\alpha_6)^{\frac{1-a}{p_3p_5(2-a)}}\right)^\frac{p_3(2-a)}{p_3(2-a)-1}\\
    &\le  F_5^{\frac{1-a}{p_3(2-a)-1}}\cdot  2^{\frac{p_3(2-a)}{p_3(2-a)-1}-1}\left( E(\alpha_5)^{\frac{1-a}{p_5(p_3(2-a)-1)}}+E(\alpha_6)^{\frac{1-a}{p_5(p_3(2-a)-1)}}\right)\\
 &   =2^\frac{1}{p_3(2-a)-1} F_5^{\frac{1-a}{p_3(2-a)-1}} \left( E(\alpha_5)^\frac{1}{p_6}+E(\alpha_6)^\frac{1}{p_6}\right).
\end{align*}
Therefore, we obtain
\beq\label{A4estb} 
\begin{aligned}
    A_3&\le \varep A_0 +c_5^{1/p_3} E(\alpha_4)^\frac{1}{p_3p_4}F_6 
    +  C'_{\varep,\alpha}F_7\left(E(\alpha_5)^\frac{1}{p_6}+E(\alpha_6)^\frac{1}{p_6}\right),
\end{aligned}
\eeq 
where 
\begin{align}
\label{Cep}
C'_{\varep,\alpha}&=\left[\varep^{-1}2^{2\lambda+1} (c_6 p_3\alpha)^{2-a}\right]^{\frac{1}{p_3(2-a)-1}},\\
\label{F78}
F_6&= F_4^{1/p_3}F_3
\text { and } 
F_7=F_5^{\frac{1-a}{p_3(2-a)-1}} F_3^\frac{p_3(2-a)}{p_3(2-a)-1}.
\end{align}

\medskip
\noindent\textit{Estimation of $A_4$.}
Similar to the estimate of $A_3$ above, by setting 
\beqs 
\alpha'_3=p_3(\alpha-\lambda),\  \alpha'_4=p_4p_3(\alpha-\lambda), 
\eeqs 
\beqs
m'_4=\frac{(\alpha_3'-1)(2-a)-(\alpha-\lambda-1)}{1-a},\ 
    \alpha'_5= p_5 m'_4,\ 
    \alpha'_6= p_5\left(m'_4+\frac{2\lambda}{1-a}\right),
\eeqs
and replacing $\psi_2^{-}$ by $\psi_1^{-}$, we have
\beq\label{A5est} 
    A_4\le \varep A_0 + c_5^{1/p_3} E(\alpha'_4)^\frac{1}{p_3p_4}F'_6    +C'_{\varep,\alpha-\lambda}F'_7\left( E(\alpha'_5)^\frac{1}{p_6}+E(\alpha'_6)^\frac{1}{p_6}\right),
\eeq 
where 
\beq\label{F783}
F'_6= F_4^{1/p_3}{\widehat F}_3,\  
F'_7=F_5^{\frac{1-a}{p_3(2-a)-1}} {\widehat F}_3^\frac{p_3(2-a)}{p_3(2-a)-1}
\text{ and }
{\widehat F}_3= \left(\int_0^T \int_{\Gamma} (\psi_1^-(x,t))^{q_3} \xi^2(t)\d S\d t\right)^{1/q_3}.
\eeq 

\medskip
\noindent\textit{Estimation of $A_5$.}
By \eqref{xiprime},  we have 
\beq\label{IIx}
A_5=\iint_{Q_T} \phi u^{\alpha} \xi \xi' \d x \d t
\le \frac{2}{T_2-T_1}\iint_{Q_T} \phi u^{\alpha} \xi  \d x \d t =\frac{2}{T_2-T_1} E(\alpha).
\eeq

Combining \eqref{MainD} with the estimates  \eqref{A1est}, \eqref{A3est}, \eqref{A4estb}, \eqref{A5est} together with $C'_{\varep,\alpha-\lambda}\le C'_{\varep,\alpha}$, and \eqref{IIx}, we obtain
\beq\label{MainD3}
\frac{\lambda}{\alpha}B_0
+(\alpha-\lambda)( c_0 \chi_*^{-2}-\varep_1)A_0
\le 4\varep A_0 + 2 C_1(\alpha-\lambda)J_5
+2c_5^{1/p_3} J_6+2C'_{\varep,\alpha} J_7
     +\frac{8\lambda}{\alpha(T_2-T_1)} E(\alpha),
\eeq
where 
\begin{align}
    \label{J56}
    J_5&=E(\alpha_1)^{1/p_1}F_1+E(\alpha_2)^{1/p_2}F_2,\quad 
    J_6=E(\alpha_4)^\frac{1}{p_3p_4}F_6 + E(\alpha'_4)^\frac{1}{p_3p_4} F'_6,\\
    \label{J7}
    J_7&=F_7\left(E(\alpha_5)^\frac{1}{p_6}+E(\alpha_6)^\frac{1}{p_6}\right)
   +F'_7\left (E(\alpha'_5)^\frac{1}{p_6} +E(\alpha'_6)^\frac{1}{p_6} \right).
\end{align}
In inequality \eqref{MainD3}, by choosing  $\varep_1=c_0\chi_*^{-2}/4$, same as in \eqref{eee},  and $\varep= c_0(\alpha-\lambda)\chi_*^{-2}/16 $  we obtain 
\beq\label{MainD2}
\begin{aligned}
 & \frac{\lambda}{\alpha}B_0
+(\alpha-\lambda)\frac{c_0 \chi_*^{-2}}{2}A_0
\le  \mathcal E_0\eqdef  2\left[ C_1(\alpha-\lambda)J_5
+c_5^{1/p_3} J_6+C_5 J_7
     +\frac{4\lambda}{\alpha(T_2-T_1)} E(\alpha)\right],
\end{aligned}
\eeq
where $C_1$ is now given by \eqref{C1ineq}, and $C_5=C'_{\varep,\alpha}$ from \eqref{Cep} now is
\beq\label{C5def}
C_5=\left[\frac{16\chi_*^2}{c_0(\alpha-\lambda)}\cdot 2^{2\lambda+1} (c_6 p_3\alpha)^{2-a}\right]^{\frac{1}{p_3(2-a)-1}}.
\eeq

\medskip\noindent\emph{Step 3.} 
For the left-hand side of \eqref{MainD2}, we connect it with the desired inequalities \eqref{newS1} and \eqref{newS2} by defining
\beqs 
B_*=\essup_{t\in(T_2,T)}\int_U \phi(x) u^{\alpha}(x,t)\d x
 \text{ and }
A_*=\int_{T_2}^T\int_U \frac{u(x,t)^{\alpha-\lambda-1}}{(1+u(x,t))^{2\lambda}}|\nabla u(x,t)|^{2-a}W_1(x)\d x \d t.
 \eeqs 
By \eqref{xiprop}, \eqref{BAo} and \eqref{MainD2}, we have 
\begin{align}
\label{Bse}
 B_*&= \essup_{t\in(T_2,T)}\int_U \phi(x) u^{\alpha}(x,t)\xi^2(t) \d x
\le B_0 \le \frac\alpha\lambda \mathcal E_0,\\
\label{Ase}
A_*&=\int_{T_2}^T\int_U \frac{u^{\alpha-\lambda-1}}{(1+u)^{2\lambda}} |\nabla u|^{2-a}W_1\xi^2 \d x \d t\le A_0\le \frac{2\chi_*^2 \mathcal E_0}{c_0(\alpha-\lambda)}\le 4\chi_*^2 \mathcal E_0.
\end{align}
The last inequality uses the properties $\alpha-\lambda\ge 1$ and $c_0=2^{a-1}\ge 1/2$.

 We estimate the right-hand side of \eqref{MainD2} next using 
\beqs
J_0\eqdef \norm {u}_{L_\phi^{\widetilde \kappa\alpha}(U\times(T_1,T))}=\left(\int_{T_1}^T\int_U \phi(x) u(x,t)^{\widetilde \kappa \alpha}\d x\d t\right)^{\frac{1}{\widetilde\kappa \alpha}}
\eeqs
and the following inequality. 
If $0<\beta<\widetilde \kappa \alpha$, then, by H\"older's inequality with powers $\widetilde\kappa \alpha/\beta$ and its H\"older conjugate, we have
\beq\label{X}
E(\beta)\le \int_{T_1}^T \int_U \phi u^{\beta}  \d x \d t \le J_0^{\beta} \left(\int_{T_1}^T \int_U \phi \d x\d t\right)^{1-\frac{\beta}{\widetilde \kappa \alpha}} \le  J_0^{\beta} (T+1) \Phi_*,
\eeq
where
\beq\label{Phistar}
\Phi_*=1+\int_U\phi \d x.
\eeq

We claim that
\beq\label{aapkap}
\alpha_i<\widetilde\kappa\alpha\text{ for }1\le i\le 6,\text{ and } 
\alpha'_j<\widetilde\kappa\alpha\text{ for }4\le j\le 6.
\eeq
Below, we quickly verify \eqref{aapkap}.

(a) From \eqref{ktil}, we have $p_1<\widetilde\kappa$, hence $\alpha_1=p_1(\alpha-\lambda-1)<p_1\alpha<\widetilde\kappa\alpha$. 

(b) We have $\widetilde\kappa>p_2$ from \eqref{ktil},  and $\alpha>\frac{p_2(\lambda(5-4a)-1)}{\widetilde\kappa-p_2}$ from \eqref{alp-Large}. These two inequalities imply $\alpha_2<\widetilde\kappa\alpha$. 

(c) From \eqref{ktil}, we have $p_3p_4<\widetilde\kappa$, hence $\alpha_4=p_3p_4\alpha<\widetilde\kappa\alpha$. Consequently, $\alpha_3=\alpha_4/p_4<\alpha_4<\widetilde\kappa\alpha$.

(d) From the formulas of $\alpha_3$  in \eqref{alpha3} and $m_4$ in \eqref{mfour}, we compute
\beq\label{a6p6}
\alpha_6=p_5\cdot \frac{(p_3\alpha-1)(2-a)-(\alpha-\lambda-1)+2\lambda}{1-a}
=p_6 \alpha  + \frac{p_5(a+3\lambda-1)}{1-a}.
\eeq 
Together with the facts $\widetilde\kappa>p_6$ from \eqref{ktil} and 
 $\alpha>\frac{p_5(a+3\lambda-1)}{(1-a)(\widetilde\kappa-p_6)}$ from \eqref{alp-Large}, this formula of $\alpha_6$ yields $\alpha_6<\widetilde\kappa\alpha$.
Also, $\alpha_5=\alpha_6-2\lambda p_5/(1-a)<\alpha_6$, hence, $\alpha_5 <\widetilde\kappa\alpha$.

Thus, (a)--(d) above prove the first part of \eqref{aapkap}.
Now, for $j=4,5,6$, one has  $\alpha'_j<\alpha_j$, then, thanks to (c) and (d) above, $\alpha'_j<\widetilde\kappa\alpha$. This proves the second part of \eqref{aapkap}. 

Applying inequality \eqref{X} to $\beta=\alpha$, $\beta=\alpha_i$,  for $i=1,2,4,5,6$, and $\beta=\alpha'_j$, for $j=4,5,6$, gives
\beq\label{Eaa}
E(\alpha)\le J_0^\alpha(T+1)\Phi_*,\quad E(\alpha_i)\le J_0^{\alpha_i}(T+1)\Phi_*, \quad E(\alpha'_j)\le J_0^{\alpha'_j}(T+1)\Phi_*.
\eeq
Next, using the fact $\xi^2(t)\le 1$, we estimate
\beq\label{FFi}
	F_i \le T^{1/q_i}\widetilde F_i\text{  for $i=1,2,\ldots,5$, and } {\widehat F}_3 \le \widetilde F_3,
\eeq
where
\beq\label{allFtil}
\begin{aligned}
\widetilde F_1&=\left(\int_U W_4(x)^{q_1} \phi(x)^{-\frac{q_1}{p_1}} \d x \right)^{1/q_1},&&
\widetilde F_2 =\left(\int_U W_3(x)^{q_2}\phi(x)^{-\frac{q_2}{p_2}} \d x \right)^{1/q_2},\\
\widetilde F_3&= \left(\int_0^T \int_{\Gamma} (\psi_1^-(x,t))^{q_3} + (\psi_2^-(x,t))^{q_3}\d S\d t\right)^{1/q_3},&&
\widetilde F_4 =\left(\int_U \phi(x)^{-\frac{q_4}{p_4}} \d x \right)^{1/q_4},\\
\widetilde F_5&=\left(\int_U W_1(x)^{-\frac{q_5}{1-a}}\phi(x)^{-\frac{q_5}{p_5}} \d x\right)^{1/q_5}.
\end{aligned}
\eeq 
Because the integrals in \eqref{wcond1} and \eqref{Psicond} are finite, all the above $\widetilde F_i$, for $1\le i\le 5$, are finite numbers.
From the definitions of $F_6$, $F_7$,  $F'_6$, $F'_7$ in \eqref{F78}, \eqref{F783} and estimates in \eqref{FFi}, it follows that 
\beq\label{F78ineq}
\begin{aligned}
    F_6,F'_6&\le T^{\frac{1}{p_3 q_4}}\widetilde F_4^\frac{1}{p_3}\widetilde F_3,\\
     F_7,F'_7&\le T^{\frac{1-a}{q_5(p_3(2-a)-1)}}\widetilde F_5^{\frac{1-a}{p_3(2-a)-1}} \widetilde F_3^\frac{p_3(2-a)}{p_3(2-a)-1}=T^\frac{p_5-1}{p_6}\widetilde F_5^{\frac{1-a}{p_3(2-a)-1}} \widetilde F_3^\frac{p_3(2-a)}{p_3(2-a)-1}.
\end{aligned}
\eeq
Above, in the estimates for $F_7$ and $F'_7$, we used the identities
\beqs
\frac{1-a}{q_5(p_3(2-a)-1)}=\frac{(p_5-1)(1-a)}{p_5(p_3(2-a)-1)}=\frac{p_5-1}{p_6}.
\eeqs
 Then using \eqref{Eaa}, \eqref{FFi}, \eqref{F78ineq} to estimate $J_5$, $J_6$ and $J_7$ in \eqref{J56} and \eqref{J7}, we  obtain 
 \begin{align*}
 J_5&\le [J_0^{\alpha_1}(T+1)\Phi_*]^{1/p_1}\cdot  T^{1/q_1}\widetilde F_1
 + [J_0^{\alpha_2}(T+1)\Phi_*]^{1/p_2}\cdot T^{1/q_2} \widetilde F_2,\\
 J_6&\le \left [J_0^\frac{\alpha_4}{p_3p_4}+J_0^\frac{\alpha'_4}{p_3p_4}\right]
 [(T+1)\Phi_*]^\frac{1}{p_3 p_4} \cdot T^\frac{1}{p_3q_4}\widetilde  F_4^\frac{1}{p_3}\widetilde F_3,\\
J_7&\le  \left[J_0^{\frac{\alpha_5}{p_6}}+J_0^{\frac{\alpha_6}{p_6}}+J_0^{\frac{\alpha'_5}{p_6}}+J_0^{\frac{\alpha'_6}{p_6}}\right][(T+1)\Phi_*]^\frac{1}{p_6}
\cdot T^\frac{p_5-1}{p_6} \widetilde F_5^\frac{1-a}{p_3(2-a)-1} \widetilde F_3^\frac{p_3(2-a)}{p_3(2-a)-1}.
 \end{align*}
For $i=1,2,4$, we have 
$
[(T+1)\Phi_*]^{1/p_i} T^{1/q_i}\le (T+1)^{1/p_i+1/q_i}\Phi_*=(T+1)\Phi_*.
$
Thus,
 \begin{align}\label{J5e}
 J_5&\le (J_0^{\alpha_1/p_1}\widetilde F_1+J_0^{\alpha_2/p_2}\widetilde F_2)(T+1)\Phi_*
 \le (T+1)\Phi_*\left[J_0^{\alpha_1/p_1}+J_0^{\alpha_2/p_2}\right]\left(\widetilde F_1+\widetilde F_2\right),\\
 J_6&\le \left [J_0^\frac{\alpha_4}{p_3p_4}+J_0^\frac{\alpha'_4}{p_3p_4}\right][(T+1)\Phi_*]^\frac{1}{p_3}\widetilde  F_4^\frac{1}{p_3}\widetilde F_3 
 \le (T+1)\Phi_*\left [J_0^\frac{\alpha_4}{p_3p_4}+J_0^\frac{\alpha'_4}{p_3p_4}\right]\left(\widetilde  F_4^\frac{1}{p_3}\widetilde F_3\right).\label{J6e}
 \end{align}
We also have 
\beqs 
(T+1)^\frac1{p_6}T^\frac{p_5-1}{p_6}\Phi_*^\frac{1}{p_6} 
\le (T+1)^{\frac1{p_6}+\frac{p_5-1}{p_6}}\Phi_*
= (T+1)^\frac{p_5}{p_6}\Phi_*\le (T+1)\Phi_*.
\eeqs 
The last inequality is due to the fact $p_5<p_6$.
Thus,
\beq \label{J7e}
J_7\le (T+1)\Phi_*  \left[J_0^{\frac{\alpha_5}{p_6}}+J_0^{\frac{\alpha_6}{p_6}}+J_0^{\frac{\alpha'_5}{p_6}}+J_0^{\frac{\alpha'_6}{p_6}}\right]\left(\widetilde F_5^\frac{1-a}{p_3(2-a)-1} \widetilde F_3^\frac{p_3(2-a)}{p_3(2-a)-1}\right).
\eeq 

For the last term in \eqref{MainD2}, we simply use the first inequality in \eqref{Eaa} and the fact $\alpha>1$ to estimate
\beq\label{lastE}
\frac{E(\alpha)}{\alpha(T_2-T_1)} \le \frac{(T+1)\Phi_*}{T_2-T_1}J_0^{\alpha}.
\eeq

We now turn to the constants on the right-hand side of \eqref{MainD2}. Regarding the number $C_1$, we re-estimate  \eqref{C1ineq} by
\beq\label{tC1}
C_1\le \widetilde C_1 \chi_*^{2(3-2a)},
\text{ where }
\widetilde C_1=\max\left\{1,(2+2^{2(\lambda+1)(1-a)+1}c_0^{-1})c_{\mathcal Z}^{2-a}\right\}.
\eeq 
Thus,
\beq\label{CC1}
C_1(\alpha-\lambda)\le \widetilde C_1 \chi_*^{2(3-2a)}\alpha.
\eeq

To estimate the constant $C_5$, we utilize formula \eqref{C5def} together with the facts that $\alpha-\lambda\ge 1$ due to \eqref{alp-Large}, $c_0\ge 1/2$  due to $c_0=2^{a-1}$ in \eqref{bint4}, and $\alpha^\frac{2-a}{p_3(2-a)-1}\le \alpha^\frac{2-a}{1-a}$ due to both $\alpha$ and $p_3$ being greater than $1$. It results in 
\beq \label{C5e}
C_5\le \left[2^{6+2\lambda}(c_6 p_3)^{2-a}\right]^{\frac 1{p_3(2-a)-1}}\chi_*^{\frac{2}{p_3(2-a)-1}}\alpha^\frac{2-a}{1-a}.
\eeq 

Therefore, by combining \eqref{J5e}, \eqref{J6e}, \eqref{J7e}, \eqref{lastE} with \eqref{CC1}, \eqref{C5e}, we estimate  each term $C_1(\alpha-\lambda)J_5$, $c_5^{1/p_3} J_6$, $C_5 J_7$ and $\frac{4\lambda}{\alpha(T_2-T_1)} E(\alpha)$ in the formula \eqref{MainD2} of $\mathcal E_0$, and obtain
\beq\label{I0fm}
\mathcal E_0 \le 2\widetilde C_2 \chi_*^\gamma \alpha^{\frac{2-a}{1-a}}(T+1)\Big(1+\frac {1}{T_2-T_1}\Big)\Phi_*\widetilde F\cdot\widetilde J,
 \eeq
 where $\gamma$ is defined by \eqref{gamdef},
\begin{align}
\widetilde C_2 &= \max\left\{  \widetilde C_1,c_5^{1/p_3}, \left[4^{3+\lambda}(c_6 p_3)^{2-a}\right]^{\frac 1{p_3(2-a)-1}}, 4\lambda\right\},\label{tC2} \\
\widetilde F &=  \left(\widetilde F_1 + \widetilde F_2\right)
+\left(\widetilde F_3\widetilde  F_4^{1/p_3}\right)
 + \left(\widetilde F_5^\frac{1-a}{p_3(2-a)-1} \widetilde F_3^\frac{p_3(2-a)}{p_3(2-a)-1}\right)+ 1,\label{tF}\\
    \widetilde J&= \left[J_0^\frac{\alpha_1}{p_1}+ J_0^\frac{\alpha_2}{p_2}\right]
    +\left[J_0^\frac{\alpha_4}{p_3p_4}+J_0^\frac{\alpha'_4}{p_3p_4}\right]
    + \left[J_0^\frac{\alpha_5}{p_6}+J_0^{\frac{\alpha_6}{p_6}}+J_0^{\frac{\alpha'_5}{p_6}}+J_0^{\frac{\alpha'_6}{p_6}}\right]
    +J_0^{\alpha}.\label{tJ}
\end{align}

\medskip\noindent\emph{Step 4.} We simplify the last two terms $\widetilde F$ and $\widetilde J$ in \eqref{I0fm} even further.
From \eqref{tF}, we estimate $\widetilde F$ simply by
\beq\label{tFe}
\widetilde F \le \left(1+ \widetilde F_1 + \widetilde F_2+\widetilde  F_4^{1/p_3}
+\widetilde F_5^\frac{1-a}{p_3(2-a)-1} \right)\left(1+\widetilde F_3+ \widetilde F_3^\frac{p_3(2-a)}{p_3(2-a)-1}  \right)
=\mathcal N (\Psi_T+\widetilde F_3),
\eeq 
where we used the simple identity $\Psi_T=1+\widetilde F_3^\frac{p_3(2-a)}{p_3(2-a)-1}$ thanks to \eqref{bcquant} and \eqref{allFtil}, and 
 \beqs \mathcal N =1+ \widetilde F_1 + \widetilde F_2+\widetilde  F_4^{1/p_3}
+\widetilde F_5^\frac{1-a}{p_3(2-a)-1}.
\eeqs
Applying inequality \eqref{ee5}, one can estimate the last $\widetilde F_3$ appearing in \eqref{tFe} by  $\widetilde F_3\le \Psi_T$. Thus, 
\beq\label{tildeFfm}
\widetilde F \le 2\mathcal N \Psi_T.
\eeq

The following are relations among the powers of the $J_0$-terms in $\widetilde J$:
\begin{align*}
\frac{\alpha_1}{p_1}=\alpha-\lambda-1< &\frac{\alpha'_4}{p_3p_4}=\alpha-\lambda
<\alpha=\frac{\alpha_4}{p_3p_4},\quad 
\frac{\alpha_1}{p_1}=\alpha-\lambda-1<\frac{\alpha_2}{p_2}=\alpha+\lambda(5-4a)-1,\\
&\frac{\alpha'_5}{p_6}<\min\left\{\frac{\alpha_5}{p_6},\frac{\alpha'_6}{p_6}\right\}
\le\max\left\{\frac{\alpha_5}{p_6},\frac{\alpha'_6}{p_6}\right\}<\frac{\alpha_6}{p_6}.
\end{align*}
Let $s$ be any of nine powers of $J_0$ in $\widetilde J$. It follows from the above comparisons that
\beq\label{smm}
\min\left\{\alpha-\lambda-1,\frac{\alpha'_5}{p_6} \right\} \le s\le \max\left\{\alpha,\alpha+\lambda(5-4a)-1,\frac{\alpha_6}{p_6}\right\}.
\eeq
For the lower bound of $s$ in \eqref{smm}, we have 
\beqs
\frac{\alpha'_5}{p_6}
=\frac{p_5 m'_4}{p_6}=\frac{1-a}{p_3(2-a)-1} \cdot \frac{[p_3(\alpha-\lambda)-1](2-a)-(\alpha-\lambda-1)}{1-a}
=\alpha-\lambda-\frac{1-a}{p_3(2-a)-1}.
\eeqs
Since $p_3>1$, it follows that
\beqs 
\frac{\alpha'_5}{p_6}\ge \alpha-\lambda-\frac{1-a}{(2-a)-1}=\alpha-\lambda-1.
\eeqs 
Combining this with the first inequality in \eqref{smm} gives
\beq\label{smin} 
s\ge \alpha-\lambda-1=\alpha-h_1>0.
\eeq 
For the upper bound of $s$ in \eqref{smm}, recall from \eqref{a6p6} that 
\beqs
\frac{\alpha_6}{p_6}
=\alpha  + \frac{p_5(a+3\lambda-1)}{p_6(1-a)}
=\alpha+\frac{a+3\lambda-1}{p_3(2-a)-1}.
\eeqs
Combining this with the second inequality in \eqref{smm} gives
\beq\label{smax}
s\le \alpha+\max\left\{0, \lambda(5-4a)-1, \frac{a+3\lambda-1}{p_3(2-a)-1}\right\}=\alpha+h_2.
\eeq 
Hence, by \eqref{smin}, \eqref{smax} and inequality \eqref{ee4}, one has $J_0^s\le J_0^{\alpha-h_1}+J_0^{\alpha+h_2}$.
Summing up this inequality over nine terms of $\widetilde J$ in \eqref{tJ} gives
\beq\label{tildeJfm}
    \widetilde J\le 9(J_0^{\alpha-h_1}+J_0^{\alpha+h_2}).
\eeq

Now, combining  \eqref{I0fm} with \eqref{tildeFfm} and \eqref{tildeJfm}, we obtain
\beq\label{I0i}
\mathcal E_0 \le 36 \widetilde C_2\alpha^{\frac{2-a}{1-a}}\chi_*^\gamma (1+T) \Big(1+\frac {1}{T_2-T_1}\Big) \Phi_* \mathcal N \Psi_T (J_0^{\alpha-h_1} +J_0^{\alpha+h_2}).
\eeq
Note in \eqref{I0i}  that $\Phi_* \mathcal N=\mathcal N_1$.

\medskip\noindent\emph{Step 5.} On the one hand,  combining \eqref{Bse} with \eqref{I0i} gives
\beqs
B_* \le \frac\alpha\lambda \mathcal E_0 
  \le 36 \widetilde C_2 \lambda^{-1}\alpha^{\frac{3-2a}{1-a}}\chi_*^\gamma(1+T)\Big(1+\frac {1}{T_2-T_1}\Big) \mathcal N_1 \Psi_T (J_0^{\alpha-h_1} +J_0^{\alpha+h_2})
\eeqs 
which proves \eqref{newS1} with $c_8=36 \widetilde C_2 \lambda^{-1}$.
On the other hand, combining \eqref{Ase}  with \eqref{I0i} gives
\beqs
A_* \le 4\chi_*^2\mathcal E_0
  \le 144\widetilde C_2 \alpha^{\frac{2-a}{1-a}}\chi_*^{2+\gamma}(1+T)\Big(1+\frac {1}{T_2-T_1}\Big)  \mathcal N_1 \Psi_T (J_0^{\alpha-h_1} +J_0^{\alpha+h_2})
\eeqs
which proves  \eqref{newS2} with $c_9=144\widetilde C_2$.
It is seen from \eqref{tC1} and \eqref{tC2} that $\widetilde C_2$ depends on $c_5$, $c_6$, $a$, $\lambda$, $c_{\mathcal Z}$ and  $p_3$, which verifies the dependence of $c_8$ and $c_9$ on $U$, $a$, $\lambda$, $c_{\mathcal Z}$ and  $p_3$.
One also has from \eqref{tC1} and \eqref{tC2} that 
$\widetilde C_2\ge \max\{\widetilde C_1,4\lambda\}\ge \max\{1,4\lambda\}$. Therefore, the  specifically defined $c_8$ and $c_9$ above are larger than $1$.
\end{proof}

\begin{remark}\label{pexist}
The set of pairs $(p_3,p_5)\in (1,\infty)^2$ that satisfy  the last condition in \eqref{ktil} is non-empty. 
Indeed, since $\widetilde\kappa>1$, one has $\widetilde\kappa(1-a)+1>2-a$, hence there exists a number $p_3$ such that
\beqs
1<p_3<\frac{\widetilde\kappa(1-a)+1}{2-a}.
\eeqs
With such a value $p_3$, one has $\displaystyle \frac{\widetilde\kappa(1-a)}{p_3(2-a)-1}>1$, which yields the existence of a number $p_5$ satisfying 
\beqs
1<p_5<\frac{\widetilde\kappa(1-a)}{p_3(2-a)-1}.
\eeqs
The last inequality implies that $(p_3,p_5)$ satisfies the last condition in \eqref{ktil}.
\end{remark}

In the next lemma, we establish an inequality between the Lebesgue norms of the solution which lays the foundation for later Moser's iteration.

Hereafter,  $r_1$ is a fixed number in $(0,1)$ that satisfies \eqref{newro}, 
and the number $r_*$ is defined by \eqref{newrs}.

\begin{lemma}\label{GLk}
 Assume the integrals in \eqref{wcond1}, \eqref{Psicond} and
 \beq\label{wcond2}
 \int_U  \left(\phi(x)^{\frac2{r_*}-1}+(1+\phi(x)^{-1})^{\frac{r_1}{1-r_1}}(1+ W_1(x)^{-1})^{\frac{r_1}{(1-r_1)^2}} \right)\d x
\eeq
are finite. Define
\beq\label{GG}
\begin{aligned}
 \mathcal N_2&=\left(1+\int_U\phi(x)\d x\right)^{\min\{1-r_1,\frac{2\lambda}{1+\lambda}\}}\left ( \int_U  \phi(x)^{\frac{2}{r_*}-1}\d x\right )^\frac{r_*}{2}\\
&\quad \cdot \left[ 1+\int_U\left(1+\phi(x)^{-1}\right)^{\frac{r_1}{1-r_1}}\left(1+ W_1(x)^{-1}\right)^{\frac{r_1}{(1-r_1)^2}}\d x \right]^{\frac{1-r_1}{r_1}}.    
\end{aligned}
\eeq
Let numbers $\widetilde \kappa>1$ and $p_i>1$ for $i=1,2,\ldots,5$  satisfy \eqref{ktil},  and let $\alpha$ be a number that satisfies \eqref{alp-Large} and 
 \beq\label{Bigalp}
\alpha > \max\left\{\frac{2(\lambda+a-1)}{r_*},\frac{1-\lambda-a}{\widetilde\kappa-1}, \frac{2\lambda}{1-r_1} \right\}.
\eeq
Let $h_1$, $h_2$ be defined by \eqref{h12}, and
\beq\label{kapmax}
    \kappa=\kappa(\alpha)\eqdef 1+\frac {r_*}2+\frac{1-\lambda-a}{\alpha}>1.
\eeq 
If $T>T_2>T_1\ge 0$  then
\beq\label{bfinest2}
\| u\|_{L_\phi^{\kappa\alpha}(U\times (T_2,T))}
\le A_\alpha^\frac 1{\alpha}\Big( \| u\|_{L_\phi^{\widetilde \kappa \alpha}(U\times(T_1,T))}^{\nu_1}+\| u\|_{L_\phi^{\widetilde\kappa \alpha}(U\times(T_1,T))}^{\nu_2}\Big)^\frac 1{\alpha},
\eeq
where  
\beq\label{tilrsdef}
\nu_1=\nu_1(\alpha)\eqdef   \frac{\alpha-h_1}{1+\frac{1}{\alpha(1+r_*/2)}} ,\quad 
\nu_2=\nu_2(\alpha)\eqdef (\alpha+h_3)\left(1+\frac{3\lambda}{\alpha}\right),
\eeq
\beq\label{Atildef2}
A_\alpha = c_{10}  \alpha^{6+\frac{5}{2(1-a)}} 
\chi_*^{2+\gamma\bar\mu}\max\{1,\mathcal N_2\}\left[(1+T)\Big(1+\frac1{T_2-T_1}\Big)\mathcal N_1 \Psi_T\right]^{\bar\mu}
\eeq
with 
\beqs
h_3=\max\{h_2, 1-a-\lambda\},\quad 
\bar\mu=1+\frac{r_*}2+\min\left\{1-r_1,\frac{2\lambda}{\lambda+1}\right\},
\eeqs
and, referring to the positive constants $c_7$, $c_8$, $c_9$ in \eqref{ppsi1}, \eqref{newS1}, \eqref{newS2}, 
\beqs 
c_{10}=3\cdot 2^{11+2\lambda}  \max\{1,c_7^{2-a}\}(\max\{c_8,c_9\})^{5/2}.
\eeqs 
\end{lemma}
\begin{proof} 
We apply the parabolic Sobolev inequality \eqref{pssi6} in Lemma \ref{RoS2} to 
\beq\label{choice3}
\varphi(x)=\phi(x), \quad W_*(x)=W_1(x),\quad \beta=2\lambda,\quad 
p = 2 - a,\quad   
s = \lambda+1,
\eeq
and the interval $[T_2, T]$ in place of $[0,T]$. 
The conditions \eqref{alcond} and, respectively,  \eqref{albe} become
\beq\label{ac3} 
\alpha\ge\lambda+1,\ \alpha >\frac{1-a-\lambda}{1-a},\ \alpha>\frac{2(\lambda+a-1)}{r_*},
\text{ and, respectively, } \alpha\ge \frac{2\lambda}{1-r_1}.
\eeq
In this case, $\lambda+1>1>\frac{1-a-\lambda}{1-a}$. Hence, the second condition in \eqref{ac3} can be neglected. The rest of conditions in \eqref{ac3} are satisfied thanks to our assumptions \eqref{alp-Large} and   \eqref{Bigalp}.

The numbers $m$ and $\theta_0$ in \eqref{powdef} and \eqref{defkappa} become
\beqs 
m=\frac{\alpha+1-\lambda-a}{2-a}\in[1,\alpha)
\text{ and }
\theta_0 =\frac{1}{1+\frac{r_*\alpha}{2(\alpha+1-\lambda-a)}}\in(0,1),
\eeqs 
while the number $\kappa$ in \eqref{defkappa} takes the specific value in \eqref{kapmax}.
By formulas \eqref{Calpha} and \eqref{EEdef} for $\mathcal E_1$, $\mathcal E_2$ and $\mathcal E_3$ with the choice \eqref{choice3}, the number ${\mathcal E}_1^\frac\beta\alpha {\mathcal E}_2 {\mathcal E}_3$ in \eqref{tilph10} becomes 
\begin{align*}
 {\mathcal E}_1^\frac\beta\alpha {\mathcal E}_2 {\mathcal E}_3&=\left(1+\int_U\phi(x)\d x\right)^{\frac{2\lambda}{\alpha}}\left ( \int_U  \phi(x)^{\frac{2}{r_*}-1}\d x\right )^\frac{r_*}{2}\\
 &\quad \cdot 
\left[ 1+\int_U\left(1+\phi(x)^{-1}\right)^{\frac{r_1}{1-r_1}}\left(1+ W_1(x)^{-1}\right)^{\frac{r_1}{(1-r_1)^2}}\d x \right]^{\frac{1-r_1}{r_1}}.
\end{align*}
Regarding the first power on the right-hand side, we have, by \eqref{Bigalp}, $\alpha>\frac{2\lambda}{1-r_1}$, and, by \eqref{alp-Large}, $\alpha>\lambda+1$, 
hence, 
\beq \label{lamin}
2\lambda/\alpha<\min\{1-r_1,2\lambda/(\lambda+1)\}.
\eeq 
Therefore,
\beq\label{PN}
 {\mathcal E}_1^\frac\beta\alpha {\mathcal E}_2 {\mathcal E}_3
 \le \mathcal N_2.
\eeq

Set  
\begin{align}
I_0&=\int_{T_2}^T\int_U \frac{u(x,t)^{\alpha-\lambda-1}}{(1+u(x,t))^{2\lambda}}|\nabla u(x,t)|^{2-a}W_1(x) \d x \d t,\notag \\
\label{Idef}
I&=I_0+\int_{T_2}^T\int_U  u(x,t)^{\alpha+1-\lambda-a}\phi(x)\d x \d t\text{ and } 
J=\essup_{t\in(T_2,T)}\norm{u(\cdot,t)}_{L_\phi^\alpha(U)}.
\end{align}
Then we have from inequality \eqref{pssi6}, using also \eqref{PN} and the fact $m<\alpha$, that 
\beq\label{beginest2}
\| u\|_{L_\phi^{\kappa\alpha}(U\times (T_2,T))}
\le (C_6 \alpha^{\frac 1{r_1}})^\frac{1}{\kappa\alpha} I^{\frac1{\kappa\alpha}}\cdot \left(J^{1-\theta_0}
+J^{1-\widehat\theta_0}\right),
\eeq
where 
\beq\label{C4}
C_6=2^{1+2\lambda+\frac{1}{r_1}}c_7^{2-a} \mathcal N_2 
\text{ and } \widehat \theta_0=\theta_0-\frac{2\lambda}{\kappa\alpha}.
\eeq 
Below, we estimate $I$ and $J$ on the right-hand side of \eqref{beginest2} in terms of $\Upsilon\eqdef  \norm {u}_{L_\phi^{\widetilde \kappa\alpha}(U\times(T_1,T))}$.

Regarding $I$, we estimate the last integral in its definition in \eqref{Idef}.
Let $\Phi_*\in[1,\infty)$ be defined by \eqref{Phistar}.
It follows from the assumption \eqref{alp-Large}, respectively, \eqref{Bigalp} that 
\beqs 
\alpha>\lambda>a-1+\lambda,\text{ respectively, } 
\alpha>\frac{1-\lambda-a}{\widetilde\kappa-1}.
\eeqs 
These imply $0<\alpha+1-\lambda-a<\widetilde\kappa\alpha$.
Then we can apply H\"older's inequality with the powers $\frac{\widetilde\kappa\alpha}{\alpha+1-\lambda-a}$ and its H\"older's conjugate to have 
\beq\label{Ifirst}
\begin{aligned}
    \int_{T_2}^T\int_U  u^{\alpha+1-\lambda-a}\phi\d x \d t
&    \le \left(\int_{T_2}^T\int_U u^{\widetilde\kappa \alpha}\phi \d x \d t\right)^\frac{\alpha+1-\lambda-a}{\widetilde\kappa\alpha} 
    \left(\int_{T_2}^T\int_U \phi\d x\d t\right)^{1-\frac{\alpha+1-\lambda-a}{\widetilde\kappa \alpha}}\\
    &\le \Upsilon^{\alpha+1-\lambda-a}\cdot (1+T)\Phi_*.
    \end{aligned}
\eeq
Combining the formula of $I$ in \eqref{Idef} with estimate  \eqref{Ifirst} gives
\beq\label{Iest3}
I\le I_0+(1+T)\Phi_*\Upsilon^{\alpha+1-\lambda-a}. 
\eeq

To estimate $J$ in \eqref{beginest2} and $I_0$ in \eqref{Iest3}, we will use Lemma \ref{caccio}.
Denote 
\begin{align}\label{mathcalM}
    \mathcal M_0& =(1+T)\Big(1+\frac1{T_2-T_1}\Big)\mathcal N_1 \Psi_T,\quad
    \mathcal M_1=2\max\{c_8,c_9\}\alpha^{1+\frac{2-a}{1-a}}\mathcal M_0,\\
\label{Sdef}
\mathcal{S}&=\mathcal M_1 (\Upsilon^{\alpha-h_1}+\Upsilon^{\alpha+h_3}).
\end{align}
In calculations below, we very often use the facts
\beq \label{simplef}
\Phi_*,\mathcal N_1,\Psi_T\ge 1\text{  and }c_8,\alpha>1.
\eeq 
In fact, these already imply $\mathcal M_0\ge 1$ and $\mathcal M_1\ge 1$.
Note also that 
$$0<\alpha-h_1\le \text{ each number in } \{\alpha+1-\lambda-a, \alpha-h_1, \alpha+h_2\}\le \alpha+h_3.$$ 
Applying inequality  \eqref{ee4} yields
\beq \label{UUU}
\Upsilon^{\alpha+1-\lambda-a},\Upsilon^{\alpha-h_1},\Upsilon^{\alpha+h_2}\le \Upsilon^{\alpha-h_1}+ \Upsilon^{\alpha+h_3}.
\eeq 

It follows from \eqref{newS1}, \eqref{newS2} and \eqref{UUU} that 
\begin{align}\label{Sest1}
J^\alpha&\le c_8\chi_*^\gamma \alpha^{\frac{3-2a}{1-a}}\mathcal M_0 (\Upsilon^{\alpha-h_1}+ \Upsilon^{\alpha+h_2})
\le  \chi_*^\gamma\mathcal S,\\
\label{Sest2}
I_0&\le c_9 \chi_*^{2+\gamma}\alpha^{\frac{2-a}{1-a}}\mathcal M_0 (\Upsilon^{\alpha-h_1}+ \Upsilon^{\alpha+h_2})
\le \chi_*^{2+\gamma}\alpha^{-1}\mathcal S.
\end{align} 
Also, by \eqref{UUU} and \eqref{Sdef},  
\beq\label{Usim}
\Upsilon^{\alpha+1-\lambda-a}\le \Upsilon^{\alpha-h_1}+ \Upsilon^{\alpha+h_3}=\mathcal M_1^{-1}\mathcal S.
\eeq 
Moreover, one has, thanks to \eqref{N1}, \eqref{mathcalM} and \eqref{simplef}, that
\beq\label{TP}
(T+1)\Phi_*\le (T+1)\mathcal N_1\le \mathcal M_0\le \alpha^{-1-\frac{2-a}{1-a}}\mathcal M_1.
\eeq 
Combining \eqref{Iest3} with \eqref{Sest2}, \eqref{Usim} and  \eqref{TP} yields  
\beq\label{Iest}
I\le \chi_*^{2+\gamma}\alpha^{-1}\mathcal S+( \alpha^{-1-\frac{2-a}{1-a}}\mathcal M_1)\cdot( \mathcal M_1^{-1}\mathcal S)
\le \chi_*^{2+\gamma}\alpha^{-1}\mathcal S+ \alpha^{-1}\mathcal S
\le 2\alpha^{-1}\chi_*^{2+\gamma}\mathcal S.
\eeq

Now, combining \eqref{beginest2},  \eqref{Sest1}, \eqref{Iest} and the facts $\chi_*\ge 1$,  $\widehat\theta_0<\theta_0$ , we obtain
\beq\label{beg33}
\begin{aligned}
\| u\|_{L_\phi^{\kappa\alpha}(U\times (T_2,T))}
&\le (C_6 \alpha^{\frac 1{r_1}})^\frac{1}{\kappa\alpha} 
(2\alpha^{-1}\chi_*^{2+\gamma} \mathcal S)^{\frac1{\kappa\alpha}}\cdot 
\left(\chi_*^{\frac{\gamma(1-\theta_0)}{\alpha}}\mathcal S^\frac{1-\theta_0}{\alpha}+\chi_*^{\frac{\gamma(1-\widehat\theta_0)}{\alpha}}\mathcal S^\frac{1-\widehat\theta_0}{\alpha}\right)\\
&\le \left (2\alpha^{-1+\frac 1{r_1}}C_6 \chi_*^{2+\gamma+\kappa\gamma(1-\widehat\theta_0)}\right)^\frac{1}{\kappa\alpha}
  \mathcal S^\frac{1}{\kappa\alpha} \left [  2(\mathcal S^{1-\theta_0}+\mathcal S^{1-\widehat\theta_0})\right]^\frac{1}{\alpha}.
\end{aligned}
\eeq
For the sake of convenience, we denote
\beqs 
C_7=2\alpha^{-1+\frac 1{r_1}}\max\{1,C_6\} \chi_*^{2+\gamma\widehat\mu_*} ,\ 
\mu_*=1+(1-\theta_0)\kappa\text{ and }\widehat\mu_*=1+(1-\widehat\theta_0)\kappa.
\eeqs 
Then it follows from \eqref{beg33} that
\beq\label{beginest3}
\| u\|_{L_\phi^{\kappa\alpha}(U\times (T_2,T))}
\le  C_7^\frac{1}{\kappa\alpha} 
\left [2 (\mathcal S^\frac{\mu_*}\kappa+\mathcal S^\frac{\widehat \mu_*}\kappa)\right]^\frac{1}{\alpha}.
\eeq
In fact, one has from \eqref{tk} and \eqref{C4} that 
\beq \label{mustar}
\mu_*=1+\frac{r_*}2 ,\quad 
\widehat\mu_*=\mu_*+\frac{2\lambda}\alpha=1+\frac{r_*}2+\frac{2\lambda}\alpha>\mu_*.
\eeq 
Using the  definition $\mathcal S$ in \eqref{Sdef}  and 
applying inequality \eqref{ee2} to estimate $(\Upsilon^{\alpha-h_1}+\Upsilon^{\alpha+h_3})^\frac{\mu_*}{\kappa}$ and $(\Upsilon^{\alpha-h_1}+\Upsilon^{\alpha+h_3})^\frac{\widehat\mu_*}{\kappa}$, we obtain
\beq \label{SSk}
\begin{aligned}
 &\mathcal S^\frac{\mu_*}\kappa+\mathcal S^\frac{\widehat \mu_*}\kappa
 =\mathcal M_1^\frac{ \mu_*}{\kappa}(\Upsilon^{\alpha-h_1}+\Upsilon^{\alpha+h_3})^\frac{\mu_*}{\kappa}
 +\mathcal M_1^\frac{\widehat \mu_*}{\kappa}(\Upsilon^{\alpha-h_1}+\Upsilon^{\alpha+h_3})^\frac{\widehat\mu_*}\kappa\\
 &\le 2^\frac{\mu_*}\kappa\mathcal M_1^\frac{ \mu_*}{\kappa} \left(   \Upsilon^{(\alpha-h_1)\frac{\mu_*}\kappa}+\Upsilon^{(\alpha+h_3)\frac{\mu_*}\kappa} \right)
 +2^\frac{\widehat\mu_*}\kappa \mathcal M_1^\frac{\widehat \mu_*}{\kappa}\left(   \Upsilon^{(\alpha-h_1)\frac{\widehat\mu_*}\kappa}+\Upsilon^{(\alpha+h_3)\frac{\widehat\mu_*}\kappa}\right).
\end{aligned}
\eeq 
We  simply use 
$2^\frac{\mu_*}\kappa\mathcal M_1^\frac{ \mu_*}{\kappa}
\le 2^\frac{\widehat\mu_*}\kappa \mathcal M_1^\frac{\widehat \mu_*}{\kappa}
\le 2^{\widehat\mu_*}\mathcal M_1^{\widehat\mu_*}.$
Let 
\beq\label{nu12} 
   \nu_3=(\alpha-h_1)\frac{\mu_*}\kappa\text{ and }
   \nu_4=(\alpha+h_3)\frac{\widehat\mu_*}{\kappa},
\eeq
which are the minimum and maximum among the last four powers of $\Upsilon$ in \eqref{SSk}.
By applying inequality  \eqref{ee4}, we  estimate 
$ \Upsilon^{(\alpha+h_3)\frac{\mu_*}\kappa},  \Upsilon^{(\alpha-h_1)\frac{\widehat\mu_*}\kappa}\le  \Upsilon^{\nu_3}+\Upsilon^{\nu_4}$.
Then one has
\beq \label{Ssum}
 \mathcal S^\frac{\mu_*}\kappa+\mathcal S^\frac{\widehat \mu_*}\kappa
 \le 2^{\widehat\mu_*}\mathcal M_1^{\widehat\mu_*}\cdot 3  (\Upsilon^{\nu_3}+\Upsilon^{\nu_4}).
\eeq 
Since $C_7\ge 1$, we have $C_7^\frac{1}{\kappa\alpha} \le C_7^\frac{1}{\alpha}$. Combining this fact with \eqref{beginest3} and \eqref{Ssum} yields
\beq\label{newJM}
\| u\|_{L_\phi^{\kappa\alpha}(U\times (T_2,T))}
\le \mathcal M_2^{\frac 1{\alpha}} \Big( \Upsilon^{\nu_3}+\Upsilon^{\nu_4} \Big)^{\frac 1{\alpha}},
\eeq
where
\beq\label{M2d} 
\mathcal M_2=C_7\cdot 2 \cdot (2^{\widehat\mu_*}\mathcal M_1^{\widehat\mu_*} \cdot 3)
=3\cdot 2^{2+\widehat \mu_*}\alpha^{-1+\frac 1{r_1}}\max\{1,C_6\} \chi_*^{2+\gamma\widehat\mu_*} \mathcal M_1^{\widehat \mu_*}.
\eeq

We simplify \eqref{newJM} further by finding upper bounds for $\mathcal M_2$, $\Upsilon^{\nu_3}$ and $ \Upsilon^{\nu_4}$ with simpler dependence on $\alpha$.
By \eqref{kapmax} and \eqref{mustar}, one has
$$ \kappa=\mu_*+\frac{1-\lambda-a}{\alpha}< \mu_*+\frac1\alpha,$$
which yields
\beq\label{kmus} 
\frac{\kappa}{\mu_*}<1+\frac{1}{\alpha \mu_*}=1+\frac{1}{\alpha(1+r_*/2)}.
\eeq 
Moreover, with the relation between $\mu_*$ and $\widehat\mu_*$ in \eqref{mustar}, one also has
$$\kappa=\widehat\mu_*-\frac{2\lambda}{\alpha}+\frac{1-\lambda-a}{\alpha}
> \widehat\mu_*-\frac{3\lambda}\alpha,$$ 
which implies
\beq\label{khatmus}
\frac{\widehat\mu_*}{\kappa}< 1+\frac{3\lambda}{\kappa\alpha}< 1+\frac{3\lambda}{\alpha}.
\eeq
Using \eqref{kmus} to find a lower bound for $\mu_*/\kappa$ in the formulas of  $\nu_3$,  and 
\eqref{khatmus} to find an upper bound for $\nu_4$ in \eqref{nu12}, and  
noticing also that $\alpha>\lambda+1=h_1$, we obtain 
$$0<\nu_1=\frac{\alpha-h_1}{1+\frac{1}{\alpha(1+r_*/2)}}< \nu_3
<\nu_4< (\alpha+h_3)(1+\frac{3\lambda}{\alpha})=\nu_2.$$
Then, by \eqref{ee4},  
\beq\label{U4}
\Upsilon^{\nu_3}\le \Upsilon^{\nu_1} + \Upsilon^{\nu_2}\text{ and }
\Upsilon^{\nu_4} \le \Upsilon^{\nu_1} + \Upsilon^{\nu_2}.
\eeq 

Next, we estimate $\mathcal M_2$ defined by \eqref{M2d}. 
Using the formula of $C_6$ in \eqref{C4} and the formula of $\mathcal M_1$ in \eqref{mathcalM}, we explicitly have
\beq \label{M2pes}
\begin{aligned}
    \mathcal M_2
    & \le 3\cdot2^{3+\widehat\mu_*+2\lambda+\frac{1}{r_1}}\alpha^{-1+\frac 1{r_1}}\max\{1,c_7^{2-a}\mathcal N_2\}\left (2\max\{c_8,c_9\}\alpha^{1+\frac{2-a}{1-a}}\mathcal M_0\right)^{\widehat\mu_*} \chi_*^{2+\gamma \widehat\mu_*}\\
   & \le  3\cdot 2^{3+2\widehat\mu_*+2\lambda+\frac{1}{r_1}}\cdot\left(\max\{1,  c_7^{2-a}\} \max\{1, \mathcal N_2\}\right) \cdot (\max\{c_8,c_9\})^{\widehat\mu_*}\\ 
   &\quad\cdot 
   \alpha^{-1+\frac 1{r_1}+(2+\frac{1}{1-a})\widehat\mu_*}\mathcal M_0^{\widehat\mu_*}\chi_*^{2+\gamma\widehat\mu_*}.
\end{aligned}
\eeq 
Combining the last formula of $\widehat\mu_*$ in \eqref{mustar} with \eqref{lamin} gives
\beq\label{barmu}
\widehat\mu_*<1+\frac{r_*}2+\min\left\{1-r_1,\frac{2\lambda}{\lambda+1}\right\}=\bar\mu.
\eeq
By the fact $0<r_*<1$, we further have
\beq\label{mup}
\widehat\mu_*<1+\frac12+(1-r_1)<\frac52.
\eeq
We use \eqref{barmu} to estimate the powers of $\mathcal M_0$ and $\chi_*$ in \eqref{M2pes}, and use \eqref{mup} for the other $\widehat\mu_*$ appearing there. Also, from the last inequality in \eqref{newro}, $1/r_1<2-a<2$.
Therefore, we  obtain
\begin{multline}\label{lastM1}
\mathcal M_2
\le 
 3\cdot 2^{3+2(5/2)+2\lambda+2}\max\{1,  c_7^{2-a}\} (\max\{c_8,c_9\})^{\frac52} \alpha^{-1+2+(2+\frac{1}{1-a})\frac52}\max\{1, \mathcal N_2\}  \mathcal M_0^{\bar\mu}\chi_*^{2+\gamma\bar\mu}\\
\le 
 3\cdot 2^{10+2\lambda}\max\{1,  c_7^{2-a}\} (\max\{c_8,c_9\})^{\frac52} \alpha^{6+\frac{5}{2(1-a)}}\max\{1, \mathcal N_2\}  \mathcal M_0^{\bar\mu}\chi_*^{2+\gamma\bar\mu}
 =\frac{A_\alpha}2.
 \end{multline}
 
Finally, it follows from \eqref{newJM}, \eqref{U4} and \eqref{lastM1} that
\beqs 
\| u\|_{L_\phi^{\kappa\alpha}(U\times (T_2,T))}
\le \left[ 2\mathcal M_2( \Upsilon^{\nu_1}+\Upsilon^{\nu_2})\right]^{\frac 1{\alpha}}
\le \left[ A_\alpha( \Upsilon^{\nu_1}+\Upsilon^{\nu_2})\right]^{\frac 1{\alpha}}
\eeqs 
which proves \eqref{bfinest2}.
\end{proof}

Next, we iterate inequality \eqref{bfinest2} to obtain the $L^\infty$-estimates. 
Though standard, the procedure requires a rigorous treatment when many factors are involved.
Moreover, only with detailed calculations below that we can obtain estimates with explicit  dependence on $\chi_*$ and the initial and boundary data.

The following technical result for sequences is recalled from \cite{CHK1}.

\begin{lemma}[\cite{CHK1}, Lemma A.2]\label{Genn}
Let $y_j\ge 0$, $\kappa_j>0$, $s_j \ge r_j>0$ and $\omega_j\ge 1$  for all $j\ge 0$.
Suppose there is $A\ge 1$ such that
\beqs
y_{j+1}\le A^\frac{\omega_j}{\kappa_j} (y_j^{r_j}+y_j^{s_j})^{\frac 1{\kappa_j}}\quad\forall j\ge 0.
\eeqs
Denote $\beta_j=r_j/\kappa_j$ and $\gamma_j=s_j/\kappa_j$.
Assume
$\bar\alpha \eqdef \sum_{j=0}^{\infty}  \frac{\omega_j}{\kappa_j}<\infty$ 
and the products $\prod_{j=0}^\infty \beta_j$, $\prod_{j=0}^\infty \gamma_j$
 converge to positive numbers  $\bar\beta$, $\bar \gamma$,  respectively.
Then
\beq\label{doub2}
\limsup_{j\to\infty} y_j\le (2A)^{G \bar\alpha} \max\{y_0^{\bar \beta},y_0^{\bar \gamma}\},
\text{ where  }G=\limsup_{j\to\infty} G_j\in(0,\infty)
\eeq
with $G_j=\max\{1, \gamma_{m}\gamma_{m+1}\ldots\gamma_{n}:1\le m\le n< j\}$.
\end{lemma}

In the case $\gamma_j\ge 1$ for all $j\ge 0$, one has
$G_j=\gamma_1\gamma_2\ldots\gamma_{j-1}$ for all $j\ge 2$, hence, the number $G$ in \eqref{doub2} actually is
\beq\label{Gprod}
G=\prod_{j=1}^\infty \gamma_j.
\eeq

\begin{assumption}\label{asmp44}
We assume the following.
\begin{enumerate}[label=\tnum]
    \item\label{a2} Let  $\widetilde\kappa$ be a fixed number in the interval $\left(1,\sqrt{1+r_*/2}\right)$, and the numbers $p_i>1$, for $i=1,2,\ldots,6$, be as in  \eqref{ktil}.
    \item\label{a3} The integrals in \eqref{wcond1}, \eqref{Psicond}, \eqref{wcond2} and $\Phi_{**}\eqdef \int_U \phi^{-1}(x) \d x$ are finite.
\end{enumerate}
\end{assumption}

Recalling  $\mathcal N_1$ and $\mathcal N_2$ from \eqref{N1} and \eqref{GG}, respectively, we denote 
\beq \label{N3}
\mathcal N_3=\max\{\mathcal N_1,\mathcal N_2\}\ge 1.
\eeq 
Then $\mathcal N_3$ involves the weight functions $\phi(x)$ and $W_k(x)$, for $0\le k\le 4$, while $\Psi_T$ in \eqref{bcquant} involves the boundary data $\psi_1(x,t)$ and $\psi_2(x,t)$.

\begin{theorem}\label{LinfU} 
Under Assumption \ref{asmp44}, let $\alpha_0$ be a number such that 

 \beq\label{alp0}
 \begin{split}
\alpha_0>\max\left\{1+\lambda, \frac{p_2(\lambda(5-4a)-1)}{\widetilde \kappa-p_2}, \frac{p_5(a+3\lambda-1)}{(1-a)(\widetilde\kappa -p_6)},
\frac{1-\lambda-a}{\widetilde\kappa-1},\frac{2\lambda}{1-r_1}, \frac{\lambda+a-1}{1+r_*/2-\widetilde\kappa^2} \right\}.
\end{split}
\eeq
Then there are positive  constants $\widehat C_0$, $\widetilde  \mu$, $\widetilde \nu$, $\omega_0 $, $\omega_1 $, $\omega_2$, $\omega_3$ independent of $\chi_*$, $\phi(x)$, $W_j(x)$ for $1\le j\le 4$, $u_0(x)$, $\psi_1(x,t)$ and $\psi_2(x,t)$ such that
if  $T>0$ and $\sigma\in (0,1)$ then
\beq\label{Li1}
\|u\|_{L^{\infty}(U\times(\sigma T,T))} 
\le \widehat C_0 \chi_*^{\omega_0}(1+T)^{\omega_1}\Big(1+\frac1{\sigma T}\Big)^{\omega_2}
\mathcal N_3^{\omega_3} \Psi_T^{\omega_2}\max\left\{ \|u\|^{\widetilde \mu}_{L_\phi^{\widetilde \kappa \alpha_0}(Q_T)},\|u\|^{\widetilde \nu}_{L_\phi^{\widetilde\kappa \alpha_0}(Q_T)}\right\}.
\eeq
\end{theorem}
\begin{proof}
In order to iterate inequality \eqref{bfinest2} in Lemma \ref{GLk}, we set, for any integer $j\ge 0$,  the numbers 
\beq\label{btj} 
\beta_j=\widetilde \kappa^j\alpha_0\text{ and }t_j=\sigma T(1-2^{-j}).
\eeq 
Since $\widetilde \kappa>1$, one has the sequence $(\beta_j)_{j=0}^\infty$ is increasing.
Note from \eqref{alp0} and, thanks to Assumption \ref{asmp44}\ref{a2}, the fact $1+r_*/2-\widetilde\kappa^2>0$  that 
\beq\label{ared}
\alpha_0>\begin{cases}
    0\ge\frac{2(\lambda+a-1)}{r_*}&\text{if } \lambda+a-1\le0,\\
    \frac{\lambda+a-1}{1+r_*/2-\widetilde\kappa^2}>\frac{\lambda+a-1}{r_*/2}=\frac{2(\lambda+a-1)}{r_*}& \text{otherwise.}
\end{cases}
\eeq
For $j\ge 0$, because $\beta_j\ge \alpha_0$  and by assumption \eqref{alp0} and its consequence \eqref{ared}, the number $\alpha=\beta_j$  satisfies \eqref{alp-Large} and \eqref{Bigalp}.
We  then can apply estimate \eqref{bfinest2} to $\alpha=\beta_{j}$, $T_2=t_{j+1}$ and $T_1=t_j$ to obtain
\begin{align}
\| u\|_{L_\phi^{\kappa(\beta_{j})\beta_{j}}(U\times (t_{j+1},T))}
&\le   A_{\beta_{j}}^\frac1{\beta_{j}}\Big( \| u\|_{L_\phi^{\beta_{j+1}}(U\times(t_j,T))}^{\widetilde r_j}+\| u\|_{L_\phi^{\beta_{j+1}}(U\times(t_j,T))}^{\widetilde s_j}\Big)^\frac 1{\beta_{j}},\label{ukbb}
\end{align}
where $A_{\beta_j}$ is defined by \eqref{Atildef2} with $\alpha=\beta_j$, and, referring to \eqref{tilrsdef}, 
\beq\label{trsj}
\widetilde r_j=\nu_1(\beta_{j})\quad \text{and}\quad\widetilde s_j=\nu_2(\beta_{j}).
\eeq

Observe from \eqref{btj}  and \eqref{alp0} that
$$\beta_j\ge \beta_0=\alpha_0>\frac{\lambda+a-1}{1+r_*/2-\widetilde\kappa^2} .$$
Thus, referring to formula \eqref{kapmax} with $\alpha=\beta_j$, one has 
\beqs
\kappa(\beta_j)=1+\frac{r_*}{2}+\frac{1-\lambda-a}{\beta_j}
>\widetilde \kappa^2
\eeqs
which implies
\beq \label{kbb}
\kappa(\beta_{j})\beta_{j}> \widetilde \kappa^2 \beta_{j}=\beta_{j+2}.
\eeq 

For $j\ge 0$, define $\mathcal Q_j=U\times (t_j,T)$ and $Y_j=\| u\|_{L_\phi^{\beta_{j+1}}(\mathcal Q_j)}.$
Thanks to \eqref{kbb}, one can apply H\"older's inequality with the powers $\kappa(\beta_{j})\beta_{j}/\beta_{j+2}$ and its H\"older conjugate to have
\begin{align*}
 Y_{j+1}
 &= \| u\|_{L_\phi^{\beta_{j+2}}(\mathcal Q_{j+1})}
\le \|u\|_{L_\phi^{\kappa(\beta_{j})\beta_{j}}(\mathcal Q_{j+1})}
\left(T\int_U\phi(x) \ dx\right)^{\frac{1}{\beta_{j+2}}-\frac{1}{\kappa(\beta_{j})\beta_{j}}}
\\
 &\le [(1+T)\Phi_*]^{\frac{1}{\beta_{j+2}}-\frac{1}{\kappa(\beta_{j})\beta_{j}}}\|u\|_{L_\phi^{\kappa(\beta_{j})\beta_{j}}(\mathcal Q_{j+1})}
 \le [(1+T)\Phi_*]^{\frac{1}{\beta_{j}}}\|u\|_{L_\phi^{\kappa(\beta_{j})\beta_{j}}(\mathcal Q_{j+1})},
\end{align*}
where $\Phi_*\ge 1$ is defined by \eqref{Phistar}.
Combining this inequality with \eqref{ukbb}, we obtain 
\beq\label{YwithQ}
Y_{j+1}\le \widehat A_j^\frac{1}{\beta_{j}}\big( Y_j^{\widetilde{r}_j}+Y_j^{\widetilde{s}_j}\big)^{\frac 1{\beta_{j}}} \text{ for all }j\ge 0, \text{ where }
\widehat A_j=  (1+T)\Phi_* A_{\beta_{j}}.
\eeq

We estimate $\widehat A_j$ next.
With $\mathcal N_3$ in \eqref{N3}, we estimate $\Phi_*\le \mathcal N_1\le \mathcal N_3$ and $\max\{1,\mathcal N_2\}\le \mathcal N_3$.
Combining these with  \eqref{Atildef2}  and \eqref{YwithQ} yields 
\begin{align*}
\widehat A_j
&\le  (1+T)\mathcal N_3 \cdot  c_{10}  \beta_j^{6+\frac{5}{2(1-a)}} 
\chi_*^{2+\gamma\bar\mu}\mathcal N_3\left[(1+T)\Big(1+\frac{2^{j+1}}{\sigma T}\Big)\mathcal N_3 \Psi_T\right]^{\bar\mu} \\
&\le c_{10}  \chi_*^{2+\gamma\bar\mu} (1+T)^{1+\bar\mu} \left[2^{j+1}\Big(1+\frac{1}{\sigma T} \Big)\right]^{\bar\mu}  \left( \widetilde \kappa^j\alpha_0\right)^{6+\frac{5}{2(1-a)}} \mathcal N_3^{2+\bar\mu}
\Psi_T^{\bar\mu}.
\end{align*}
Simply using the facts that all the numbers $c_{10}$, $\widetilde\kappa$, $\alpha_0$, $\chi_*$, $\mathcal N_3$, $\Psi_T$   are greater or equal to $1$, 
$j+1\ge 1$ and $\widetilde \kappa^j\le \widetilde \kappa^{j+1}$,
we obtain
$\widehat A_j
\le A_{T,\sigma,\alpha_0}^{j+1}$,
where
\beq \label{AAA}
 A_{T,\sigma,\alpha_0}=c_{10} 2^{\bar\mu} \chi_*^{2+\gamma\bar\mu} (1+T)^{1+\bar\mu} \Big(1+\frac{1}{\sigma T} \Big)^{\bar\mu}  \left( \widetilde \kappa\alpha_0\right)^{6+\frac{5}{2(1-a)}} \mathcal N_3^{2+\bar\mu}
\Psi_T^{\bar\mu}.
 \eeq 
 Thus,
\beq\label{YY}
Y_{j+1}\le A_{T,\sigma,\alpha_0}^\frac{j+1}{\beta_{j}}\big( Y_j^{\widetilde{r}_j}+Y_j^{\widetilde{s}_j}\big)^{\frac 1{\beta_{j}}} \text{for all }j\ge 0.
\eeq

With \eqref{YY}, we apply Lemma \ref{Genn} to the sequence $(Y_j)_{j=0}^\infty$.
We check the requirements there.
We have
\beqs
L_0\eqdef \sum_{j=0}^\infty \frac {j+1} { \beta_j}= \frac 1{ \alpha_0}\sum_{j=0}^\infty \frac {j+1} {\widetilde \kappa^j}=\frac{1}{\alpha_0(1-\widetilde \kappa^{-1})^2}\in(0,\infty).
\eeqs
Also, by \eqref{trsj} and \eqref{tilrsdef} with $\alpha=\beta_j$, one has 
\beqs
\sum_{j=0}^\infty\ln\frac{\widetilde r_j}{\beta_j}=\sum_{j=0}^\infty\ln\left(1 - \frac{h_1}{\widetilde \kappa^j \alpha_0}\right)-\sum_{j=0}^\infty\ln\left (1+\frac{1}{\widetilde\kappa^j \alpha_0(1+ r_*/2)}\right)=L_1\in\R,
\eeqs
\beqs 
\sum_{j=0}^\infty\ln{\frac{\widetilde s_j}{\beta_j}}
=\sum_{j=0}^\infty\ln\left(1+\frac{h_3}{\widetilde\kappa^j\alpha_0}\right)
+\sum_{j=0}^\infty\ln \left(1+\frac{3\lambda}{\widetilde\kappa^j\alpha_0}\right)
=L_2\in\R.
\eeqs 
Thus,  
\beq \label{mnutil}
\widetilde \mu \eqdef   \prod_{j=0}^\infty \frac{\widetilde r_j}{\beta_j}=e^{L_1}
\text{ and }
\widetilde \nu \eqdef \prod_{j=0}^\infty \frac{\widetilde s_j}{\beta_j}=e^{L_2}
 \text{ are  positive numbers.}
\eeq 
Then, thanks to \eqref{YY}, \eqref{mnutil} and Lemma \ref{Genn} with the particular case \eqref{Gprod}, we obtain
\beq\label{limY}
\limsup_{j\to\infty} Y_{j}\le (2A_{T,\sigma,\alpha_0})^\omega\max\{Y_0^{\widetilde\mu}, Y_0^{\widetilde\nu}\},
\eeq
where 
\beqs 
\omega=\mathcal G\sum_{j=0}^{\infty}\frac {j+1}{\beta_j}=\mathcal G L_0\text{ with  }
\mathcal G=\prod_{k=1}^\infty \frac{\widetilde s_k}{\beta_k}\in(0,\infty).
\eeqs 

For $j\ge 1$, by applying H\"older's inequality with the powers $2$ and $2$, and then using the fact $t_{j-1}<\sigma T$, we have
\beq\label{bbeta2}
\begin{aligned}
&\|u\|_{L^{\beta_{j}/2}(U\times(\sigma T,T))}= \left(\int_{\sigma T}^T \int_U \left(|u(x,t)|^{\beta_j/2}\phi(x)^{1/2}\right)\cdot\phi(x)^{-1/2} \d x \d t\right)^{2/\beta_j} \\
&\le \left(\int_{\sigma T}^T \int_U |u(x,t)|^{\beta_j} \phi(x) \d x \d t\right)^{1/\beta_j}
\left (\int_{\sigma T}^T \int_U \phi(x)^{-1} \d x \d t\right)^{1/\beta_j}\le Y_{j-1}
\left (T\Phi_{**}\right)^{1/\beta_j},
\end{aligned}
\eeq
where $\Phi_{**}$ is the finite number from part \ref{a3} of Assumption \ref{asmp44}.
By taking the limit superior, as $j\to\infty$, of \eqref{bbeta2} and 
 utilizing \eqref{limY},  we obtain 
 \beq \label{limbeta} 
\begin{aligned}
&\|u\|_{L^\infty(U\times(\sigma T,T))}
=\lim_{j\to\infty}\|u\|_{L^{\beta_{j}/2}(U\times(\sigma T,T))} 
\le \limsup_{j\to\infty}Y_{j-1} \\
&\le (2A_{T,\sigma,\alpha_0})^\omega\max\{Y_0^{\widetilde\mu}, Y_0^{\widetilde\nu}\}
=  (2A_{T,\sigma,\alpha_0})^\omega\max\Big\{ \|u\|^{\widetilde \mu}_{L_\phi^{\widetilde \kappa \alpha_0}(Q_T)},\|u\|^{\widetilde \nu}_{L_\phi^{\widetilde\kappa \alpha_0}(Q_T)}\Big\}
.
\end{aligned}
\eeq 
With $A_{T,\sigma,\alpha_0}$ in \eqref{AAA}, we have
\beq\label{AT2}
(2A_{T,\sigma,\alpha_0})^\omega
\le \widehat C_0 \chi_*^{\omega_0}(1+T)^{\omega_1} \Big(1+\frac1{\sigma T}\Big)^{\omega_2} 
\mathcal N_3^{\omega_3}  \Psi_T^{\omega_2},
\eeq
where 
$\widehat C_0=\left[2^{1+\bar\mu}c_{10} ( \widetilde \kappa \alpha_0)^{6+\frac{5}{2(1-a)}}\right]^\omega$, and 
\beq \label{omo}
\omega_0=(2+\gamma\bar\mu)\omega,\quad 
\omega_1=(1+\bar\mu)\omega,\quad 
\omega_2=\bar\mu\omega,\quad 
\omega_3=(2+\bar\mu)\omega.
\eeq 
Then the desired estimate \eqref{Li1} follows from \eqref{limbeta}, \eqref{AT2} and \eqref{omo}.
\end{proof}

In the following, Theorem \ref{LinfU} is used in conjunction with Theorem \ref{Labound} to obtain  the $L^\infty$-estimates for the solution $u(x,t)$ in terms of the initial  and boundary data.

\begin{theorem}\label{thm45}
Under Assumption \ref{asmp44}, let $r$ be a number that satisfies \eqref{newrr}, $\widetilde r$ be defined by \eqref{newtr},  and  $\alpha_0$ be a number that satisfies 
\beq\label{alast}
\begin{aligned}
&\alpha_0 >\max  \left \{  \frac{1-\lambda-a}{\widetilde \kappa -1} ,  \frac{p_2(\lambda(5-4a)-1)}{\widetilde \kappa-p_2}, 
\frac{p_5(a+3\lambda-1)}{(1-a)(\widetilde\kappa -p_6)},\frac{2\lambda}{1-r_1}, \frac{\lambda+a-1}{1+r_*/2-\widetilde\kappa^2},\right.\\  
&\left.  \frac{2(2-a)(r+a+\lambda-1)}{\widetilde\kappa r_* (1-a)}, 
\frac{1}{\widetilde \kappa}\left[\frac{2\lambda(2+r_*)}{(1-a)r_*}+\frac{2}{r_*}\right],
\frac{1}{\widetilde \kappa}\left[ 2\max\{r,\widetilde r\}\left(1+\frac{2}{r_*}\right)+\frac{4\lambda}{r_*}\right]\right\}.
\end{aligned}
\eeq
Let $\widetilde  \mu$, $\widetilde \nu$, $\omega_0$, $\omega_1$, $\omega_2$, $\omega_3$ be the same constants as in Theorem \ref{LinfU}. Denote $\beta_1=\widetilde \kappa\alpha_0$.
Assume all ${\mathcal K}_j$, for $0\le j\le 6$, from \eqref{Kz} to \eqref{Kf} with $\alpha=\beta_1$ are finite numbers.
Let the numbers $\mu_*$, $\gamma_*$ and $\bar C_{\beta_1}$ be as in Theorem \ref{Labound} for $\alpha=\beta_1$.
 Let 
\beqs 
V_{0}=1+\int_U \phi(x) u_0(x)^{\beta_1}\d x\text{ and }
 M(t) =1+\int_\Gamma \left[ (\psi_1^-(x,t))^{\frac{\beta_1+r}{r+\lambda}}+ (\psi_2^-(x,t))^{\frac{\beta_1+r}{r}} \right]\d S.
\eeqs
Assume $T>0$ satisfies 
\beq\label{Tsm0}
\int_0^T M(\tau) \d\tau< \frac{\beta_1}{2\mu_* \chi_*^{\gamma_*}\bar C_{\beta_1}  } V_{0}^{-\mu_*/\beta_1}  .
\eeq
Let $\varep$ be any positive number satisfying $\varepsilon<\min\{1,T\}$.
Then one has
\beq\label{Lb1}
\|u\|_{L^{\infty}(U\times(\varepsilon,T))}
\le \widehat C_1\chi_*^{\omega_0}(1+T)^{\omega_1}\varepsilon^{-\omega_2}  
V_0^{\frac{\widetilde \nu}{\beta_1}} \Psi_T^{\omega_2}
\max\left\{ \|\delta^{-\frac{1}{\mu_*}}\|_{L^{\beta_1}(0,T)}^{\widetilde \mu}, 
\|\delta^{-\frac{1}{\mu_*}}\|_{L^{\beta_1}(0,T)}^{\widetilde \nu}\right\}, 
\eeq
where $\widehat C_1$ is a positive constant depending on the $\mathcal K_j$ and $\mathcal N_3$, but independent of $T$, $\varep$, $\chi_*$, $u_0(x)$, $\psi_1(x,t)$ and $\psi_2(x,t)$,
and 
\beq\label{newVt}
\delta(t)=1-2\mu_*\chi_*^{\gamma_*}\beta_1^{-1}\bar C_{\beta_1}V_0^{\mu_*/\beta_1}\int_0^t M(\tau) \d \tau \in(0,1) \text{ for }t\in[0,T].
\eeq
Consequently,  one has  
\beq\label{Lb2}
\|u\|_{L^{\infty}(U\times(\varepsilon,T))}\le \widehat C_1 \chi_*^{\omega_0}\varepsilon^{-\omega_2}(1+T)^{\omega_1 + \widetilde \nu/\beta_1} \delta(T)^{-\frac{\widetilde \nu}{\mu_*}}   
\Big(1+\norm{u_0}_{L_\phi^{\beta_1}(U)}\Big)^{\widetilde \nu}
\Psi_T^{\omega_2}  .
\eeq
\end{theorem}
\begin{proof}
We apply Theorem \ref{LinfU} to $\sigma T=\varepsilon<1$. Since $\widetilde \kappa\in (1,\sqrt{1+r_*/2})$ and $r_*\in (0,1)$, we have $\widetilde\kappa<2$, hence the 7th number on the right-hand side of \eqref{alast} can be bounded from below by
 \beq\label{kaplam}
 \frac{1}{\widetilde \kappa}\left[\frac{2\lambda(2+r_*)}{(1-a)r_*}+\frac{2}{r_*}\right]>  \frac{1}{2}(4\lambda+2) >\lambda +1.
 \eeq
Using assumption \eqref{alast} with the numbers on its right-hand side from the 1st to the 5th places and then the 7th place together with estimate \eqref{kaplam}, we have that $\alpha_0$ satisfies \eqref{alp0}. We also use the simple estimate
$(1+\frac1{\sigma T})^{\omega_2}=(1+\varep^{-1})^{\omega_2}\le (2\varep^{-1})^{\omega_2}$.
Then it follows  estimate \eqref{Li1} of Theorem \ref{LinfU}  that
\beq\label{Li2}
\|u\|_{L^{\infty}(U\times(\varepsilon,T))}
\le 2^{\omega_2}\widehat C_0 \chi_*^{\omega_0} (1+T)^{\omega_1} \varepsilon^{-\omega_2} \mathcal N_3^{\omega_3}  \Psi_T^{\omega_2} \max\left\{ \|u\|^{\widetilde \mu}_{L_\phi^{\beta_1}(Q_T)},\|u\|^{\widetilde \nu}_{L_\phi^{\beta_1}(Q_T)}\right\}.
\eeq

We estimate the last two norms in \eqref{Li2} further by applying  Theorem \ref{Labound} to $\alpha=\beta_1$.
The requirement $\beta_1=\widetilde\kappa\alpha_0>\alpha_*$ where $\alpha_*$ is defined by  \eqref{alstar}, is equivalent to 
\beq\label{bka}
\begin{aligned}
\alpha_0>\frac{\alpha_*}{\widetilde\kappa}= \frac1{\widetilde\kappa}\max & 
\left\{  \frac{2\lambda r_1}{1-r_1}, 
\frac{2(2-a)(r+a+\lambda-1)}{r_*(1-a)},\frac{2\lambda(2+r_*)}{(1-a)r_*}+\frac{2}{r_*}, \right.\\
&\quad \left. 
2r\left(1+\frac2{r_*}\right)+\frac{4\lambda}{r_*},
2\widetilde r\left(1+\frac2{r_*}\right)+\frac{4\lambda}{r_*}\right\}.
\end{aligned}
\eeq
Since $r_1\in(0,1) $ and $\widetilde\kappa>1$, we have 
\beqs 
\frac{2\lambda r_1}{\widetilde\kappa(1-r_1)}<\frac{2\lambda}{1-r_1}.
\eeqs 
Hence, assumption \eqref{alast} with the particular numbers on the 4th, 6th--8th places on its right-hand side implies \eqref{bka}. Moreover, thanks to \eqref{Tsm0}, the condition \eqref{tsmall1} is also satisfied for $\alpha=\beta_1$. 
Consequently, using estimate \eqref{uwest} from Theorem \ref{Labound} with $\alpha=\beta_1$ yields
\beq\label{Lb}
\|u\|_{L_\phi^{\beta_1}(Q_T)}
=\left(\int_0^T \|u(\cdot,t)\|_{L_\phi^{\beta_1}(U)}^{\beta_1} \d t\right)^{\frac{1}{\beta_1}}
\le \left(\int_0^T V_0 \delta(t)^{-\frac{\beta_1}{\mu_*}}\d t\right)^{\frac{1}{\beta_1}} 
=V_0^{\frac{1}{\beta_1}} \|\delta^{-\frac{1}{\mu_*}}\|_{L^{\beta_1}(0,T)}.
\eeq 

Since  $V_0\ge 1$ and $\widetilde\mu<\widetilde\nu$, which is due to \eqref{tilrsdef}, \eqref{trsj} and \eqref{mnutil},   we can estimate $V_0^\frac{\widetilde\mu}{\beta_1}\le V_0^\frac{\widetilde\nu}{\beta_1}$.
Combining this with  \eqref{Li2} and \eqref{Lb}, we obtain inequality \eqref{Lb1} with $\widehat C_1=2^{\omega_2}\widehat C_0 \mathcal N_3^{\omega_3}$.

Observing that  $\beta_1>\alpha_0>1$, see \eqref{kaplam}, we apply inequality \eqref{ee3} to $p=1/\beta_1<1$
to have 
\beq\label{Voe}
V_0^\frac{1}{\beta_1}\le  1+\norm{u_0}_{L_\phi^{\beta_1}(U)}.
\eeq 
Thanks to condition \eqref{Tsm0} again, the number $\delta(t)$ in \eqref{newVt} belongs to the interval $(0,1)$ and is decreasing in $[0,T]$. 
Then  one has
\beq\label{Lb3}
\|\delta^{-\frac{1}{\mu_*}}\|_{L^{\beta_1}(0,T)}
\le \delta(T)^{-\frac{1}{\mu_*}}T^\frac{1}{\beta_1}
\le \delta(T)^{-\frac{1}{\mu_*}}(1+T)^\frac{1}{\beta_1}.
\eeq
Note that the last upper bound satisfies $\delta(T)^{-\frac{1}{\mu_*}}(1+T)^\frac{1}{\beta_1}\ge 1$.
Combining \eqref{Lb1} with \eqref{Voe}, \eqref{Lb3} and the fact $\widetilde\mu<\widetilde \nu$ yields 
\beqs
\|u\|_{L^{\infty}(U\times(\varepsilon,T))}
\le \widehat C_1\chi_*^{\omega_0} (1+T)^{\omega_1} \varepsilon^{-\omega_2} 
\left (1+\norm{u_0}_{L_\phi^{\beta_1}(U)}\right)^{\widetilde \nu} \Psi_T^{\omega_2} \left[\delta(T)^{-\frac{1}{\mu_*}}(1+T)^\frac{1}{\beta_1}\right]^{\widetilde \nu} 
\eeqs
which implies the desired estimate \eqref{Lb2}.
\end{proof}

\appendix
\section{}\label{Appex}

This Appendix contains the proofs for Lemmas \ref{RoS1}, \ref{Rtrace} and \ref{RoS2} of Section \ref{Prem}.

\begin{proof}[Proof of Lemma \ref{RoS1}]
Denote
\beqs
    I_*=\int_U \frac{|u|^{\alpha-s}}{(1+|u|)^\beta}|\nabla u|^p W_*\d x  
\text{ and }    J_0=\int_U (1+|u|)^{\alpha}\varphi  \d x.
\eeqs
For a relation between $J_*$ and $J_0$, we  note, with the use of inequality \eqref{ee2} to estimate $(1+|u|)^\alpha$,  that
\beq\label{Joest}
  J_0  \le  \int_U 2^\alpha(\varphi + |u|^{\alpha}\varphi)  \d x
  \le 2^\alpha {\mathcal E}_1(1+J_*). 
\eeq
Moreover, for any number $z>0$, raising  inequality \eqref{Joest} to the power $z/\alpha$ and then applying inequality  \eqref{ee2}, we have
\beq\label{group2u} 
  J_0^{\frac{z}{\alpha}}
\le (2^\alpha {\mathcal E}_1)^{\frac{z}{\alpha}} (1+J_*)^{\frac{z}{\alpha}}\le (2^\alpha {\mathcal E}_1)^{\frac{z}{\alpha}} 2^{\frac{z}{\alpha}}\left(1+J_*^{\frac{z}{\alpha}}\right)\\
= 2^{z(1+\frac1\alpha)}  
{\mathcal E}_1^{\frac{z}{\alpha}}\left(1+J_*^{\frac{z}{\alpha}}\right). 
\eeq

We apply Lemma \ref{WS1} with the following  weight function
\beq\label{subW}
W(x)=(1+|u(x)|)^{-\beta}W_*(x). 
\eeq 
Then $G_2$ in \eqref{Gdef} becomes
\beq\label{newG2}
     G_2=\left(\int_U (1+|u|)^{\frac{r_1\beta}{1-r_1}}W_*^{-\frac{r_1}{1-r_1}}\d x\right)^{\frac{1-r_1}{r_1}}.
\eeq
It follows inequality \eqref{S10} that
\beq \label{SSW}
\int_U |u|^{\alpha+r} \omega \d x 
 \le \varep I_*+D_{1,m,\theta} \Phi_1 J_*^{1+r/\alpha}
 +\varep^{-\frac\theta{1-\theta}}D_{2,m,\theta} \Phi_2 J_*^{1+\mu_1/\alpha},
\eeq 
where $\Phi_1$ and $\Phi_2$ are defined in \eqref{Gdef} with the new $G_2$ in \eqref{newG2}.

Under condition \eqref{abeta}, we estimate, by using H\"older's inequality with powers $\frac{\alpha(1-r_1)}{r_1\beta}$ and $\frac{\alpha(1-r_1)}{\alpha(1-r_1)-r_1\beta}$, 
\beq\label{G2est}
G_2=\left(\int_U \left\{ (1+|u|)^{\frac{r_1\beta}{1-r_1}}\varphi^{\frac{r_1\beta}{\alpha(1-r_1)}}\right\}
     \cdot 
     \left\{ \varphi^{-\frac{r_1\beta}{\alpha(1-r_1)}}W_*^{-\frac{r_1}{1-r_1}}\right\} 
     \d x\right)^{\frac{1-r_1}{r_1}}
     \le  J_0^{\beta/\alpha} {\widetilde G}_2.
\eeq
Letting $z=\frac{\beta\theta}{1-\theta}$ in \eqref{group2u}, we have
\beq\label{group1u} 
  J_0^{\frac{\beta\theta}{\alpha(1-\theta)}}
\le 2^{\widehat\beta_1}  
{\mathcal E}_1^{\frac{\beta\theta}{\alpha(1-\theta)}}\left(1+J_*^{\frac{\beta\theta}{\alpha(1-\theta)}}\right). 
\eeq
It follows from  \eqref{G2est} and \eqref{group1u} that
\beq\label{G22}
G_2^\frac\theta{1-\theta}\le J_0^{\frac{\beta\theta}{\alpha(1-\theta)}}
{\widetilde G}_2^\frac\theta{1-\theta}
\le 2^{\widehat\beta_1}  
{\mathcal E}_1^{\frac{\beta\theta}{\alpha(1-\theta)}}{\widetilde G}_2^\frac\theta{1-\theta}\left(1+J_*^{\frac{\beta\theta}{\alpha(1-\theta)}}\right).
\eeq
Utilizing the estimate \eqref{G22} in the formula of $\Phi_2$ in \eqref{Gdef}, we obtain 
\beq\label{Ph22}
\Phi_2=G_2^\frac{\theta}{1-\theta}G_3
\le 2^{\widehat\beta_1} {\mathcal E}_1^{\frac{\beta\theta}{\alpha(1-\theta)}} {\widetilde G}_2^\frac\theta{1-\theta}G_3 
\left(1+J_*^{\frac{\beta\theta}{\alpha(1-\theta)}}\right)
=2^{\widehat\beta_1}  
\widetilde\Phi_2 \left(1+J_*^{\frac{\beta\theta}{\alpha(1-\theta)}}\right).
\eeq
Combining inequalities \eqref{SSW} and \eqref{Ph22} gives
\begin{align}
\notag 
&\int_U |u|^{\alpha+r} \omega \d x  \le \varep I_* 
  +  D_{1,m,\theta} \Phi_1 J_*^{1+r/\alpha} +\varep^{-\frac\theta{1-\theta}}D_{2,m,\theta}  \cdot 2^{\widehat\beta_1} 
\widetilde\Phi_2 \left(1+J_*^{\frac{\beta\theta}{\alpha(1-\theta)}}\right) \cdot  J_*^{1+\mu_1/\alpha} \\
&=  \varep I_* 
  +D_{1,m,\theta} \Phi_1 J_*^{1+r/\alpha}  
  + 2^{\widehat\beta_1}  
\varep^{-\frac\theta{1-\theta}}D_{2,m,\theta} \widetilde\Phi_2  
\left(J_*^{1+\mu_1/\alpha}
+ J_*^{1+\widehat\mu_1/\alpha}\right).\label{RSi1}
\end{align}
Denote by $q$ any of the last three powers of $J_*$ in  \eqref{RSi1}.
With $\widehat\mu_1>\mu_1$, we have $1+\widehat\mu_2/\alpha\le q\le 1+\widehat\mu_3/\alpha.$
Applying inequality  \eqref{ee4} gives
\beq\label{JqJJ}
J_*^q\le J^{1+\widehat\mu_2/\alpha}+J_*^{1+\widehat\mu_3/\alpha}.
\eeq
Utilizing this inequality to estimate the last three  $J_*$-terms in \eqref{RSi1}, we obtain \eqref{RSi2} with the first part of the minimum.
Now, simply using $J_*^q\le (1+J_*)^{1+\widehat\mu_3/\alpha}$ in place of \eqref{JqJJ},
we  obtain \eqref{RSi2} again with the second part of the minimum.
\end{proof}
\begin{proof}[Proof of Lemma \ref{Rtrace}]
We use the same notation as in the proof of Lemma \ref{RoS1}.
Applying inequality \eqref{trace1} in Lemma \ref{trace} to the function $W(x)$ in \eqref{subW} results in
\beq \label{mainTrace}
\begin{aligned}
\int_\Gamma |u|^{\alpha+r} \d S 
&\le 3\varepsilon I_* 
+ z_1 \Phi_3 J_*^{1+r/\alpha}  + \varep^{-\frac\theta{1-\theta}}z_2 \Phi_4 J_*^{1+\mu_1/\alpha}
+ \varepsilon^{-\frac 1{p-1}}z_4 \Phi_6
  J_*^{1+\widetilde r/\alpha} \\
&\quad + \varep^{-(\frac 1 {p-1}+\frac p{p-1}\cdot \frac{\widetilde\theta}{1-\widetilde\theta} )}z_5 \Phi_7
  J_*^{1+\widetilde\mu_1/\alpha}.
\end{aligned}
\eeq 
We need to estimate $\Phi_4$, $\Phi_6$ and $\Phi_7$ appearing in \eqref{mainTrace}.

Note from the formula of $\Phi_4$ in \eqref{Ph34} and estimate \eqref{G22} of $G_2^\frac\theta{1-\theta}$ that
\beq\label{Ph4est}
\Phi_4\le 2^{\widehat\beta_1}  
{\widetilde G}_2^\frac\theta{1-\theta} 
{\mathcal E}_1^{\frac{\beta\theta}{\alpha(1-\theta)}}\left(1+J_*^{\frac{\beta\theta}{\alpha(1-\theta)}}\right)\cdot G_4 \\
=2^{\widehat\beta_1}  
{\widetilde\Phi}_4\left(1+J_*^{\frac{\beta\theta}{\alpha(1-\theta)}}\right). 
\eeq

Regarding $G_5$ in the formula of $\Phi_6$ in \eqref{Ph67},  we have, with the substitution \eqref{subW},
\begin{align*}
G_5&=\left(\int_U  \varphi^{-1}\left\{ (1+|u|)^{-\beta}W_*\right\}^{-\frac{1}{(p-1)(1-\widetilde\theta) (1+\widetilde\mu_1/\alpha)}}\d x\right)^{1+\widetilde\mu_1/\alpha}\\
&=\left(\int_U  \varphi^{-1}(1+|u|)^{\beta_*}W_*^{-\beta_*/\beta}\d x\right)^{1+\widetilde\mu_1/\alpha}. 
\end{align*}
Under assumption \eqref{bcond}, applying H\"older's inequality with powers $\alpha/\beta_*$ and $\alpha/(\alpha-\beta_*)$ yields
\begin{align*}
    G_5 &=\left(\int_U  \left\{(1+|u|)^{\beta_*}\varphi^{\beta_*/\alpha}\right\}\cdot\left\{\varphi^{-1-\beta_*/\alpha} W_*^{-\beta_*/\beta}\right\}\d x\right)^{1+\widetilde\mu_1/\alpha} \\
    &\le \left(\int_U (1+|u|)^\alpha\varphi \d x\right)^{\frac{\beta_*}{\alpha}(1+\widetilde\mu_1/\alpha)}\left(\int_U \varphi^{-\frac{\alpha+\beta_*}{\alpha-\beta_*}} W_*^{-\frac{\beta_*}{\beta}\cdot \frac{\alpha}{\alpha-\beta_*}}\d x\right)^{\frac{\alpha-\beta_*}{\alpha}(1+\widetilde\mu_1/\alpha)}.
\end{align*}
Thus,
\beq\label{Gfiveest}
    G_5\le J_0^{\frac{\beta_*}{\alpha}(1+\widetilde\mu_1/\alpha)}{\widetilde G}_5.
    \eeq  
Note from the definition of $\beta_*$ in \eqref{bcond} that
\beq\label{bebe}
\beta_*(1-\widetilde\theta) (1+\widetilde\mu_1/\alpha)= \frac{\beta}{p-1}.
\eeq
Then $\Phi_6$ defined in \eqref{Ph67} can be estimated as following, with the aid of \eqref{Gfiveest}, \eqref{bebe} and then inequality \eqref{group2u} for $z=\beta/(p-1)$,
\beq\label{fisixest}
\begin{aligned}
    \Phi_6
    &\le G_1^{\widetilde\theta}\cdot J_0^{\frac{\beta_*}{\alpha}(1+\widetilde\mu_1/\alpha)(1-\widetilde\theta)}{\widetilde G}_5^{1-\widetilde\theta}
    =G_1^{\widetilde\theta}{\widetilde G}_5^{1-\widetilde\theta} J_0^{\frac{\beta}{\alpha(p-1)}}\\
       &\le G_1^{\widetilde\theta}{\widetilde G}_5^{1-\widetilde\theta} \cdot  2^{\frac{\beta}{p-1}(1+\frac1\alpha)} {\mathcal E}_1^{\frac{\beta}{\alpha(p-1)}}
    (1+J_*^{\frac{\beta}{\alpha(p-1)}})
    =  2^{\widehat\beta_2} \widetilde\Phi_6 \left( 1+J_*^{\frac{\beta}{\alpha(p-1)}}\right).
\end{aligned}
\eeq

Next, using \eqref{G2est} and \eqref{Gfiveest} to estimate $\Phi_7$ defined in \eqref{Ph67}, we have
\beq\label{P7o}
\Phi_7
\le (J_0^\frac\beta\alpha {\widetilde G}_2)^{\frac{\widetilde\theta}{1-\widetilde\theta}}  \cdot  (J_0^{\frac{\beta_*}{\alpha}(1+\widetilde\mu_1/\alpha)}{\widetilde G}_5)
= {\widetilde G}_2^{\frac{\widetilde\theta}{1-\widetilde\theta}}{\widetilde G}_5 J_0^{\frac{\beta\widetilde\theta}{\alpha(1-\widetilde\theta)}+\frac{\beta_*}{\alpha}(1+\widetilde\mu_1/\alpha)}.
\eeq
Regarding the last power of $J_0$ in \eqref{P7o}, we have 
\beqs 
\beta_*\left(1+\frac{\widetilde\mu_1}{\alpha}\right)+\frac{\beta\widetilde\theta}{1-\widetilde\theta}
=\frac{\beta}{(1-\widetilde\theta)(p-1)}+\frac{\beta\widetilde\theta}{1-\widetilde\theta}
=\frac{\beta}{(1-\widetilde\theta)}\left(\frac1{p-1}+\widetilde\theta\right)=\widehat\beta_3.
\eeqs 
Applying inequality \eqref{group2u}  to $z=\widehat\beta_3$ to estimate $J_0^{\widehat\beta_3/\alpha}$ in \eqref{P7o}, we obtain
\beq\label{fisevenest}
\begin{aligned}
\Phi_7
&\le  {\widetilde G}_2^{\frac{\widetilde\theta}{1-\widetilde\theta}}{\widetilde G}_5  J_0^{\frac{\widehat\beta_3}{\alpha}}
\le {\widetilde G}_2^{\frac{\widetilde\theta}{1-\widetilde\theta}}{\widetilde G}_5 \cdot  2^{\widehat\beta_3(1+\frac1\alpha)} {\mathcal E}_1^{\frac{\widehat\beta_3}{\alpha}} \left(1+ J_*^{\frac{\widehat\beta_3}{\alpha}}\right)
= 2^{\widehat\beta_3(1+\frac1\alpha)}  \widetilde\Phi_7\left(1+ J_*^{\frac{\widehat\beta_3}{\alpha}}\right).
\end{aligned}
\eeq

Then using estimates \eqref{Ph4est} of $\Phi_4$, \eqref{fisixest} of $\Phi_6$, and \eqref{fisevenest}  of $\Phi_7$  in \eqref{mainTrace}, we have
\beq\label{trace02}
\begin{aligned}
&\int_\Gamma |u|^{\alpha+r} \d S 
\le 3\varepsilon I_* 
+ z_1 \Phi_3 J_*^{1+r/\alpha}
+2^{\widehat\beta_1}  
 \varep^{-\frac\theta{1-\theta}}z_2 \widetilde\Phi_4 \left ( J_*^{1+\mu_1/\alpha}+ J_*^{1+\widehat\mu_1/\alpha}\right)\\
&\quad + 2^{\widehat\beta_2} \varepsilon^{-\frac 1{p-1}}z_4  \widetilde\Phi_6 \left(J_*^{1+\widetilde r/\alpha}+J_*^{1+\widetilde r/\alpha+\frac{\beta}{\alpha(p-1)}}\right)\\
&\quad+ 2^{\widehat\beta_3(1+\frac1\alpha)} \varep^{-(\frac 1 {p-1}+\frac p{p-1}\cdot \frac{\widetilde\theta}{1-\widetilde\theta} )}z_5  \widetilde\Phi_7
\left(J_*^{1+\widetilde\mu_1/\alpha}
+ J_*^{1+\widetilde\mu_1/\alpha+\frac{\widehat\beta_3}{\alpha}}\right).
\end{aligned}
\eeq
Observe that the maximum and minimum of the powers of $J_*$ in \eqref{trace02} are $1+\widehat\mu_4/\alpha$ and $1+\widehat\mu_5/\alpha$, respectively.
Then using the same arguments as in the proof of Lemma \ref{RoS1} when we derived \eqref{RSi2} from \eqref{RSi1}, we obtain  \eqref{tru3} from \eqref{trace02}.
\end{proof}

\begin{proof}[Proof of Lemma \ref{RoS2}]
In this proof, we denote
\begin{align*}
 {\widetilde E}_1&=\left(\int_U \varphi(x)^{-\frac{r_1}{1-r_1}}\d x\right)^\frac{1-r_1}{r_1},&& 
\widehat I_*=\int_0^T\int_U \frac{|u(x,t)|^{\alpha-s}}{(1+|u(x,t)|)^\beta}|\nabla u(x,t)|^p W_*(x)\d x\d t,\\
J&=\int_0^T\int_U |u(x,t)|^{\alpha-s+p} \varphi(x)  \d x\d t,&& 
\widehat J_*=\essup_{t\in(0,T)} \|u(\cdot,t)\|_{L_\varphi^{\alpha}(U)}.
\end{align*} 

In the proof of the above Lemma \ref{WS2}, that is the proof of \cite[Lemma 2.4]{CHK5}, we change $W(x)$ to  
\beq \label{subt}
W(x,t)=(1+|u(x,t)|)^{-\beta}W_*(x).
\eeq 
Then the constant $K_2$ in \cite[formula  (2.77)]{CHK5} becomes $K_2(t)$ and, subsequently, the constant $K_4$ in \cite[formula  (2.80)]{CHK5} becomes $K_4(t)$, and the constant $\Phi_8=K_4$ at the end of the proof  becomes $\Phi_8(t)$.
In fact, $\Phi_8(t)$ has the same formula as $\Phi_8$ in \eqref{Jss} but with $W(x,t)$ in \eqref{subt} replacing $W(x)$, that is,
\beq \label{Ph81}
\Phi_8(t)={\mathcal E}_2\left[ 
 J_1(t)^{1-r_1}
+{\widetilde E}_1^{r_1}\right]^\frac{1}{r_1},\text{ where } 
J_1(t)=\int_U (1+|u(x,t)|)^{\frac{\beta r_1}{1-r_1}} W_*(x)^{-\frac{r_1}{1-r_1}}\d x.
\eeq 
Then the \cite[inequality  (2.83)]{CHK5} still holds with $K_4$ being replaced by
\beqs 
\essup_{t\in(0,T)} K_4(t)=\essup_{t\in(0,T)} \Phi_8(t).
\eeqs  
Thus, we obtain inequality \eqref{ppsi1} where  $W(x)$ is replaced with $W(x,t)$ given in \eqref{subt} and $\Phi_8$ is replaced with
\beq \label{Ph82}
{\widetilde\Phi}_8
=\essup_{t\in(0,T)} \Phi_8(t).    
\eeq 
Explicitly, we have
\beq\label{ppsi5}
\begin{aligned}
\|u\|_{L^{\kappa \alpha}_\varphi(U\times(0,T))}
&\le  (c_7^p m^\frac1{r_1} {\widetilde\Phi}_8)^\frac{1}{\kappa\alpha}  
(\widehat I_* + J)^\frac{1}{\kappa\alpha}
 \cdot  \widehat J_*^{1-\theta_0}
\end{aligned}
\eeq
with, see \eqref{Ph81} and \eqref{Ph82}, 
\beqs
{\widetilde\Phi}_8={\mathcal E}_2\left[ 
 \essup_{t\in(0,T)} \left(J_1(t) ^{1-r_1}\right)
+{\widetilde E}_1^{r_1}\right]^\frac{1}{r_1}.
\eeqs 
Using inequality \eqref{ee2}, we estimate
\beq\label{ph81}
{\widetilde\Phi}_8\le 2^\frac1{r_1}{\mathcal E}_2\left[ 
 \essup_{t\in(0,T)} \left(J_1(t) ^\frac{1-r_1}{r_1}\right)
+{\widetilde E}_1\right].
\eeq 
Thanks to \eqref{albe} and the fact $r_1<1$ from \eqref{rone}, we have $\alpha> \frac{\beta r_1}{1-r_1}$. 
Then, by H\"older's inequality  with powers $\frac{\alpha(1-r_1)}{\beta r_1}$ and $\frac{\alpha(1-r_1)}{\alpha(1-r_1)-\beta r_1}$, we have
\beq
J_1(t)^\frac{1-r_1}{r_1}
= \left(\int_U
\left\{ (1+|u(x,t)|)^{\frac{\beta r_1}{1-r_1}}\varphi^{\frac{\beta r_1}{\alpha(1-r_1)}} \right\}
\cdot 
\left\{ \varphi^{-\frac{\beta r_1}{\alpha(1-r_1)}} W_*^{-\frac{r_1}{1-r_1}}\right\}
\d x \right)^\frac{1-r_1}{r_1}  \le J_0(t)^\frac{\beta}{\alpha}{\widetilde E}_2,\label{ph82}
\eeq
where $J_0(t)=\int_U (1+|u(x,t)|)^{\alpha}\varphi(x)\d x$ and 
\beq\label{E2} 
{\widetilde E}_2=\left(\int_U  \varphi(x)^{-\frac{\beta r_1}{\alpha(1-r_1)-\beta r_1}} W_*(x)^{-\frac{r_1\alpha}{\alpha(1-r_1)-\beta r_1}}\d x \right)^{\frac{\alpha(1-r_1)-\beta r_1}{\alpha r_1}}.
\eeq 
We estimate, by applying estimate \eqref{Joest} to $u=u(\cdot,t)$ and raising to the power $\beta/\alpha$, 
\beq\label{ph83}    
J_0(t)^\frac{\beta}{\alpha}
\le \left\{ 2^\alpha {\mathcal E}_1\left (1+ \int_U |u(x,t)|^{\alpha}\varphi \d x \right)\right\}^\frac{\beta}{\alpha}
 \le 2^{\beta}{\mathcal E}_1^{\frac{\beta}{\alpha}}\left\{1+\left(\int_U |u(x,t)|^{\alpha}\varphi \d x \right)^\frac{\beta}{\alpha}\right\}.
\eeq 
For the last inequality, we used \eqref{ee3} with $p=\beta/\alpha$ which is less than $1$ thanks to  \eqref{albe}.
Thus, combining \eqref{ph81}, \eqref{ph82} and \eqref{ph83}, together with the fact $\mathcal E_1\ge 1$, yields
\begin{align*}
{\widetilde\Phi}_8
&\le  
2^\frac1{r_1} {\mathcal E}_2 \left\{ 2^{\beta}{\mathcal E}_1^{\frac{\beta}{\alpha}} \left[ 
1+  \left(\essup_{t\in(0,T)}\int_U |u(x,t)|^{\alpha}\varphi(x)\d x \right)^\frac{\beta}{\alpha}\right]
\cdot {\widetilde E}_2+{\widetilde E}_1
\right\}\\
&\le 
2^{\beta+\frac1{r_1}}{\mathcal E}_1^{\frac{\beta}{\alpha}} {\mathcal E}_2 ({\widetilde E}_2+{\widetilde E}_1)
 (1+  \widehat J_*^\beta).
\end{align*}

We now have from this estimate of ${\widetilde\Phi}_8$ and \eqref{ppsi5} that
\beq\label{ppsi3}
\|u\|_{L^{\kappa \alpha}_\varphi(U\times(0,T))}
\le  \left\{2^{\beta+\frac1{r_1}} c_7^p m^\frac{1}{r_1}{\mathcal E}_1^{\frac{\beta}{\alpha}} ( {\widetilde E}_2+{\widetilde E}_1)
{\mathcal E}_2 \right\}^\frac{1}{\kappa\alpha}
(1+  \widehat J_*^\beta)^\frac{1}{\kappa\alpha} (\widehat I_*  +J)^\frac{1}{\kappa\alpha}
 \cdot  \widehat J_*^{1-\theta_0}.
\eeq 
Observe that 
\beq\label{E35}
{\widetilde E}_1\le {\mathcal E}_3.
\eeq
With $\alpha>\beta r_1/(1-r_1)$, we estimate ${\widetilde E}_2$ given in \eqref{E2} by 
\beq\label{E22}
{\widetilde E}_2
\le \left(1+\int_U  (1+\varphi ^{-1})^{\frac{\beta r_1}{\alpha(1-r_1)-\beta r_1}} (1+W_* ^{-1})^{\frac{r_1\alpha}{\alpha(1-r_1)-\beta r_1}}\d x \right)^{\frac{\alpha(1-r_1)-\beta r_1}{\alpha r_1}}.    
\eeq 
Consider the powers under the integral in \eqref{E22} as  functions in $\alpha\in[\beta/(1-r_1),\infty)$, that is,
$$\alpha\mapsto \frac{\beta r_1}{\alpha(1-r_1)-\beta r_1}\text{ and } 
\alpha\mapsto \frac{r_1 \alpha}{\alpha(1-r_1)-\beta r_1}\text{ for }\alpha\ge \beta/(1-r_1).$$
Then they are decreasing functions. Hence, they attain their maximum values when $\alpha=\beta/(1-r_1)$. For an upper bound of the last power in \eqref{E22}, we just drop the $(-\beta r_1)$ part, i.e.,
\beqs 
\frac{\alpha(1-r_1)-\beta r_1}{\alpha r_1}\le \frac{1-r_1}{r_1}.
\eeqs 
Therefore, it follows from \eqref{E22} that
\beq\label{E23}
{\widetilde E}_2
\le \left(1+\int_U (1+ \varphi ^{-1})^{\frac{\beta r_1}{\beta-\beta r_1}} (1+W_* ^{-1})^{\frac{r_1\beta/(1-r_1)}{\beta-\beta r_1}}\d x 
\right)^{\frac{1-r_1}{r_1}}
={\mathcal E}_3.    
\eeq 
One has $\alpha>1$ thanks to \eqref{ag1} and the first two inequalities in the assumption \eqref{alcond}. 
Also, $\kappa>1$  thanks to \eqref{defkappa}.
Then, in estimate \eqref{ppsi3}, applying inequality  \eqref{ee3} with $p=1/(\kappa\alpha)<1$ to estimate $( 1+ \widehat J_*^\beta )^\frac{1}{\kappa\alpha}$ gives
\beq \label{JJJ}
 ( 1+ \widehat J_*^\beta )^\frac{1}{\kappa\alpha}\widehat J_*^{1-\theta_0}
\le (1+ \widehat J_*^\frac{\beta}{\kappa\alpha})\widehat J_*^{1-\theta_0}=\widehat J_*^{1-\theta_0}+\widehat J_*^{1-\widehat\theta_0}.
\eeq 
Thus, combining \eqref{ppsi3} with \eqref{E35}, \eqref{E23} and \eqref{JJJ}, we obtain 
\beqs
\|u\|_{L^{\kappa \alpha}_\varphi(U\times(0,T))}
\le  \left\{2^{\beta+\frac1{r_1}} c_7^p m^\frac{1}{r_1}{\mathcal E}_1^{\frac{\beta}{\alpha}} \cdot ( 2\mathcal E_3)\cdot 
{\mathcal E}_2 \right\}^\frac{1}{\kappa\alpha}
(\widehat I_*  +J)^\frac{1}{\kappa\alpha}
(\widehat J_*^{1-\theta_0}+\widehat J_*^{1-\widehat\theta_0}),
\eeqs 
which proves \eqref{pssi6}.
\end{proof}

\medskip
\noindent\textbf{Data availability.} 
No new data were created or analyzed in this study.

\medskip
\noindent\textbf{Methods.} 
No Artificial Intelligence Generated Content (AIGC) tools are used in developing any portion of this paper. 

\medskip
\noindent\textbf{Funding.} No funds were received for conducting this study. 

\medskip
\noindent\textbf{Conflict of interest.}
There are no conflicts of interests.

\bibliography{paperbaseall}{}
 \bibliographystyle{abbrv}

\end{document}